\documentclass{gtart_h}  

\def\ifplaintex{\expandafter\ifx\csname documentclass\endcsname\relax}

\ifplaintex 
\else
\fi

\expandafter\ifx\csname beginpicture\endcsname\relax
\expandafter\ifx\csname documentclass\endcsname\relax
\input pictex \else
\input prepictex \input pictex \input postpictex \fi\fi

\def\gt{{\mathsurround=0pt\it $\cal G\mskip-2mu$eometry \&\ 
$\cal T\!\!$opology}}        

\def\gtp{{\mathsurround=0pt\it $\cal G\mskip-2mu$eometry \&\ 
$\cal T\!\!$opology $\cal P\!$ublications}}  

\def\lognumber#1{\def\thelognumber{#1}}
\def\volumenumber#1{\def\thevolumenumber{#1}}
\def\papernumber#1{\def\thepapernumber{#1}}
\def\volumeyear#1{\def\thevolumeyear{#1}}

\def\pagenumbers#1#2{\def\startpage{#1}\def\finishpage{#2}}
\def\published#1{\def\publishdate{#1}}
\def\proposed#1{\def\theproposer{#1}}
\def\seconded#1{\def\theseconders{#1}}
\def\received#1{\def\receiveddate{#1}}
\def\revised#1{\def\reviseddate{#1}}
\def\accepted#1{\def\accepteddate{#1}}

\long\def\asciiabstract#1{\long\def\theasciiabstract{#1}}

\let\\\par\let\thelognumber\relax
\let\thevolumenumber\relax\let\thepapernumber\relax
\let\thevolumeyear\relax\let\thesamplenumber\relax\let\startpage\relax
\let\finishpage\relax\let\publishdate\relax\let\receiveddate\relax
\let\reviseddate\relax\let\accepteddate\relax\let\theasciititle\relax
\let\theasciiauthors\relax
\let\theasciiabstract\relax
\let\theasciiemail\relax\let\theshortauthors\relax\let\theshorttitle\relax

\long\def\maketitlep{   

\count0=\startpage

\gt\hfill      
\beginpicture
\setcoordinatesystem units <0.33truein, 0.33truein> point at 2.2 0.9
\setplotsymbol ({$\cal G$})
\plotsymbolspacing=9truept
\circulararc 315 degrees from 0 1 center at 0 0
\setplotsymbol ({$\cal T$})
\circulararc 315 degrees from 1 -1 center at 1 0
\endpicture
\break
{\small\ifx\thesamplenumber\relax 
Volume \else Sample
\fi\thevolumenumber\ (\thevolumeyear)
\startpage--\finishpage\nl
Published: \publishdate}
\vglue 0.5truein plus 0.4fil minus 0.1truein

{\parskip=0pt\leftskip 0pt plus 1fil\def\\{\par\smallskip}{\ifplaintex\large
\else\Large\fi\bf\thetitle}\par\medskip}   

\vglue 0pt plus 0.1fil 

{\parskip=0pt\leftskip 0pt plus 1fil\def\\{\par}{\sc\theauthors}
\par\medskip}

\vglue 0pt plus 0.1fil 

{\small\parskip=0pt\let\newline\\
{\leftskip 0pt plus 1fil\def\\{\par}{\sl\theaddress}\par}
\expandafter\ifx\theemail\relax    
\relax\else\vglue 5pt plus 0.02fil minus 2pt\def\\{\stdspace{\rm 
and}\stdspace} 
\cl{Email:\stdspace\tt\theemail}\fi
\ifx\theurl\relax                  
\relax\else\vglue 5pt plus 0.02fil minus 2pt\def\\{\stdspace{\rm 
and}\stdspace}
\cl{URL:\stdspace\tt\theurl}\fi\par}

\vglue 7pt plus 0.3fil minus 3pt

{\bf Abstract}
\vglue 5pt plus 0.1fil minus 2pt

\theabstract

\vglue 7pt plus 0.3fil minus 3pt

{\bf AMS Classification numbers}\quad Primary:\quad \theprimaryclass

Secondary:\quad \thesecondaryclass

\vglue 5pt plus 0.3fil minus 2pt

{\bf Keywords:}\quad \thekeywords

\vglue 10pt plus 0.5fil minus 5pt

{\small  Proposed: \theproposer\hfill Received: \receiveddate\nl
Seconded: \theseconders\hfill 
\ifx\reviseddate\relax                         
Accepted: \accepteddate                        
\else
Revised: \reviseddate                          
\fi}
\eject
}       

\let\maketitlepage\maketitlep
\let\maketitle\maketitlepage

\font\phead=cmsl9 scaled 950
\font=cmsl9 scaled 1050
\font\pnum=cmbx10 scaled 913
\font\lnum=cmbx10 
\font\pfoot=cmsl9 scaled 950
\font=cmsl9 scaled 1050
\ifplaintex
\headline{\vbox to 0pt{\vskip -4.5mm\line{\small\phead\ifnum
\count0=\startpage ISSN 1364-0380 (on line)
1465-3060 (printed) \hfill {\pnum\folio}\else\ifodd\count0\def\\{ }%
\ifx\theshorttitle\relax\thetitle\else\theshorttitle\fi\hfill{\pnum\folio}
\else\def\\{ and }{\pnum\folio}\hfill\ifx\theshortauthors\relax\theauthors
\else\theshortauthors\fi\fi\fi}\vss}}
\footline{\vbox to 0pt{\vglue 0mm\line{\small\pfoot\ifnum\count0=\startpage
\copyright\ \gtp\hfill\else
\gt, Volume \thevolumenumber\ (\thevolumeyear)\hfill\fi}\vss
}}
\else
\makeatletter
\def\@oddhead{{\small\lhead\ifnum\count0=\startpage ISSN 1364-0380 (on line)
1465-3060 (printed) \hfill {\lnum\number\count0}\else\ifodd\count0
\def\\{ }\ifx\theshorttitle\relax \thetitle \else\theshorttitle\fi\hfill
{\lnum\number\count0}\else\def\\{ and }{\lnum\number\count0}
\hfill\ifx\theshortauthors\relax 
\theauthors\else\theshortauthors\fi\fi\fi}}\def\@evenhead{\@oddhead}
\def\@oddfoot{\small\lfoot\ifnum\count0=\startpage\copyright\ \gtp\hfill\else
\gt, Volume \thevolumenumber\ (\thevolumeyear)\hfill\fi}
\def\@evenfoot{\@oddfoot}
\makeatother
\fi

\newwrite\gtoutfile
\long\gdef\makeheadfile{  
{\def\\{, }\def\s{ }
\immediate\openout\gtoutfile head.xxx
\immediate\write\gtoutfile{Proxy-for: \ifx\theasciiauthors\relax
\theauthors\else\theasciiauthors\fi\s<\ifx\theasciiemail\relax\theemail\else\theasciiemail\fi>}
\immediate\write\gtoutfile{\noexpand\\}
\immediate\write\gtoutfile{Authors: \ifx\theasciiauthors\relax
\theauthors\else\theasciiauthors\fi}
{\def\\{ }\immediate\write\gtoutfile{Title: \ifx\theasciititle\relax
\thetitle\else\theasciititle\fi}}
\immediate\write\gtoutfile{Subj-class: GT or SG or MG etc}
\immediate\write\gtoutfile{MSC-class: \theprimaryclass\ifx\thesecondaryclass\relax\else, \thesecondaryclass\fi}
\immediate\write\gtoutfile{Journal-ref: Geom. Topol. \thevolumenumber
(\thevolumeyear) \startpage-\finishpage}
\immediate\write\gtoutfile{Comments: Published by Geometry and Topology at}
\immediate\write\gtoutfile{\s\s http://www.maths.warwick.ac.uk/gt/GTVol\thevolumenumber/paper\thepapernumber.abs.html}
\immediate\write\gtoutfile{\noexpand\\}
\immediate\write\gtoutfile{}
\ifx\theasciiabstract\relax
\immediate\write\gtoutfile{\theabstract}\else
\immediate\write\gtoutfile{\theasciiabstract}\fi
\immediate\write\gtoutfile{}
\immediate\write\gtoutfile{\noexpand\\}
\immediate\write\gtoutfile{}
\immediate\closeout\gtoutfile}}  

\def\maketitlepage{\maketitlep\makeheadfile}
\let\maketitle\maketitlepage

\lognumber{481}
\received{29 July 2004}
\volumenumber{9}\papernumber{13}\volumeyear{2005}
\pagenumbers{403}{482}   
\revised{17 February 2005}
\published{9 March 2005}
\accepted{22 February 2005}
\proposed{Walter Neumann}
\seconded{Martin Bridson, David Gabai}

\usepackage{amssymb, amsmath}
 
\let\sf\sl
\let\tilde\widetilde

\newfont{\cyr}{wncyr10}
\newcommand{\hh}{\ensuremath{\mathbb H}}  
  
\newcommand{\rr}{\ensuremath{\mathbb R}}
\newcommand{\rb}{\ensuremath{\bar{\mathbb R}}}

\renewcommand{\ss}{\ensuremath{\mathbb S}}
\newcommand{\zz}{\ensuremath{\mathbb Z}}
\newcommand{\cc}{\ensuremath{\mathbb{C}}}

\newcommand{\ff}{\ensuremath{\mathcal{F}}}

\newcommand{\al}{\ensuremath{\alpha}}

\newcommand{\be}{\ensuremath{\beta}}
\newcommand{\g}{\ensuremath{\gamma}}
\newcommand{\G}{\ensuremath{\mathcal{G}}}
\newcommand{\de}{\ensuremath{\delta}}

\newcommand{\s}{\ensuremath{\sigma}}  
\newcommand{\la}{\ensuremath{\lambda}}
\newcommand{\T}{\ensuremath{\mathcal{T}}}

\newcommand{\p}{\ensuremath{\partial}}  
\newcommand{\jx}{\ensuremath{{X\!\Join\! X}}}
\newcommand{\becross}{\ensuremath{\be^{\scriptscriptstyle\times\!}}}
\newcommand{\dcross}{\ensuremath{d^{\scriptscriptstyle\times\!}}}

\newcommand{\simplus}{\ensuremath{
\mspace{3mu}\hbox{$\sim$}\mspace{-14.5mu}\raise-0.5ex
\hbox{$\scriptscriptstyle +$}\mspace{8mu}
}}

\newcommand{\ssimplus}{\ensuremath{
\mspace{3mu}\hbox{$\sim$}\mspace{-11.5mu}\raise-0.5ex
\hbox{$\scriptscriptstyle +$}\mspace{8mu}
}}

\newcommand{\sj}{\ensuremath{
{
\raise0.1ex\hbox{$\circ\mspace{-9mu}*$}
}}}

\newcommand{\starx}{\ensuremath{
{\raise0.2ex\hbox{$*$}
X}}}

\newcommand{\starxb}{\ensuremath{
{\raise0.2ex\hbox{$*$}
\bar{X}}}}

\newcommand{\starpx}{\ensuremath{
{\raise0.2ex\hbox{$*$}
(\partial X)}}}

\newcommand{\starxxo}{\ensuremath{
\raise-0.3ex\hbox{$\scriptscriptstyle x_0$}\mspace{-12mu}
   {\raise0.22ex\hbox{$*$}
X}}}

\newcommand{\starxbxo}{\ensuremath{
\raise-0.3ex\hbox{$\scriptscriptstyle x_0$}\mspace{-12mu}
   {\raise0.22ex\hbox{$*$}
\bar{X}}}}

\newcommand{\sy}{\ensuremath{
{\raise0.2ex\hbox{$\circ\mspace{-9mu}*$}
Y}}}

\newcommand{\sx}{\ensuremath{
{\raise0.2ex\hbox{$\circ\mspace{-9mu}*$}
\mspace{-1.5mu}
X}}}

\newcommand{\sxb}{\ensuremath{
{\raise0.2ex\hbox{$\circ\mspace{-9mu}*$}
\mspace{-1.5mu}\bar{X}}}}

\newcommand{\sxxo}{\ensuremath{
\raise-0.3ex\hbox{$\scriptscriptstyle x_0$}\mspace{-12mu}
   {\raise0.22ex\hbox{$\circ$}\mspace{-11.5mu}\raise0.22ex\hbox{$*$}
X}}}

\newcommand{\sxbxo}{\ensuremath{
\raise-0.3ex\hbox{$\scriptscriptstyle x_0$}\mspace{-12mu}
   {\raise0.22ex\hbox{$\circ$}\mspace{-11.5mu}\raise0.22ex\hbox{$*$}
\bar{X}}}}

\newcommand{\hsxb}{\ensuremath{
{\raise0.12ex\hbox{$\scriptscriptstyle\smallsmile$}
   \mspace{-12.2mu}\raise0.3ex\hbox{$*$}   
\bar{X}}}}

\newcommand{\spx}{\ensuremath{
{\raise0.2ex\hbox{$\circ\mspace{-9mu}*$}
\p X}}}

\newcommand{\sxxone}{\ensuremath{
\raise-0.3ex\hbox{$\scriptscriptstyle x_1$}\mspace{-10mu}
   {\raise0.22ex\hbox{$\circ$}\mspace{-9mu}\raise0.22ex\hbox{$*$}
X}}}

\newcommand{\diamondx}{\ensuremath{
\raise0.2ex\hbox{$\diamond$}
X}}

\newcommand{\diamondxb}{\ensuremath{
\raise0.2ex\hbox{$\diamond$}
\bar{X}}}

\newcommand{\diamondxbxo}{\ensuremath{
\raise-0.3ex\hbox{$\scriptscriptstyle x_0$}
   {\mspace{-10mu}\raise0.3ex\hbox{$\diamond$}
\bar{X}}}}

\newcommand{\diamondxxo}{\ensuremath{
\raise-0.3ex\hbox{$\scriptscriptstyle x_0$}
   {\mspace{-12mu}\raise0.3ex\hbox{$\diamond$}
X}}}

\newcommand{\No}{N\raise4pt\hbox{\tiny o}\kern+.2em}
\newcommand{\Isom}{\ensuremath{\mathrm{Isom}}}
\newcommand{\Geod}{\ensuremath{\mathrm{Geod}}}
\newcommand{\fx}{\ensuremath{\mathcal{F}(X)}}

\newcommand{\inv}{\ensuremath{^{-1}}} 
\newcommand{\dhat}{\ensuremath{\hat{d}}}

\newcommand{\se}{\ensuremath{\subseteq}}

\newcommand{\ds}{\ensuremath{\displaystyle}}

\newcommand{\lbr}{\ensuremath{ | \! [ }}
\newcommand{\rbr}{\ensuremath{ ] \! | }}

\newcommand{\fs}{\mathcal{F}}

\newtheorem{ttt}{Theorem}  
\newtheorem{lemma}[ttt]{Lemma}

\newtheorem{ppp}[ttt]{Proposition}  
 
\theoremstyle{definition}
\newtheorem{ddd}[ttt]{Definition}  

\begin{document}  
\title{Flows and joins of metric spaces}
\author{Igor Mineyev}

\address{Department of Mathematics, University of Illinois at Urbana-Champaign\\
250 Altgeld Hall, 1409 W Green Street, Urbana, IL 61801, USA}

\email{mineyev@math.uiuc.edu}

\begin{abstract}
We introduce the functor~$\sj$ which assigns to every 
metric space~$X$ its {\sf symmetric join}~$\sx$.
As a set, $\sx$ is a union of intervals connecting ordered
pairs of points in~$X$.
Topologically, $\sx$ is a natural quotient of the usual join of $X$
with itself. We define an $\Isom(X)$--invariant metric 
$d_*$ on~$\sx$. 

Classical concepts known for $\hh^n$ and negatively curved
manifolds are defined in a precise way for any {\sf hyperbolic complex}~$X$,
for example for a Cayley graph of a Gromov hyperbolic group.
We define a {\sf double difference}, a {\sf cross-ratio} 
and {\sf horofunctions}
in the compactification $\bar{X}=X\sqcup\p X$. They are continuous, $\Isom(X)$--invariant,
and satisfy sharp identities. We characterize 
the {\sf translation length} of a hyperbolic isometry $g\in\Isom (X)$.

For any hyperbolic complex~$X$, the symmetric join
$\sxb$ of $\bar{X}$ and the (generalized) metric $d_*$
on it are defined. The {\sf geodesic flow
space} $\ff(X)$ arises as a part of~$\sxb$.
$(\ff(X),d_*)$ is an analogue of (the total space of) the unit tangent bundle
on a simply connected negatively curved manifold. 
This flow space is defined for any hyperbolic complex~$X$ and has sharp properties.
We also give a construction of the {\sf asymmetric join} $X\sj\, Y$ of two metric spaces.

These concepts are canonical, i.e.\ functorial in~$X$, and involve no ``quasi"-language.
Applications and relation to the Borel conjecture and others are discussed.
\end{abstract}

\asciiabstract{%
We introduce the functor * which assigns to every metric space X its
symmetric join *X.  As a set, *X is a union of intervals connecting
ordered pairs of points in X.  Topologically, *X is a natural quotient
of the usual join of X with itself. We define an Isom(X)-invariant
metric d* on *X.  Classical concepts known for H^n and negatively
curved manifolds are defined in a precise way for any hyperbolic
complex X, for example for a Cayley graph of a Gromov hyperbolic
group.  We define a double difference, a cross-ratio and horofunctions
in the compactification X-bar= X union bdry X.  They are continuous,
Isom(X)-invariant, and satisfy sharp identities.  We characterize the
translation length of a hyperbolic isometry g in Isom(X).  For any
hyperbolic complex X, the symmetric join *X-bar of X-bar and the
(generalized) metric d* on it are defined.  The geodesic flow space
F(X) arises as a part of *X-bar.  (F(X),d*) is an analogue of (the
total space of) the unit tangent bundle on a simply connected
negatively curved manifold.  This flow space is defined for any
hyperbolic complex X and has sharp properties.  We also give a
construction of the asymmetric join X*Y of two metric spaces.  These
concepts are canonical, i.e. functorial in X, and involve no
`quasi'-language.  Applications and relation to the Borel conjecture
and others are discussed.}



{\small\maketitlepage}

\setcounter{section}{-1}

\section{Introduction}
Let $X$ be a proper geodesic hyperbolic metric space in the sense of Gromov.
In~\cite{MMS}, the following discrete--continuous dichotomy was shown for
a non-elementary closed
subgroup $H<\mathrm{Isom}(X)$ acting cocompactly on~$X$:\break either
\begin{itemize}
\item  $H$ has a proper non-elementary vertex-transitive action on a hyperbolic graph
of bounded valency, or
\item there is a finite-index open subgroup $H^*< H$ and a compact normal subgroup
$K\triangleleft H$ contained in~$H^*$ such that $H^*/K$ is a connected simple Lie group
of rank one.
\end{itemize}
This says, less formally, that to understand general hyperbolic spaces 
it suffices to study hyperbolic graphs and Lie groups. While the theory of Lie groups
and symmetric spaces is quite developed, the hyperbolic graphs and groups
introduced by Gromov~\cite{Gr2} are relatively recent phenomena.
By its very nature of being discrete a hyperbolic graph lacks a nice local structure, and therefore
the tools of differential geometry. 

In this paper we fill in the gaps in the discrete spaces and 
the blanks in the discrete spaces theory.
The main philosophical point is that hyperbolic groups, despite being discrete,
do give rise to many concepts that were known on the~continuous side.
This paper introduces several {\em sharp} geometric concepts first for arbitrary
metric spaces, and then for hyperbolic spaces in the sense of Gromov.
These concepts in particular eliminate the need of the ``quasi"-language
when talking about hyperbolic groups and spaces.

We introduce the notion of {\sf symmetric join}.
If $X$ is a set, the symmetric join of~$X$, denoted $\sx$, is the ``obvious'' union of
formal intervals connecting ordered pairs of points in~$X$;
the interval connecting a point to itself is required to degenerate.
When $X$ is a topological space we define a natural topology on~$\sx$.
When~$X$ is a metric space, we define a metric~$d_*$ on $\sx$ with natural 
properties. The symmetric join is therefore an example of a {\sf metric join}.
Even though $\sx$ is an abstract union of intervals, the construction
of~$d_*$ is very explicit. The metric~$d_*$ is canonical and $\Isom(X)$--invariant.

In~\cite{MY}, the Baum--Connes conjecture was proved for hyperbolic groups
and their subgroups by constructing a strongly bolic metric~$\dhat$
on any hyperbolic group. Having a strongly bolic metric is not sufficient
for the constructions of the present paper.
We show that (a modified version of) $\dhat$ has stronger properties
and use it to define the {\sf double difference} $\left<\cdot,\cdot|\cdot,\cdot\right>$
in~$X$ (see~\ref{ss_dhat-everywhere}). The main property used is that $\dhat$
and $\left<\cdot,\cdot|\cdot,\cdot\right>$ behave ``exponentially well" at infinity
(Theorem~\ref{t_dhat}). We show that the double difference continuously extends
to $\bar{X}$ and gives rise to a continuous {\sf cross-ratio} in~$\bar{X}$
(section~\ref{s_cross-ratio}). It is the use of the metric $\dhat$ that allows
things to extend continuously to the boundary.
This generalizes the work of Otal~\cite{Ot4} who defined and used the cross-ratio
for negatively curved manifolds.

The ``Hyperbolic groups'' article by Gromov~\cite{Gr2} was an inspiration for many
mathematicians over the last years, including the author of this paper.
Gromov outlined a construction of a metric space~$\hat{G}$ with $\rr$--, $\zz_2$--,
and $\Gamma$--actions~\cite[8.3.C]{Gr2}. He considers the set of all biinfinite
geodesics in the Cayley graph and then identifies those geodesics that connect
the same pairs of points in $\p\Gamma$. 
Math\'eus~\cite{Matheus1991} and Champetier~\cite{Champetier1994}
provided further details of the Gromov's construction. The identification of geodesics
is by quasiisometries, so $\hat{G}$ is rather a quasigeodesic flow; $\rr$ acts
on the $\rr$--orbits in~$\hat{G}$ by quasiisometric homeomorphisms.

In~\cite{Furman2002} Furman takes $\p^2\Gamma\times\rr$ as a model set for
the flow space, so geodesics are unique by definition. He uses boundedness of cocycles
from \cite{Coornaert1993} to construct a geodesic current, i.e.\ an invariant measure
on $\p^2 \Gamma$, then provides a $\Gamma$--action
on~$\p^2\Gamma\times\rr$ and a $\Gamma$--invariant cross-ratio on $\p\Gamma$.
Both the action and the cross-ratio are measure-theoretic, that is defined
up to subsets of measure~0. B\"uhler~\cite{Buehler2002}
considers the space of all geodesics in a hyperbolic space~$X$, as in
the Gromov's construction, and uses the amenability of the $\Gamma$--action
on $\p\Gamma$ to make a measurable choice of geodesic, i.e.\ a choice of geodesic
for almost every pair of points in $\p\Gamma$. Since the geodesics in~$X$ are chosen
in a measurable fashion, and they usually do not depend continuously on their endpoints,
there is no obvious way to topologically identify the union of such geodesics
with $\p^2\Gamma\times\rr$. In both \cite{Furman2002} and~\cite{Buehler2002}
the space considered is a measure space rather than a metric space.
Bourdon~\cite{Bou} presented geodesic flows with sharp properties
in the case of $\mathrm{CAT}(-1)$ spaces.

The present paper provides a new approach to constructing a geodesic flow~$\fs(X)$
for an arbitrary hyperbolic complex~$X$, for example when $X$ is a Cayley
graph of a hyperbolic group (see~\ref{ss_hyp-compl} for definitions).
First we enlarge $\sx$ to the symmetric join $\sxb$
of the compactification~$\bar{X}=X\sqcup\p X$.
We put the metric~$\dhat$ on~$X$ and show that
the metric~$d_*$ canonically determined by $\dhat$ extends to $\sxb$
(with the obvious infinite values allowed at infinite points). The use of $\dhat$
is essential here. Then $\ff(X)\se\sxb$ is by definition the union of lines in $\sxb$
that connect pairs of points in~$\p X$, equipped with the restricted metric~$d_*$.

Our construction of symmetric join $\sxb$ is more general than the flow space,
since it allows for lines to connect points in~$X$ as well as points in~$\p X$.
But even when restricted to $\ff(X)$ it provides strong properties
(see Theorem~\ref{geodesic-flow}).
Generally speaking, the outcomes of this construction are {\em continuous}, rather
than measurable as in~\cite{Furman2002, Buehler2002},
and {\em sharp}, rather than defined ``up to a bounded amount" as
in~\cite{Gr2, Matheus1991, Champetier1994, Paulin1996}.
Continuity is important for future topological applications.
The construction is also more general since we consider an arbitrary hyperbolic complex
$X$ with the action of the full simplicial isometry group $\Isom(X)$, i.e.\ the group
acting on $X$ does not have to be discrete.
The metric $d_*$, both on~$\ff(X)$ and on~$\sx$, is $\mathrm{Isom}(X)$--invariant.
Each $\rr$--orbit in $\ff(X)$ is an isometric copy of~$\rr$.
$\rr$ acts on $(\ff(X),d_*)$ by bi-Lipschitz homeomorphisms, and on each $\rr$--orbit
by isometries in the standard way.
The $\rr$--orbits in $\fs(X)$ converge synchronously and uniformly exponentially.

Symmetric join is a unified concept relating topology and geometry:
it combines the usual notion of topological join with the notion
of geodesic flow.
Symmetric join plays the role of a Riemannian structure on
a manifold, and it is canonically assigned to every metric space.
When $X$ is a hyperbolic complex,  $\sxb$ provides a link between
the local and global structures of~$X$. This is important, for example, in the study
of the topology of manifolds; the manifolds can be chosen to be smooth or not.

For any hyperbolic group~$\Gamma$ acting on a hyperbolic complex~$X$,
for example its Cayley graph, $\fs(X)$ provides a convenient model space.
It plays the role of (the total space of) the unit tangent bundle over a manifold,
even though {\em no manifold was given}.
There are other models provided by the construction:
$\sx$, $\sxb$, $\ff(X)$, $X\sj X$, $X\sj \bar{X}$, $\bar{X}\sj \bar{X}$, 
etc (see section~\ref{s_asym-join}).
They are all equipped with canonical induced $\Gamma$--actions. Their geometry is
closely related to the geometry of~$\Gamma$ but, unlike~$\Gamma$, the spaces are
not discrete. It is an interesting question how one can use these models to generalize
the Farrell--Jones theory \cite{FJ1986, FJ1987, FJ1989, FJ1989-PNASUSA, FJ1990}
to manifolds with hyperbolic fundamental groups. Also, the asymmetric join of two manifolds
might be used in place of a cobordism. One would probably need
to generalize the theory of Chapman~\cite{Chapman1979}, Ferry~\cite{Ferry1997}
and Quinn~\cite{Quinn1979}, and to come up with an ``h--join theorem''
that would play the role of the h--cobordism theorem.

We define continuous {\sf horofunctions} both in~$\bar{X}$ and in~$\sxb$
(Theorem~\ref{t_horofunction}).
They depend only on a point at infinity, rather than on a ray converging
to the point (section~\ref{s_horofunctions-in-barX}).
For each hyperbolic $g\in\Isom(X)$, we define the translation length $\hat{l}(g)$
in terms of metric~$\dhat$ on~$X$. $\hat{l}(g)$ is indeed realized as the shift of
the axis of $g$ in $\ff(X)$ (see section~\ref{s_translation} and
Theorem~\ref{geodesic-flow}(i)).

Symmetric join is used to provide a notion of join of {\em two} metric
spaces $Y$ and $Y'$, called the {\sf asymmetric join} $Y \sj\, Y'$
(section \ref{s_asym-join}). The asymmetric join is therefore another
example of a metric join\label{i_metric-join2}.  In the case when the
two spaces are given an action by the same group, for example when $Y$
and $Y'$ are the universal covers of two manifolds with the same
fundamental group, we describe a canonical way to put a metric on $Y
\sj\, Y'$ (see~\ref{ss_asym-join}). This situation occurs, for
example, in the Borel conjecture that asserts that two closed
aspherical manifolds with the same fundamental group are homeomorphic.

If one thinks of the intervals in $Y \sj\, Y'$ as parameterized by $[-\infty, \infty]$, 
the set of slices at various times $t\in [-\infty, \infty]$ is a {\sf sweep-out},
or rather sweep-between, from $Y$ to $Y'$. Topologically slices are the same for all
$t\in \rr$ (they are all homeomorphic to $Y\times Y'$),
but the metric on $Y \sj\, Y'$ makes it interesting:
the slices indeed approach $Y$ as $t\to-\infty$ and $Y'$ as $t\to\infty$,
in a metrically controlled way (Gromov--Hausdorff convergence on bounded subsets).
This is of interest in particular in relation with the Borel conjecture.
If there exists a homeomorphism between manifolds
$M$ and $M'$, then it must be present in each slice of $M\sj M'$, diagonally.
Since the construction is equivariant, it also
allows working with the universal coverings (this gives even more freedom):
denote $Y:={\tilde M}$ and $Y':={\tilde M'}$.
One needs to find an equivariant copy of $\tilde M$ in each slice
of $Y\sj Y'$ (see section \ref{s_asym-join}). 
Another advantage of using universal coverings is that when $\pi_1(M)=\pi_1(M')$
is hyperbolic, the ideal boundary of $Y$ and $Y'$ can also be used, since
the join of the compactifications $\bar{Y}\sj\bar{Y'}$ is defined 
as well.

The constructions of this paper require sharp estimates carefully written out. 
At the first reading, the reader might want to take a look at~\ref{ss_sj-as-top} and
theorems~\ref{t_dhat}, \ref{dd}, \ref{cr}, \ref{t_extended-d*}, \ref{t_horofunction} 
and~\ref{geodesic-flow} which constitute the main results of this paper.
Sections \ref{ss_dd-gp}--\ref{ss_metric-complexes} deal with 
arbitrary metric spaces and simplicial complexes. After that we work in the category
of (uniformly locally finite) hyperbolic complexes.

One interesting open problem is to prove the group-theoretic rigidity
conjecture, that is the Borel conjecture in the case of hyperbolic
fundamental groups. This implies the Poincar\'e
conjecture~\cite{FJ1990a,FJ1989}, and it can be viewed as a
group-theoretic (or discrete) analog of the Mostow rigidity theorem.

As we mentioned above,  the symmetric join plays the role of a Riemannian structure.
Note also that there are examples of closed manifolds with hyperbolic
fundamental groups that do not admit a Riemannian structure of negative curvature
\cite{DJ1991, CD1995}, and our symmetric join construction applies in those cases.
Moreover, it applies to all dimensions, and to all PL manifolds that do not admit
smooth structure.

Another interesting question is the Cannon's
conjecture~\cite{Ca, CFKW1997, CS1998, CC1992}
that a hyperbolic group~$\Gamma$ with $\p\Gamma$ homeomorphic to~$S^2$ admits
a proper cocompact action on~$\hh^3$. Note that the boundary of a hyperbolic group
is usually very much {\em not} a smooth manifold, even if it is a manifold
topologically, so one would probably need to use metric geometry rather than
the Riemannian one. 
The essence of all these questions is establishing a link between the local and global
structures of hyperbolic metric spaces, and the symmetric flow provides such a link.

The paper is organized as follows.
Preliminary definitions are given in section~\ref{s_preliminaries}.
In section~\ref{s_symm-join} the symmetric join of a 
metric space $X$ is described. In section~\ref{s_metric-on-sj} we define 
a metric $d_*$ on the symmetric join.
In section~\ref{s_metric-and-topology} we discuss
the metric and the topology of~$X$.
In section~\ref{s_metric-complexes} we define metric
complexes and hyperbolic complexes.
Section~\ref{s_metric-and-dd} describes a metric
$\dhat$ and the double difference on a hyperbolic complex~$X$.
In section~\ref{s_cross-ratio} cross-ratio on $X$ is defined.
In section~\ref{s_sj-of-barX} the symmetric join is extended
to~$\bar{X}$. Section~\ref{s_top-of-sxb} deals with the topology
of the extended symmetric join $\sxb$. In section~\ref{s_horofunctions-in-barX}
we define horofunctions and horospheres in~$\sxb$.
Sections~\ref{s_synchronous} and~\ref{s_translation}
prove some technical results about convergence of lines
and translation length. In section~\ref{s_the-geod-flow}
the geodesic flow space is defined and its properties are summarized.
Section~\ref{s_asym-join} defines asymmetric join and gives general
remarks.  
The paper ends with an index of symbols and terminology \ref{index}.

The author benefited a lot from helpful conversations with Alex Furman,
Misha Gromov, Tadeusz Januszkiewicz, Lowell Jones, Misha Kapovich,
Yair Minsky, Frank Quinn and Shmuel Weinberger. Misha Gromov comments
that the existence of a continuous cross-ratio and a geodesic flow with sharp
properties can be also deduced from the discussion in~\cite{Gr2}.

The author would like to thank the hospitality of
MSRI, Berkeley, in the summer of 2002, of Max-Planck-Institut, Bonn in the summer of 2003,
and of IAS, Princeton, in the year 2003--04, where he was supported by 
NSF grant DMS-0111298.
This project is partially supported by NSF CAREER grant DMS-0228910.

\section{Preliminaries}\label{s_preliminaries}
\subsection{The double difference and Gromov product}
\label{ss_dd-gp}
Let $(X, d)$ be a metric space.
The {\sf double difference in} $X$ is
the function\newline
$\left<\cdot,\cdot|\cdot,\cdot\right>\co X^4\to\rr$ defined by
\begin{equation}\label{dd-def}
\left<a,a'|b,b'\right>:=
  \frac{1}{2} \Big( d(a,b)- d(a',b)- d(a,b')+ d(a',b') \Big).
\end{equation}
This notion appeared in a paper by Paulin~\cite{Paulin1996}
(under the name ``cross ratio").

The following properties are immediate from the definition.
\begin{itemize}
\item  $\left<a,a'| b,b'\right>= \left<b,b'| a,a'\right>$.
\item  $\left<a,a'|b,b'\right>= -\left<a',a| b,b'\right>= -\left<a,a'| b',b\right>$.
\item  $\left<a,a|b,b'\right>=0$, $\left<a,a'|b,b\right>=0$.
\item  $\left<a,a'| b,b'\right>+ \left<a',a''| b,b'\right>= \left<a,a''| b,b'\right>$.
\item  $\left<a,b|c,x\right>+\left<b,c|a,x\right>+\left<c,a|b,x\right>=0$.
\end{itemize}

The function
\begin{equation}\label{d_gr-pr}
\left<a|b\right>_c:= \left<a,c|c,b\right>=
\frac{1}{2} \big( d(a,c)+ d(b,c)- d(a,b)\big)
\end{equation}
is the {\sf Gromov product}\label{i_grpr} with respect to metric~$d$.
The triangle inequality implies $\left<b|c\right>_a\ge 0$. 
The two concepts are related by the formula
$$\left<a,b|x,y\right>= \left<b|x\right>_a -\left<b|y\right>_a.$$

\subsection{Generalized metrics}\label{ss_gener-metr}
We will deal with points at infinity, so it is convenient to extend the class of metric spaces.
A {\sf generalized metric space} is a topological space $Y$ with a function
$ d \co Y\times Y\to[0,\infty]$ such that $d(x,y)=d(y,x)$,
$ d(x,y)=0$ iff $x=y$, and $ d(x,z)\le  d(x,y)+ d(y,z)$, for all
$x,y,z\in Y$. Here we use the conventions $r+\infty=\infty$ and $r\le\infty$
for every $r\in[0,\infty]$. For $x\in X$, the {\sf finite component} of $x$
is the set $\{y\in Y\ |\  d(x,y)<\infty\}$. Obviously, $Y$ is the disjoint union of
its finite components, and the restriction of~$ d$ to each finite component is a metric
in the usual sense. Moreover, we require as a part of definition that the topology
on each finite component $V$ of $Y$ coincides with the topology induced by the restriction 
of~$ d$ to~$V$. The function~$ d$ is called a {\sf generalized metric} on~$Y$.
Note that, for a given $ d$, there might be many topologies on~$Y$ that make $(X, d)$ a generalized
metric space.

The main examples to keep in mind:
\begin{itemize}
\item Any metric space is a generalized metric space.
\item $\bar\rr:=[-\infty,\infty]$ with the topology of a closed interval and the obvious
generalized metric that we will denote $|\cdot|$; that is, $|x-y|$ denotes the (generalized)
distance between $x$ and~$y$.
\item Given a hyperbolic graph~$\G$, $(\bar{\G},  d)$ is a generalized metric space,
where $\bar\G=\G\sqcup \p\G$ is the compactification of $\G$ by its ideal boundary 
and $d$ is the obvious extension of the word metric on~$\G$ to $\bar\G$:
we have $ d(x,y)=\infty$ for any $x\in\p\G$ and any $y\in\bar{\G}\setminus\{x\}$.
\end{itemize}

A map $X\to Y$ between generalized metric spaces is called an 
{\sf isometric embedding}\label{i_isom-imbedding}
if it preserves the generalized metric. Note that an isometric embedding must not be
continuous.

\subsection{$^+$equivalence, $^\times$equivalence, $^{\times +}$equivalence}
\label{ss_equivalence}
In this section we introduce convenient equivalence relations of functions that will be
used later in the paper.

For subsets $U\se\bar\rr$ and $V\se\rr$, addition and multiplication can be defined 
in the obvious way: 
$$U+V := \{u+v\ |\ u\in U,\ v\in V\}\quad\mbox{and}\quad 
UV :=\{uv\ |\ u\in U,\ v\in V\}.$$
If $U$ or $V$ consists of just one element, then we write the element instead of the set
notation.

Let $S$ be any set, and $\varphi$ and $\psi$ be functions from~$S$ to~$\bar\rr$.
We say that $\varphi$ and $\psi$ are 
{\sf $^+$equivalent}, denoted $\varphi\simplus \psi$,
\label{i_simplus} if
there exists $B\in [0,\infty)$ such that 
$\varphi(s)\in \psi(s)+ \left[-B, B\right] $
for all~$s\in S$.
They are {\sf $^\times$equivalent} if
there exists $A\in[1,\infty)$ such that 
$\varphi(s)\in [1/A, A]\, \psi(s)$
for all~$s\in S$; and they are {\sf $^{\times +}$equivalent} if
there exist $A\in [1,\infty)$ and $B\in [0,\infty)$ such that 
$\varphi(s)\in [1/A, A]\, \psi(s)+ \left[-B,B\right]$
for all~$s\in S$.
It is left to the reader to check that these are indeed equivalence relations.

We will say that a map $f\co (X_1,d_1)\to (X_2,d_2)$ between 
metric spaces is a {\sf $^+$map}\label{i_+map} if
$d_1(x,y)$ and $d_2(f(x),f(y))$ are $^+$equivalent
as functions of $(x,y)\in X_1\times X_1$.
Similarly for {\sf $^\times$map} and {\sf $^{\times+}$map}.

A $^+$map $f\co (X_1,d_1)\to (X_2,d_2)$
is called a {\sf $^+$isometry}\label{i_+isometry} if there exists a $^+$map
$f\co (X_2,d_2)\to (X_1,d_1)$ such that 
$f\circ g$ and $g\circ f$ are $^+$equivalent to the identity maps.
A $^{\times+}$map $f\co (X_1,d_1)\to (X_2,d_2)$
is called a {\sf $^{\times+}$isometry}\label{i_times+isometry} if there exists a $^{\times+}$map
$f\co (X_2,d_2)\to (X_1,d_1)$ such that 
$f\circ g$ and $g\circ f$ are $^+$equivalent to the identity maps.
$^{\times +}$equivalence of metric spaces 
is the same thing as {\sf quasiisometry}.\label{i_quasiisometry}
Our definitions of equivalences are more general since they allow
considering functions with negative or infinite values.

\begin{lemma}\label{l_plus-eq}
 Suppose that $\varphi$, $\varphi'$, $\psi$, $\psi'$ take values in~$\rr$.
\begin{itemize}
\item [\rm(a)] If $\varphi\simplus \psi$ and $\varphi'\simplus \psi'$, then
$(\varphi+ \varphi')\simplus (\psi+ \psi')$ and $(\varphi- \varphi')\simplus (\psi- \psi')$.
\item [\rm(b)] If $\varphi\simplus \psi$, then $|\varphi|\simplus |\psi|$.
\end{itemize}
\end{lemma}
\noindent The proof follows directly from definitions.

\subsection{Dealing with $\rb$}\label{ss_dealing-with-rb}\ 
Denote $\bar{\rr}:=[-\infty,\infty]$.\label{i_rbar}
Throughout the paper we will consider the notions of addition, subtraction,
multiplication, taking absolute values and distances in $\rb$, generalizing
the corresponding notions in~$\rr$.
For example, 
\begin{gather*}
r+(\pm\infty):=\pm\infty\quad\mbox{for}\quad r\in\rr,
\qquad\qquad (\pm\infty)\cdot l:=\pm\infty\quad\mbox{for}\quad l\in (0,\infty],\\
\infty-\infty:=0,\quad (-\infty)-(-\infty):=0,\quad
    |\infty-\infty|:=0,\quad |(-\infty)-(-\infty)|:=0.
\end{gather*}
Also $e^{\infty}:=\infty$ and $e^{-\infty}:=0$.

\subsection{The smooth-out}\label{ss_theta}
Call a function $\theta\co\rb\to\rb$ {\sf non-expanding} if
$|\theta(t)-\theta(s)|\le |t-s|$ for all $s,t\in\rb$.
Non-expanding functions might be discontinuous: let $\theta$ map
$[-\infty,\infty)$ to $[-\infty,0)$ and $\infty$ to $\infty$.

\begin{lemma}\label{l_theta}
Let $\theta\co \rb\to\rb$ be a continuous non-expanding non-decreasing 
function whose image is an interval $[\al,\be]\se\rb$.
Then the function $\theta'\co \rb\to\rb$,
$$\theta'(t):= \int_{-\infty}^{\infty} \theta(r+t)\frac{e^{-|r|}}{2}\,dr$$
is a well-defined continuous non-expanding non-decreasing map from $\rb$ onto $[\al,\be]$.
If, in addition, $\al <\be$, then $\theta'$ is an increasing homeomorphism
of $\rb$ onto $[\al,\be]$.
\end{lemma}
\proof
If there exists $t\in\rr$ with $\theta(t)=\infty$, then since $\theta$ is non-expanding,
$\theta(\rr)=\{\infty\}$, and since $\theta$ is continuous, $\theta(\rb)=\{\infty\}$.
Then $\theta'(\rb)=\{\infty\}$ and the lemma holds. Similarly for the case $\theta(t)=-\infty$,
so from now on we will assume $\theta(\rr)\se \rr$. Since $\theta$ is non-expanding,
then for each $t\in\rr$, $\theta'(t)$ is a well-defined real number. Since $\theta$ is non-decreasing,
$\theta(-\infty)=\al$, $\theta(\infty)=\be$,
$$\theta'(-\infty)= \int_{-\infty}^{\infty} \theta(r-\infty)\frac{e^{-|r|}}{2}\,dr= 
 \theta(-\infty) \int_{-\infty}^{\infty} \frac{e^{-|r|}}{2}\,dr= \al,$$
and similarly, $\theta'(\infty)=\be$.

$\theta'$ is non-decreasing: if $t\le t'$, then for all $r\in\rr$, $r+t\le r+t'$,
$\theta(r+t)\le \theta(r+t')$, hence $\theta'(t)\le \theta'(t')$. This implies that
$\theta'(\rb)\se [\al,\be]$.

Since $\theta$ is non-expanding,
\begin{eqnarray*}
&& |\theta'(t')-\theta'(t)|\le
\int_{-\infty}^\infty \left| \theta(r+t')-\theta(r+t)\right| \frac{e^{-|r|}}{2}\,dr\\
&& \qquad \le \int_{-\infty}^\infty \left|(r+t')-(r+t)\right| \frac{e^{-|r|}}{2}\,dr\le |t'-t|
\end{eqnarray*}
shows that $\theta'$ is non-expanding, therefore it is continuous on~$\rr$.

Now we show the continuity of~$\theta'$ at~$-\infty$. 
First assume that $\al\in\rr$. For any $\varepsilon>0$ pick $T\in\rr$ such that
$\theta(T)\le \al+\varepsilon/2$ and $R\in[0,\infty)$ such that
$\int_R^\infty re^{-r}\,dr= e^{-R}(R+1)\le \varepsilon$.
Since $\theta$ is non-expanding,
$$\theta(r+t)\le \theta(t)+ |\theta(t+r)-\theta(t)|\le \theta(t)+|r|.$$
We claim that $\theta'([-\infty,T-R])\se [\al,\al+\varepsilon]$.
Indeed, for any $t\in[-\infty,T-R]$, since $\theta$ is non-decreasing,
\begin{eqnarray*}
&& \theta'(t)= \int_{-\infty}^R \theta(r+t) \frac{e^{-|r|}}{2}\,dr
  +\int_R^\infty \theta(r+t) \frac{e^{-|r|}}{2}\,dr\\
&& \le \int_{-\infty}^R \theta\big(R+(T-R)\big) \frac{e^{-|r|}}{2}\,dr
  +\int_R^\infty \big(\theta(t)+|r|\big) \frac{e^{-|r|}}{2}\,dr\\
&&\le \int_{-\infty}^R \theta(T) \frac{e^{-|r|}}{2}\,dr
  +\int_R^\infty \big(\theta(T)+r\big) \frac{e^{-r}}{2}\,dr\\
&& = \theta(T)+ \int_R^\infty r\, \frac{e^{-r}}{2}\,dr
  \le (\al+\varepsilon/2)+\varepsilon/2=\al+\varepsilon.
\end{eqnarray*}
This shows the continuity of $\theta'$ at~$-\infty$ when $\al\in\rr$.
When $\al=-\infty$ the argument is similar:
for any $N$, pick $T$ so that $\theta(T)\le N-1$ and let $R=0$, so that
$\int_0^\infty r e^{-r}\,dx=1$, then deduce that
$\theta'\big([-\infty,T]\big)\se [-\infty,N]$. If $\al=\infty$, then since $\theta$
is non-decreasing, $\theta(\rb)=\{\infty\}$; we dealt with this obvious
case before. The same argument applies to~$\infty$, so $\theta'$ is continuous on~$\rb$.
This implies that $\theta'$ maps~$\rb$ onto $[\al,\be]$.

Now assume $\al<\be$. It remains to show that $\theta'$ is strictly increasing,
and for that it suffices to show that $\theta'$ is such  on~$\rr$.
Take any $t,t'\in\rr$ with $t<t'$. Suppose that for all $r\in\rr$,
$\theta(r+t)=\theta(r+t')$. Then
$\theta\big((t'-t)+t\big)=\theta(t')=\theta(t)$ and inductively, for any positive
integer~$n$,
$$\theta\big( n(t'-t)+t\big)= \theta\big((n-1)(t'-t)+t'\big)=
 \theta\big((n-1)(t'-t)+t\big)=\theta(t).$$
Since $\theta$ is non-decreasing and continuous, this implies that
$\theta$ is a constant function, which is a contradiction with our assumption
$\al<\be$. Therefore there exists $r_0\in\rr$ such that
$\theta(r_0+t')>\theta(r_0+t)$. Since $\theta(r+t')-\theta(r+t)\ge 0$ is continuous in $r$,
$$\theta'(t')-\theta'(t)=
  \int_{-\infty}^\infty \big(\theta(r+t')-\theta(r+t)\big) \frac{e^{-|r|}}{2}\,dr >0.
\eqno{\qed}$$
The lemma says in other words that  the ``operator"
$\theta\mapsto\theta'$ makes homeomorphisms out of non-constant
continuous monotone functions.

\subsection{Auxiliary functions $\theta$ and $\theta'$}\label{ss_aux-fun}
For all $\al,\be\in\rb$ with $\al\le \be$ 
we define a specific function $\theta[\al,\be;\cdot] \co  \bar{\rr}\to [\al,\be]$
that mimics projection of one geodesic into another. 
$$\theta[\al,\be;t]:=
\begin{cases}
\al              & \text{if \ $t\in[-\infty,\al]$}\\
\ t  & \text{if \ $t\in [\al,\be]$}\\
\be                & \text{if \ $t\in[\be, \infty]$.}
\end{cases} $$
In other words, $\theta[\al,\be;\cdot]$ is the obvious extension of the identity
map\newline
$[\al,\be]\to[\al,\be]$. $\theta$ satisfies the following properties.
\begin{itemize}
\item $\theta$ is continuous in three variables.
\item $\theta[\al,\be;\cdot]$ maps
   $\rb$ onto the interval $[\al,\be]$.
\item $\theta[\al,\be;\cdot]$ is non-expanding, that is,  
      $|\theta[\al,\be;t_1]-\theta_0[\al,\be;t_2]|\le |t_1- t_2|.$
\item $\theta[\al,\be;t]$ is increasing in $\al$ and in $\be$.
\item $\theta[-\infty, \infty;\cdot\,]\co  \rb\to\rb$ is the identity map.
\item $\theta$ is {\sf shift-invariant}:\label{i_shift-invariance}
 $\theta[\al+s,\be+s;t+s]=\theta[\al,\be;t]+s$
for all $s\in\rr$.
\end{itemize}

Define a new function $\theta'[\cdot,\cdot\,;\cdot\,]$ as in Lemma~\ref{l_theta}:
$$\theta'[\al,\be;t]:= 
 \int_{-\infty}^{\infty} \theta[\al,\be;r+t]\frac{e^{-|r|}}{2}\,dr.$$
Then $\theta'$ satisfies all the above properties  of $\theta$ and, in addition,
\begin{itemize}
\item for $\al<\be$, $\theta'[\al,\be;\cdot\,]$ is a non-expanding homeomorphism
of $\rb$ onto $[\al,\be]$.
\end{itemize}

The following lemma says that $\theta$ and $\theta'$ are close, i.e.\ the smooth-out
does not change functions much.
\begin{lemma}\label{theta'theta}
$\theta'[\al,\be;t]=\theta[\al,\be;t]+\left(e^{-|t-\al|}-e^{-|t-\be|}\right)/2$ \ for all
$\al,\be\in\rb$, $\al\le\be$, $t\in\rr$.
In particular, $|\theta'[\al,\be;t]- \theta[\al,\be;t]|\le 1$.
\end{lemma}
\begin{proof}
By the definitions of $\theta'$ and $\theta$,
$$\theta'[\al,\be;t]=\int_{-\infty}^{\al-t}\al\frac{e^{-|r|}}{2}\,dr+
\int_{\al-t}^{\be-t} (r+t)\frac{e^{-|r|}}{2}\,dr+
\int_{\be-t}^{\infty}\be\frac{e^{-|r|}}{2}\,dr.$$
Now the statement is proved by brute force calculus in each of 
the following cases:
$t\in(-\infty,\al]$, $t\in[\al,\be]$, $t\in[\be,\infty)$.
\end{proof}

The derivative of $\theta'[\alpha,\beta;\cdot]$ is obtained
by direct calculation:
\begin{equation}\label{derivative}
\frac{\p }{\p t}\theta'[\alpha,\beta;t]= 
\begin{cases}
(e^{t-\al}- e^{t-\beta})/2  &\text{if $t\in (-\infty,\al]$},\\
1- (e^{\al- t}+ e^{t-\beta})/2  &\text{if $t\in [\al,\beta]$},\\
(-e^{\al- t}+ e^{\beta- t})/2  &\text{if $t\in [\beta,\infty)$}.
\end{cases}
\end{equation}
This extends by the same formula to a function
$$\ds\frac{\p}{\p t} \theta'[\al,\beta;\cdot]\co [-\infty,\infty]\to [0,1].$$ 
The function is increasing on $\left[-\infty, \frac{\al+ \beta}{2}\right]$ and decreasing
on $\left[\frac{\al+\beta}{2},\infty\right]$.

\begin{lemma}\label{l_lower-bound}
Suppose that $\al,\beta\in\rb$, $\al\le\be$, $t\in \rr$ and $\epsilon\in [0,(1-e^{\al-\beta})/2]$,
then 
$$\theta'[\al,\beta;t]\in[\al+\epsilon,\beta-\epsilon]
\quad \Rightarrow\quad
\frac{\p }{\p t}\theta'[\alpha,\beta;t]\ge  \epsilon.$$
\end{lemma}

\proof
Fix $\al$, $\beta$ and $\epsilon$.
If $\beta-\al< 2\epsilon$, then the statement is obvious,
so we assume $\beta-\al\ge 2\epsilon$.
Then we can denote $t_\epsilon$ the number such that
$\theta'[\al,\beta;t_\epsilon]=\beta-\epsilon$.
By Lemma~\ref{theta'theta} and the definition of $\theta$,
$$\theta'[\al,\beta;\beta]=\beta- (1-e^{\al-\beta})/2\le
\beta- \epsilon=\theta'[\al,\beta;t_\epsilon].$$
Since $\theta'[\al,\beta;\cdot]$ is increasing, we have
$\beta\le t_\epsilon$. Then by Lemma~\ref{theta'theta},
$$\beta-\epsilon=\theta'[\al,\beta;t_\epsilon]=
\beta+(e^{\al-t_\epsilon}-e^{\beta-t_\epsilon})/2, \text{ \ hence \ }
  (-e^{\al-t_\epsilon}+e^{\be-t_\epsilon})/2\ge\epsilon.$$
(The last inequality indeed holds when $\be=-\infty$ or $\be=\infty$.)
Then by~(\ref{derivative}),
$$\frac{\p }{\p t}\theta'[\alpha,\beta;t_\epsilon]= 
(-e^{\al-t_\epsilon}+ e^{\beta- t_\epsilon})/2\ge\epsilon.$$
Similarly, 
$$\frac{\p }{\p t}\theta'[\alpha,\beta;s_\epsilon]\ge \epsilon,$$
where $s_\epsilon$ is defined by
$\theta'[\al,\beta;s_\epsilon]=\al+\epsilon$.
Now it is a calculus exercise to see from~(\ref{derivative})
that the minimum of the function
$\ds\frac{\p }{\p t}\theta'[\alpha,\beta;\cdot]$
on the interval $[s_\epsilon,t_\epsilon]$ is attained
at the endpoints,
so 
$$\frac{\p }{\p t}\theta'[\alpha,\beta;t]\ge \epsilon\quad
\text{for all \ } t\in [s_\epsilon,t_\epsilon].\eqno{\qed}$$

\section{$\sx$: the symmetric join of $X$}\label{s_symm-join}
\subsection{Symmetric join as a topological space}\label{ss_sj-as-top}
Given a topological space $X$, its usual topological {\sf join} (with itself) is the ``obvious'' union
of intervals connecting pairs of points in $X$. Formally,  
it is the topological space
${X\!\Join\! X}:= X^2\times\rb/\sim$\label{i_join-of-X},
where $(x,y,-\infty)\sim (x,y',-\infty)$ and $(x,y,\infty)\sim (x',y,\infty)$
for all $x,x',y,y'\in X$.
We will call the further quotient 
$$\sx:= \jx /\sim\label{i_sx}$$ 
by the equivalence
relation $(x,x,s)\sim (x,x,t)$ for all $x\in X$, $s,t\in\rb$, 
the {\sf symmetric join}\label{i_symm-join-of-X}
 of $X$.
That is, each interval connecting a point $x\in X$ to itself degenerates to a point in~$\sx$.
The topology on $\sx$ is induced by the double quotient
$X^2\times \rb\twoheadrightarrow X\!\Join\! X\twoheadrightarrow \sx$.

Applying the two quotients above in the reverse order gives an equivalent, and more convenient
for our purposes, definition of symmetric join: 
\begin{itemize}
\item denote
$\diamondx:= \bigsqcup\{[a,b]\ |\ a,b\in X\}$,\label{i_diamondx}
where $[a,b]$ is a topological copy of $\rb$ if $a\not= b$ and
a point if $a=b$; 
\item let $\sx$ be the quotient of $\diamondx$
by the equivalence relation identifying the left endpoint of each $[a,b]$ with $[a,a]$
and the right endpoint with $[b,b]$. 
\end{itemize}
The topology on $\sx$ is induced by the obvious double surjection
\begin{equation}
X^2\times \rb\twoheadrightarrow \diamondx\twoheadrightarrow \sx
\end{equation}
in which the lines $\{a\}\times\{b\}\times\rb$ map onto intervals~$[a,b]$.

The canonical embedding $X\hookrightarrow \sx$, $a\mapsto [a,a]$, is a homeomorphism
onto its image. We will identify $X$ with its image in~$\sx$. 
The sets $\{a\}\times \{b\}\times \rb$ in $X^2\times\rb$, and their images in $\jx$,
$\diamondx$, and $\sx$, will be called {\sf lines}. In particular, each point $a\in X$
is the line $[a,a]$. Denote 
$$\starx:= \sx \setminus X\se \sx.$$
This is the {\sf open symmetric join}\label{i_open-symm-join-of-X}
 of~$X$. 
The topology on $\starx$ is induced from its inclusion into $\sx$.
Let $\Delta$ be the diagonal in~$X^2$.
The above double surjection
restricts to a {\em bijection} $(X^2\setminus\Delta)\times\rr\to\starx$, so $\starx$ is topologically
$(X^2\setminus\Delta)\times\rr$.

When $X$ is a metric space, we will define actions 
by $\rr$, $\zz_2$ and $\mathrm{Isom}(X)$ on $\sx$, and equip $\sx$ with a $\zz_2$-- and
$\mathrm{Isom}(X)$--invariant metric  which, as we will see later, also behaves nicely with
respect to the $\rr$--action. That is, the symmetric join of~$X$ will become
an example of a {\sf metric join}\label{i_metric-join}.

\subsection{Parametrizations of $\sx$}
\label{ss_parametrizations}
Given a metric space~$(X, d)$, first we want to put a metric on each line
in $\sx$. We will do this by identifying each line with a subinterval
of~$\rr$.

Define the functions $\left<\cdot,\cdot|\cdot,\cdot\right>$
and $\left<\cdot|\cdot\right>_{\cdot}$ in~$X$ as in~\ref{ss_dd-gp},
using the metric on~$X$, and fix a basepoint $x_0\in X$.
For every pair~$(a,b)\in X^2$, denote
$$\al:= -\left<b|x_0\right>_a\qquad\mbox{and}\qquad \be:=\left<a|x_0\right>_b,$$
these are the {\sf end-point coordinates}\label{i_end-point-coord} of the pair $(a,b)$
(with respect to $x_0$).
Let $[\![a,b]\!]=[\![a,b]\!]_{x_0}$\label{i_interval[[]]}
 be a copy of the interval 
$[\al, \be]\se \rr$, and $]\!]a,b[\![$ its interior. 
Note that the interval $[\al, \be]\se\rr$ always contains 0
and its length is
$$\be-\al=\left<a|x_0\right>_b+ \left<b|x_0\right>_a=  d(a,b)\ge 0.$$
It is convenient to think of $[\![a,b]\!]$ as a formal geodesic connecting $a$ to~$b$;
it indeed degenerates to a point when $a=b$.
The convenience of this parametrization of $[\![a,b]\!]$
will be clear later when we define 
the $\Isom(X)$--action and the projection of~$X$ to $[\![a,b]\!]$;
the basepoint~$x_0$ will project to $0\in [\![a,b]\!]$, and $a$ and $b$ will project
to the endpoints of~$[\![a,b]\!]$.
Also this parametrization will be needed to extend things to infinity
when $X$ is a hyperbolic complex.

Define functions
\begin{eqnarray}\label{def-[[]]-[[]]'}
&& [\![a,b;\cdot]\!]=[\![a,b;\cdot]\!]_{x_0}\co  \rb\to [\![a,b]\!]\quad \mbox{by}\quad
[\![a,b;t]\!]:= \theta[\al, \be; t]\\
\mbox{and}&& [\![a,b;\cdot]\!]'= [\![a,b;\cdot]\!]'_{x_0}\co \rb\to [\![a,b]\!]\quad \mbox{by}\quad
[\![a,b;t]\!]':= \theta'[\al, \be; t],
\end{eqnarray}
where $\theta$ and $\theta'$ are the functions from~\ref{ss_aux-fun}. 
By the definition of $\theta$, 
$[\![a,b;\cdot]\!]$ is just the obvious extension of the identity map of
$[\![a,b]\!]$, so it is a non-expanding surjection. 
A real number $t$ will be called {\sf appropriate}\label{i_appropriate}
for $[\![a,b]\!]$ if $t\in\big[\al, \be\big]$.
A number $t$ appropriate for $[\![a,b]\!]$ will be called the 
{\sf $x_0$--coordinate}\label{i_x0-coordinate}
of the point
$[\![a,b;t]\!]$ in $[\![a,b]\!]$, or simply the {\sf coordinate} when~$x_0$ is understood.
If $(\al,\be)$ are the endpoint coordinates of $(a,b)$, then
\begin{equation}\label{abt'}
[\![a,b;t]\!]'= [\![a,b;\theta'[\al,\be;t] ]\!].
\end{equation}
This holds just because $[\![a,b;\cdot]\!]$ is identity on appropriate numbers.

$[\![a,b;\cdot]\!]'$ is the smooth-out of $[\![a,b;\cdot]\!]$. By Lemma~\ref{l_theta},
$[\![a,b;\cdot]\!]'$ is a non-expanding surjection, and it is an increasing homeomorphism
when $a\not=b$. Denote
\begin{equation}\label{i_diamondxxo}
\diamondxxo:= \bigsqcup \big\{ [\![a,b]\!]\ |\ a,b\in X \big\},\qquad\qquad
\sxxo:= \diamondxxo/\sim,
\end{equation}
where the equivalence relation $\sim$ identifies the left endpoint
of each interval $[\![a,b]\!]$ with the point $[\![a,a]\!]$ and
the right endpoint with $[\![b,b]\!]$.
Thus we view $X$ as embedded in $\sxxo$ via $a\mapsto [\![a,a]\!]$,
and abusing notations we view $[\![a,b]\!]$ as a subset both
of $\diamondxxo$ and of $\sxxo$.

We denote $|\cdot|$ the standard metric on each $[\![a,b]\!]$,
that is $|x-y|$ is the distance between $x$ and $y$ in $[\![a,b]\!]$.

\subsection{The models $\sx$ and $\sxxo$}\label{ss_two-models-sx-sxxo}
The above maps $[\![a,b;\cdot]\!]'$ induce surjections
\begin{equation}\label{i_[[...]]'}
[\![\cdot,\cdot\,;\cdot]\!]'\co  X^2\times\rb\to \diamondxxo\qquad\mbox{and}\qquad
[\![\cdot,\cdot\,;\cdot]\!]'\co  X^2\times\rb\to\sxxo,
\end{equation}
and passing to quotients, bijections
\begin{equation}\label{bijection-sx-sxxo}
[\![\cdot,\cdot\,;\cdot]\!]'\co  \diamondx\to\diamondxxo\qquad
\mbox{and}\qquad [\![\cdot,\cdot\,;\cdot]\!]'\co \sx\to \sxxo.
\end{equation}
To summarize:
\begin{itemize}
\item $\sx$ and $\sxxo$ are two different models, or parametrizations,
  of the symmetric join;
\item $\sx$ has a natural topology, and each line is topologically parameterized by~$\rb$;
\item $\sxxo$ is just a set, but lines in $\sxxo$ are
metric spaces which are isometrically parameterized by closed subintervals of~$\rr$;
\item the two models are identified by the bijection $[\![\cdot,\cdot\,;\cdot]\!]'$
in~(\ref{bijection-sx-sxxo}).
\end{itemize}
We induce the topology on $\sxxo$ from $\sx$ by this bijection. 
Equivalently, the topology on $\sxxo$ is induced by the surjection
$[\![\cdot,\cdot\,;\cdot]\!]'\co X^2\times\rb\to\sxxo$.
This topology is consistent
with the metric topology on each line in $\sxxo$ because each 
$[\![a,b;\cdot]\!]'\co \rb\to [\![a,b]\!]$ is a homeomorphism for $a\not= b$.

The bijection in~(\ref{bijection-sx-sxxo}) allows us to identify 
$\sx$ with $\sxxo$, so we will often use the simpler notation
$\sx$ instead of $\sxxo$.

\subsection{The models $\starx$ and $\starxxo$}
\label{ss_two-models-starx-starxxo}
There are similar models for the open symmetric join: denote
\begin{equation}\label{i_starxxo}
\starxxo:= \bigcup\big\{]\!]a,b[\![\se\sxxo\ |\ (a,b)\in X^2\big\}=
 (\sxxo)\setminus X\se \sxxo.
\end{equation}
Define the topology on $\starxxo$ by its inclusion into $\sxxo$.

Let $\Delta$ be the diagonal of $X^2$.
We can assume $(a,b)\in X^2\setminus\Delta$ above, since $]\!]a,b[\![$ is empty otherwise.
Let $]\!]a,b;\cdot [\!['\co \rr\to ]\!]a,b[\![$ be the restriction of $[\![a,b;\cdot]\!]'$
to the interiors of the intervals; this restriction is a homeomorphism.
In the quotient this gives {\em bijections}
\begin{equation}
]\!]\cdot,\cdot\,;\cdot{[\!['}\co  (X^2\setminus\Delta)\times\rr\to\starxxo\qquad
\mbox{and}\qquad
]\!]\cdot,\cdot\,;\cdot{[\!['}\co  \starx\to\starxxo
\end{equation}
which induce the same topology on~$\starxxo$ as the one described above.

Our goal is to equip the symmetric join with a metric and three actions.
In the process we will be using the two models interchangeably.

\subsection{Change of basepoint in $\sx$}\label{ss_change-of-var}
The parametrization $[\![a,b;\cdot]\!]= [\![a,b;\cdot]\!]_{x_0}$
of each line $[\![a,b]\!]= [\![a,b]\!]_{x_0}$ depends on the choice
of~$x_0\in X$, but the isometry type of $[\![a,b]\!]$ does not.
Another choice of basepoint, $x_1$
gives another interval $[\![a,b]\!]_{x_1}:=[-\left<b|x_1\right>_a, \left<a|x_1\right>_b]$
of length~$d(a,b)$ and another identity map
$$[\![a,b;\cdot]\!]_{x_1}\co 
  [-\left<b|x_1\right>_a, \left<a|x_1\right>_b]\to [\![a,b]\!]_{x_1}.$$

We will always identify two such parametrizations by the unique isometry between
$[\![a,b]\!]_{x_0}$ and $[\![a,b]\!]_{x_1}$ for all~$a,b\in X$.
Since the left endpoints\newline
$[\![a,b;-\left<b|x_0\right>_a]\!]_{x_0}$ and
$[\![a,b;-\left<b|x_1\right>_a]\!]_{x_1}$
must be identified, the explicit formula is
\begin{equation}\label{change-basepoint}
[\![a,b;t]\!]_{x_1}= [\![a,b;t+\left<b|x_1\right>_a -\left<b|x_0\right>_a]\!]_{x_0}=
  [\![a,b;t+ \left<a,b|x_1,x_0\right>]\!]_{x_0},\quad t\in\rb.
\end{equation}

\subsection{The $\rr$--action on~$\sx$}\label{ss_r-action}
Let $X$ be an arbitrary metric space. 
The shift action
$$+\co \rr\times\rb\to \rb,\qquad (r,t)\mapsto r^+ t,\quad
\mbox{where}\quad r^+t:=r+t$$
with the convention $r\pm\infty=\pm\infty$,
induces the obvious $\rr$--action on 
$X^2\!\times\rb$ by translations in the $\rb$--coordinate,
and, passing to the quotient, the $\rr$--action on $\sx$:
$$+\co \rr\times (\sx)\to \sx,\qquad (r,z)\mapsto r^+ z.$$
This action preserves lines and fixes $X$ pointwise.
The $\rr$--orbits in~$\sx$ are exactly the points of $X$ and 
the interiors of the lines in $\sx$.
The bijection $[\![\cdot, \cdot\,;\cdot]\!]'$ transfers this further to the action 
$$+\co \rr\times (\sxxo)\to (\sxxo), \quad (r,z)\mapsto r^+z.$$
The explicit formula for the action is \ $r^+[\![a,b;t]\!]':=[\![a,b;r+t]\!]'$.

\begin{lemma}\label{l_r-action}
For each $x\in[\![a,b]\!]\se \sxxo$, the $\rr$--orbit map $\rr\to([\![a,b]\!],|\cdot|)$,
   $r\mapsto r^+ x$,
   is non-expanding, therefore continuous.
   In addition, if $x\in]\!]a,b[\![$, then the orbit map
   is an increasing homeomorphism of $\rr$ onto $]\!]a,b[\![$.
\end{lemma}
\begin{proof}
{\bf (a)}\qua $x=[\![a,b;t]\!]'$ for some $t\in\rb$.
Since $[\![a,b;\cdot]\!]'$ is non-expanding,
$$|r_1^+x- r_2^+x| = \big|[\![a,b;r_1+t]\!]'- [\![a,b;r_1+t]\!]'\big|\le
|(r_1+t)-(r_2+t)|=|r_1- r_2|,$$
i.e.\ the orbit map is non-expanding.
If $x\in ]\!]a,b[\![$ then $a\not= b$, $t\in\rr$ and
the orbit map $r\mapsto ]\!]a,b;r+t[\!['$ is an increasing homeomorphism 
$\rr\to]\!]a,b[\![$ because it is the composition of increasing homeomorphisms
$\rr\to\rr\to ]\!]a,b[\![$, $r\mapsto r+t\mapsto ]\!]a,b;r+t[\!['$.
\end{proof}

\subsection{The $\zz_2$--action on $\sx$}\label{ss_z2-action}
The map $\star\co \diamondx\to\diamondx$ given by $[\![a,b;t{]\!]}^\star:= [\![b,a;-t]\!]$
for $a,b\in X$ and an appropriate $t$ is a well-defined involution of $\diamondx$.
Note that $[\![a,a{]\!]}^\star=[\![a,a]\!]$. 
Also if $[\![a,b;t]\!]$ is the left endpoint of $[\![a,b]\!]$ (which is identified with
$[\![a,a]\!]$ in $\sx$) then necessarily $t=-\left<b|x_0\right>_a$,
so $[\![a,b;t{]\!]}^\star= [\![b,a; \left<b|x_0\right>_a]\!]$ is the right endpoint of $[\![b,a]\!]$
(which is also identified with $[\![a,a]\!]$). Therefore the same formula gives an involution
$\star\co \sx\to\sx$ in the quotient,  and $\star$ fixes~$X$ pointwise.

\subsection{The $\mathrm{Isom}(X)$--action on $\sx$}\label{ss_isom-action}
$\mathrm{Isom}(X)$ acts on $\diamondxxo$ by
\begin{eqnarray}
\label{Isom-action}
&& g\, [\![a,b;t]\!]\ :=\ [\![ga,gb; t+\left<a,b|x_0,g\inv x_0\right>]\!],\\
\nonumber && g\in \mathrm{Isom}(X),\quad (a,b)\in X^2,\quad
t\in [-\left<b|x_0\right>_a, \left<a|x_0\right>_b].
\end{eqnarray}
One checks that if $t$ runs through $[-\left<b|x_0\right>_a, \left<a|x_0\right>_b]$ then
$t+\left<a,b|x_0,g\inv x_0\right>$ runs through 
$[-\left<gb|x_0\right>_{ga}, \left<ga|x_0\right>_{gb}]$, so $g$ maps each line
$[\![a,b]\!]$ isometrically onto the line $[\![ga,gb]\!]$.

This is indeed an action since for $\mathrm{id}\in \mathrm{Isom}(X)$,
$$\mathrm{id}\, [\![a,b;t]\!]=
[\![\mathrm{id}\, a,\mathrm{id}\, b; t+\left<a,b|x_0, \mathrm{id}\, x_0\right>]\!]= [\![a,b;t]\!]$$ 
 and for $f,g\in \mathrm{Isom}(X)$,
\begin{eqnarray*}
&& f\big(g\,[\![a,b;t]\!]\big) =
   f\,[\![ga,gb; t+\left<a,b|x_0,g\inv x_0\right>]\!]\\
&& = [\![fga,fgb; t+\left<a,b|x_0,g\inv x_0\right>+
                 \left<ga,gb|x_0,f\inv x_0\right>]\!]\\
&& = [\![fga,fgb; t+\left<a,b|x_0,g\inv x_0\right>+
                 \left<a,b|g\inv x_0,g\inv f\inv x_0\right>]\!]\\
&& = [\![fga,fgb; t+\left<a,b|x_0,g\inv f\inv x_0\right>]\!]\\
&& = [\![fga,fgb; t+\left<a,b|x_0,(fg)\inv x_0\right>]\!] = (fg) [\![a,b; t]\!].
\end{eqnarray*}
Since $g$ sends the left endpoint of $[\![a,b]\!]$ to the left endpoint
of $[\![ga,gb]\!]$, and similarly for the right endpoints, the above action descends to an action
on $\sxxo$, given by the same formula~(\ref{Isom-action}), and therefore
to an action on~$\sx$.

\begin{lemma}\label{l_actions-prop}
Let $X$ be an arbitrary metric space.
\begin{itemize}
\item [\rm(a)] The $\mathrm{Isom}(X)$-- and $\rr$--actions on $\sx$ commute.
\item [\rm(b)] The $\mathrm{Isom}(X)$-- and $\zz_2$--actions on $\sx$ commute.
\item [\rm(c)] The $\zz_2$--action anticommutes with the $\rr$--action on~$\sx$:\newline
$(r^+x)^\star=(-r)^+x^\star$ for $x\in \sx$ and $r\in\rr.$
\item [\rm(d)] All the three actions on $\sx$ map lines onto lines
  and are independent of~$x_0$. The $\zz_2$-- and $\rr$--actions fix~$X$ pointwise.
\item [\rm(e)] The $\mathrm{Isom}(X)$--action on $\sx$ is an extension of
the $\mathrm{Isom}(X)$--action on~$X$.
\end{itemize}
\end{lemma}
\begin{proof}
{\bf (a)}\qua The difficulty is that the $\Isom(X)$-- and $\rr$--actions are defined
with respect to different parametrizations, $[\![\cdot, \cdot\,;\cdot]\!]$
and $[\![\cdot, \cdot\,;\cdot]\!]'$, respectively. But the relation~(\ref{abt'})
between the two parametrizations and the shift invariance of $\theta'$
will suffice for the proof. Pick any point $[\![a,b;t]\!]'\in \sx$ and $g\in\Isom(X)$
and denote $s:=\left<a,b|x_0,g\inv x_0\right>$.
By direct calculation, if $(\al,\be)$ are the endpoint coordinates of $(a,b)$, then
$(\al+s,\be+s)$ are the endpoint coordinates of $(ga, gb)$. Hence
for $r\in\rr$ and $g\in\Isom(X)$,
\begin{eqnarray*}
&& r^+\big(g\,[\![a,b;t]\!]'\big)= r^+\big(g\,[\![a,b;\theta'[\al,\be;t]]\!]\big)=
  r^+[\![ga,gb;\theta'[\al,\be;t]+s]\!]\\
&&= r^+[\![ga,gb;\theta'[\al+s,\be+s;t+s] ]\!]=
  r^+[\![ga,gb;t+s ]\!]'= [\![ga,gb;r+t+s ]\!]'\\
&& = [\![ga,gb;\theta'[\al+s,\be+s;r+t+s] ]\!]=
  [\![ga,gb;\theta'[\al,\be;r+t]+s ]\!]\\
&& = g\,[\![a,b;\theta'[\al,\be;r+t] ]\!]=
  g\,[\![a,b;r+t]\!]'= g\big(r^+[\![a,b;t]\!]'\big).
\end{eqnarray*}

{\bf (d)}\qua
The involution $\star$ is independent of $x_0$ because for another $x_1\in X$,
by the change of basepoint formula~(\ref{change-basepoint}),
\begin{eqnarray*}
&& [\![a,b;r]\!]_{x_1}^\star= [\![a,b;r+\left<a,b|x_1,x_0\right>]\!]_{x_0}^\star= 
[\![b,a;-r-\left<a,b|x_1,x_0\right>]\!]_{x_0} \\
&& = [\![b,a;-r-\left<a,b|x_1,x_0\right>+\left<b,a|x_0,x_1\right>]\!]_{x_1}=
 [\![b,a;-r]\!]_{x_1}.
\end{eqnarray*}
Similarly, one checks using~(\ref{change-basepoint}) and the $\Isom(X)$--invariance of 
the double difference that 
$$g[\![a,b;r]\!]_{x_1}= [\![ga,gb; r+ \left<a,b|x_1, g\inv x_1\right>]\!]_{x_1},$$
i.e.\ the $\Isom(X)$--action does not depend on the choice of~$x_0$.
For the $\rr$--action, use~(\ref{change-basepoint}) and the shift invariance of $\theta'$ to prove
$r^+[\![a,b;t]\!]_{x_1}=[\![a,b;r+t]\!]_{x_1}$.

(b), (c), (e) and the rest of (d) directly follow from definitions.
\end{proof}

\subsection{Projecting $X$ to lines in $\sx$}
\label{ss_projX}
Our metric join construction will work for an arbitrary metric space, 
but for inspiration consider first the classical hyperbolic space~$\hh^n$, or more generally,
a $\mathrm{CAT}(-1)$--space. Let $a,a',b,b'\in\p\hh^n$
and denote $[b,b'|a]$ the nearest-point projection
of $a$ to the geodesic from $b$ to $b'$. 
Then the double difference $\left<a,a'|b',b\right>$
indeed makes sense when $a,a',b,b'$ are in the boundary of $\hh^n$, and
it equals the directed distance from $[b,b'|a]$ to $[b,b'|a']$
along the geodesic from~$b$ to~$b'$, or, symmetrically,
the directed distance from $[a,a'|b]$ to $[a,a'|b']$
along the geodesic from~$a$ to~$a'$ (see for example~\cite[1.3]{Bourdon1996}).
We will use this observation in the constructions that follow
(though the endpoints of geodesics will be in the metric space rather than in the ideal boundary).

Again we fix a basepoint $x_0\in X$.
\begin{ddd}\label{d_projX}
Let $a,a',b\in X$.
The {\sf coordinate of the projection} of $b\in X$ to the line
$[\![a,a']\!]\se \sxxo$ is
   $\left< a,a'|b \right>:= \left<a,a'|b,x_0\right>$
and the {\sf projection} of~$b$ to the line $[\![a,a']\!]$ is 
   $[\![a,a'|b]\!]:= [\![a,a'; \left<a,a'|b\right> ]\!]$.
\end{ddd}
\noindent Since $[\![a,a']\!]$ is identified with an interval in~$\rr$, $\left< a,a'|b \right>$
is essentially the same thing as $[\![a,a'|b]\!]$. We only use two different notations
to emphasize that $\left< a,a'|b \right>$ means a real number and $[\![a,a'|b]\!]$
represents a point in the metric space $[\![a,a']\!]$.

The parametrization $[\![a,a'\,;\cdot]\!]$ and the above definitions are chosen so that
in particular,
$$[\![a,a'|a]\!]= [\![a,a';\left<a,a'|a\right>]\!]= [\![a,a';\left<a,a'|a,x_0\right>]\!]=
  [\![a,a';-\left<a'|x_0\right>_a]\!] =a$$
(recall that we do not distinguish between $a\in X$
and the left endpoint of $[\![a,a']\!]\se \sxxo$)
and similarly $[\![a,a'|a']\!]=a'$, that is $a$ and $a'$ project to themselves.
Note also that $[\![a,a'|x_0]\!]= [\![a,a'; \left<a,a'|x_0\right> ]\!]= [\![a,a';0]\!]$;
in other words, the basepoint $x_0$ always project to the ``origin" $[\![a,a';0]\!]$
of~$[\![a,a']\!]$.

The following lemma gives an equivalent description of the projection.
\begin{lemma}\label{l_projection-other-def}
For $a,a',b\in X$, $[\![a,a'|b]\!]$ is the unique point~$c$ in $[\![a,a']\!]$ that satisfies
$$|a-c|-|c-a'|=  d(a,b)- d(b,a').$$
\end{lemma}
\proof
\begin{eqnarray*}
&& \big| a- [\![a,a'|b]\!]\big|= \big| [\![a,a'|a]\!]- [\![a,a'|b]\!]\big|\\
&& \qquad = \big| \left<a,a'|a\right>- \left<a,a'|b\right>\big|
 = |\left<a,a'|a,b\right>|= \left<a'|b\right>_a,
\end{eqnarray*}
and similarly, $\big| [\![a,a'|b]\!]- a'\big|= \left<a|b\right>_{a'}$,
hence 
$$\big| a- [\![a,a'|b]\!]\big|- \big| [\![a,a'|b]\!]- a'\big|=
\left<a'|b\right>_a+ \left<a|b\right>_{a'}= d(a,b)- d(a',b).\eqno{\qed}$$

For any $a,a',b,b'\in X$,
$$\left<a,a'|b'\right>-\left<a,a'|b\right>=
  \left<a,a'|b',x_0\right>- \left<a,a'|b,x_0\right>= \left<a,a'|b',b\right>,$$
i.e.\ the meaning of $\left<a,a'|b',b\right>$ now is the difference of the coordinates
of the projections of $b'$ and $b$ to $[\![a,a']\!]$.

\begin{lemma}
\label{l_projection} 
\begin{itemize}
\item [\rm(a)]  The projection map $[\![\cdot,\cdot|\cdot]\!]$ is independent of the choice
    of basepoint~$x_0$.
\item [\rm(b)]  $[\![\cdot,\cdot|\cdot]\!]$ is $\mathrm{Isom}(X)$--invariant, i.e.\
   $$g\,[\![a,a'|b]\!]=[\![ga,ga'|gb]\!]\quad\mbox{for}\ \ a,a',b\in X,\ \ g\in \mathrm{Isom{(X)}}.$$
\item [\rm(c)] The projection map relates to the $\zz_2$--action on~$\sx$ by the formula
   $${[\![}a,a'|b{]\!]}^\star=[\![a',a|b]\!],\quad a,a',b\in X.$$
\end{itemize}
\end{lemma}
\proof
(a) and (b) follow from Lemma~\ref{l_projection-other-def}, and (c) follows from definitions:
$${[\![}a,a'|b{]\!]}^\star= {[\![}a,a'|\left<a,a'|b\right>{]\!]}^\star=
[\![a',a|-\left<a,a'|b\right>]\!]= [\![a',a|\left<a',a|b\right>]\!]= 
[\![a',a|b]\!].\eqno{\qed}$$

\subsection{Projection and change of basepoint}
\label{ss_pr-ch-of-basepoint}
The change of basepoint formula 
(\ref{change-basepoint}) impies that 
$$[\![a,b;0]\!]_{x_1}=
  [\![a,b;\left<a,b|x_1 \right>]\!]=[\![a,b|x_1]\!].$$
In other words, $[\![a,b;\cdot]\!]_{x_1}$ is the isometric orientation-preserving
reparametrization of $[\![a,b]\!]$ whose origin $[\![a,b;0]\!]_{x_1}$ is the projection of~$x_1$
to~$[\![a,b]\!]$.

\section{The metric~$d_*$ on $\sx$}\label{s_metric-on-sj}
\subsection{The cocycle $\becross$ in $\sx$}\label{ss_the-cocycle}
\begin{ddd}\label{d_cocycle-inX}
Let $X$ be any metric space.
For $u\in X$ and $x=[\![a,a';s]\!]\in\sx$, where $a,a'\in X$ and $s$ is appropriate, 
let
\ $\ell(u,x):= \left<a|a'\right>_u + \left|s-\left<a,a'|u\right>\right|$\label{ell-ux-def}\newline
(see Figure~\ref{f_ell}).

For $u\in X$, $x,y\in\sx$, let
\ $\be_u^{\scriptscriptstyle\times\!}(x,y):=\ell(u,x)-\ell(u,y)$\label{becross-def}.
\end{ddd}
\noindent This gives a function $\becross\co  X\times(\sx)^2\to \rr$ of three variables $(u,x,y)$.
The definition of $\ell$ mimics the case of a tree: if the imaginary triangle $\{u,a,a'\}$ was
degenerate to a tripod, then $\ell(u,x)$ would be exactly the distance between $u$ and $x$.
Note also that when $x\in X$, $\ell(u,x)=d(u,x)$, so the restriction of $\becross$
to $X\times X^2$ is the {\sf distance cocycle}:
 $\becross_u(x,y)=d(u,x)- d(u,y)$.
\setlength{\unitlength}{1cm} 
\begin{figure}[ht!]
  \begin{center}
   \begin{picture}(7,4)

\put(0,0){\circle*{0.1}}
\put(-0.4,0.1){\footnotesize$u$}
\put(5,4){\circle*{0.1}}
\put(5.2,4){\footnotesize$a$}
\put(6,0){\circle*{0.1}}
\put(6.2,0){\footnotesize$a'$}

\qbezier[100](0,0)(5,0)(5,4)


\qbezier(0,0)(2.3,0)(4.7,0)
\qbezier[20](4.7,0)(5.4,0)(6,0)


\qbezier[20](5,4)(5,3.5)(5.1,2.7)
\qbezier(5.1,2.7)(5.2,2)(5.47,1.2)
\qbezier[20](5.47,1.2)(5.68,0.6)(6,0)

\put(4.7,0){\circle*{0.1}}
\put(4.18,1.4){\circle*{0.1}}
\put(5.47,1.2){\circle*{0.1}}
\put(5.7,1.2){\footnotesize$[\![a,a'|u]\!]$}
\put(5.1,2.7){\circle*{0.1}}
\put(5.4,2.7){\footnotesize$x=[\![a,a';s]\!]$}
   \end{picture}
  \end{center}
\caption{$\ell(u,x)$}
\label{f_ell}
\end{figure}
\begin{ttt}\label{t_becross}$\phantom{99}$
\begin{itemize}
\item [\rm(a)] The above functions $\ell\co X\times(\sx)\to [0,\infty)$ and
   $\becross\co  X\times(\sx)^2\to \rr$ are well-defined, independent of~$x_0$,
and Lipschitz in the first variable $u\in X$ for each fixed~$(x,y)$.
\item [\rm(b)] $\be_u^{\scriptscriptstyle\times\!}$ satisfies the cocycle condition \ 
$\be_u^{\scriptscriptstyle\times\!}(x,y)+ \be_u^{\scriptscriptstyle\times\!}(y,z)=\be_u^{\scriptscriptstyle\times\!}(x,z)$.
\item [\rm(c)] $\be_u^{\scriptscriptstyle\times\!}$ is $\zz_2$--invariant in each variable:
$\be_u^{\scriptscriptstyle\times\!}(x,y)=\be_u^{\scriptscriptstyle\times\!}(x^\star,y)=\be_u^{\scriptscriptstyle\times\!}(x,y^{\star})$.
\item [\rm(d)] $\becross$ is $\Isom(X)$--invariant:
  $\be_{gu}^{\scriptscriptstyle\times\!}(gx,gy)=\be_u^{\scriptscriptstyle\times\!}(x,y)$ for $g\in\Isom(X)$.
\end{itemize}
\end{ttt}
\begin{proof}
\noindent {\bf (a)}\qua 
It suffices to show (a) for $\ell$.
Recall from section~\ref{s_symm-join} that $\sx$ is a quotient of $\diamondx$, that is
each $x\in X\se \sx$ can be represented as $[\![a,a';s]\!]$ in many ways.
It is to be shown that the formula for $\ell$ does not depend on the representations
of~$x\in\sx$. 

Assume $[\![a,x;s]\!]$ is the right endpoint of $[\![a,x]\!]\se\diamondx$ and $s$
is appropriate, then necessarily
$s=\left<a|x_0\right>_x$. $[\![a,x;s]\!]$ and $[\![x,x;0]\!]$
represent the same point in $X\se\sx$.
\begin{eqnarray*}
&&\ell(u,[\![a,x;s]\!])=
\left<a|x\right>_u+ \big|\left<a|x_0\right>_x- \left<a,x|u\right> \big|\\
&&  = \left<a|x\right>_u+ |\left<a,x|x,x_0\right>- \left<a,x|u,x_0\right>|
 = \left<a|x\right>_u+ |\left<a,x|x,u\right>|\\
&& =  \left<a|x\right>_u+ \left<a|u\right>_x= d(x,u)
 =\left<x|x\right>_u+ \big|0-\left<x,x|u\right>\big|= \ell(u,[\![x,x;0]\!]),
\end{eqnarray*}
so $\ell$ is well-defined on the level of the quotient~$\sx$.

$s$ is the $x_0$--coordinate of $x$ in $[\![a,a']\!]$, then by~(\ref{change-basepoint})
for another basepoint $x_1$, the $x_1$--coordinate of $x$ is
$s+\left<a,a'|x_0,x_1\right>$. The identity
$$|s-\left<a,a'|u,x_0\right>|= |(s+\left<a,a'|x_0,x_1\right>)- \left<a,a'|u,x_1\right>|$$
shows that $\ell$ is independent of~$x_0$.

By Lemma~\ref{l_projection-other-def},
$\ell(u,x)$ is uniquely determined by the distances between points $u$, $a$, $a'$,
$[\![a,a'|u]\!]$ and $x$. These distances are preserved under isometries of $X$,
hence $\ell$ is $\Isom(X)$--invariant.

By the triangle inequality,
\begin{eqnarray*}
&& |\ell(u,x)-\ell(v,x)|=
  \left|\left<a|a'\right>_u - \left<a|a'\right>_v\right|+
  \Big|\left|s-\left<a,a'|u\right>\right|- \left|s-\left<a,a'|v\right>\right|\Big|\\
&& \le d(u,v)+\left|\left<a,a'|v\right>- \left<a,a'|u\right>\right|=
   d(u,v)+\left|\left<a,a'|v,u\right>\right|\le 2 d(u,v),
\end{eqnarray*}
so $\ell(u,x)$ is Lipschitz in the variable~$u$.

Parts (b) and (c) of the lemma are straightforward from definitions,
(d) follows from definitions and the change of basepoint formula~(\ref{change-basepoint}).
\end{proof}

\begin{lemma}\label{l_beuxy}
Let $u,a,a'\in X$, $x=[\![a,a';s]\!]$, $y=[\![a,a';t]\!]$ with appropriate $s$ and $t$, then
\begin{itemize}
\item [\rm(a)] $|\be_u^{\scriptscriptstyle\times\!}(x,y)|\le |s-t|\le d(a,a')$, \ $\be_a^{\scriptscriptstyle\times\!}(x,y)= s-t$, 
\item [\rm(b)] $|\be_u^{\scriptscriptstyle\times\!}(a,a')|\le d(a,a')$, \ $\be_a^{\scriptscriptstyle\times\!}(a,a')=- d(a,a')$.
\end{itemize}
\end{lemma}
\proof
From the definition of $\becross$,
\begin{eqnarray*}
&&\be_u^{\scriptscriptstyle\times\!}(x,y)= |s- \left<a,a'|u\right>| - |t- \left<a,a'|u\right>|\le |s - t|\le d(a,a'),\\
&&\be_a^{\scriptscriptstyle\times\!}(x,y)= |s- \left<a,a'|a\right>| - |t- \left<a,a'|a\right>|\\
&& \qquad = (s- \left<a'|x_0\right>_a) - (t-\left<a'|x_0\right>_a)= s-t.
\end{eqnarray*}
If $x=a$ and $y=a'$, then $s=-\left<a'|x_0\right>_a$ and $t=\left<a|x_0\right>_{a'}$,
so 
\begin{align*}
&\be_u^{\scriptscriptstyle\times\!}(a,a')\le |s - t|=|-\left<a'|x_0\right>_a- \left<a|x_0\right>_{a'}|=
  d(a,a'),\\
&\be_a^{\scriptscriptstyle\times\!}(a,a')= s-t=-\left<a'|x_0\right>_a- \left<a|x_0\right>_{a'}= -d(a,a').\tag*{\qed}
\end{align*}

\subsection{The pseudometric~$\dcross$ in~$\sx$}\label{ss_d^times}
For $x,y\in \sx$ define 
\begin{equation}\label{def-dcross}
\dcross(x,y):=
\sup_{u\in X} |\be_u^{\scriptscriptstyle\times\!}(x,y)|.
\end{equation}
\begin{ttt}\label{t_d^times}$\phantom{99}$
\begin{itemize}
\item [\rm(a)] The function $\dcross$ above is a well-defined $\mathrm{Isom}(X)$--invariant pseudometric
on $\sx$ independent of $x_0$.
\item [\rm(b)] The inclusion of each line 
  $([\![a,b]\!],|\cdot|)\hookrightarrow (\sx,\dcross)$, $a,b\in X$,
  is an isometric embedding.
\item [\rm(c)] The canonical embedding $(X, d)\hookrightarrow (\sx,\dcross)$
  is an isometric embedding.
\end{itemize}
\end{ttt}
\begin{proof}
{\bf (a)}\qua
Let $x=[\![a,a';s]\!]$, $y=[\![b,b';t]\!]$ for appropriate $s$ and $t$.
By the cocycle condition and Lemma~\ref{l_beuxy},
$$|\be_u^{\scriptscriptstyle\times\!}(x,y)|\le
|\be_u^{\scriptscriptstyle\times\!}(x,a)|+ |\be_u^{\scriptscriptstyle\times\!}(a,b)|
+ |\be_u^{\scriptscriptstyle\times\!}(b,y)|\le  d(a,a')+ d(a,b)+ d(b,b')$$
for all $u\in X$, so $\dcross(x,y)\in [0,\infty)$.

The product $\left<a|a'\right>_u$ is independent of~$x_0$
because it is expressed in terms of the metric~$d$,
the projection $[\![a,a'|u]\!]$ is independent by Lemma~\ref{l_projection-other-def}.
The quantity\break $|s-\left<a,a'|u\right>|$ is the distance between $x$ and 
$[\![a,a'|u]\!]$ in $[\![a,a']\!]$, so it is independent of $x_0$, and similarly for all the terms
in the definition of $\becross$, so this shows the independence of $\becross$.
The same argument provides $\Isom(X)$--invariance of $\dcross$.

Pick any triple $x,y,z\in \sx$ and any $\varepsilon >0$. 
By the definition of~$\dcross$, there is $u\in X$ such that
$\dcross(x,z)-\varepsilon\le |\be_u^{\scriptscriptstyle\times\!}(x,z)|$,
therefore
$$\dcross(x,z)\le |\be_u^{\scriptscriptstyle\times\!}(x,z)|+\varepsilon
\le |\be_u^{\scriptscriptstyle\times\!}(x,y)|+
  |\be_u^{\scriptscriptstyle\times\!}(y,z)|+\varepsilon
\le \dcross(x,y)+ \dcross(y,z)+\varepsilon.$$
Since this holds for each $\varepsilon >0$, the triangle inequality for~$\dcross$ follows.
Since $\be_u^{\scriptscriptstyle\times\!}(x,x)=0$ for all $u\in X$,
then $\dcross(x,x)=0$, so $\dcross$ is a pseudometric.

\noindent {\bf (b)}\qua If $x$ and $y$ lie on the same line $[\![a,a']\!]$ in $\sxxo$, then
$x=[\![a,a';s]\!]$ and $y=[\![b,b';t]\!]$ for some appropriate $s$ and $t$.
By Lemma~\ref{l_beuxy}(a), $\dcross(x,y)=|s-t|$.

\noindent {\bf (c)}\qua If $a,b\in X$, then by part (b),
$\dcross(a,b)$ equals the length of the interval $[\![a,b]\!]$. But $[\![a,b]\!]$
was chosen so that its length is $ d(a,b)$.
\end{proof}

\subsection{The metric $d_*=\sj d$ on $\sx$}\label{ss_metric-d*}
Define a function $\sj d\co (\sx)^2\to[0,\infty)$ by
\begin{equation}\label{def-dstar}
\sj d(x,y):=
  \int_{-\infty}^{\infty} \dcross(r^+x,r^+y)\,\frac{e^{-|r|}}{2}\,dr,
\qquad x,y\in \sx,
\end{equation}
where $r^+$ comes from the $\rr$--action on~$\sx$ described in~\ref{ss_r-action}.
For simplicity we will use the notation $d_*$ instead of $\sj d$.

\begin{ttt}\label{p_d*}
The function $d_*$ above is a well-defined $\mathrm{Isom}(X)$--invariant
metric on $\sx$ independent of the choice of~$x_0$.
\end{ttt}
\proof
By Theorem~\ref{t_d^times}(b), $\dcross$ induces the original topology on each line, 
therefore by the triangle inequality for $\dcross$, the restriction of $\dcross$
to each product $[\![a,a']\!]\times[\![b,b']\!]$ is continuous.
By Lemma~\ref{l_r-action}, for fixed $x$ and $y$,
$\dcross(r^+x,r^+y)$ is continuous in~$r$. 
First assume $x,y\in\starx$.
By Lemma~\ref{l_r-action} and Theorem~\ref{t_d^times}(b),
the $\rr$--orbit maps $\rr\to(\sx,\dcross)$ are non-expanding, 
hence
\begin{equation}\label{dtimesrxry}
0\le \dcross(r^+x,r^+y)\le \dcross(r^+x,x)+ \dcross(x,y)+  \dcross(y,r^+y)
 \le  \dcross(x,y)+2|r|
\end{equation}
for all $r\in \rr$.
If $x\in X$, then $\dcross(r^+x,x)= \dcross(x,x)=0$, and the same
inequality as above holds, and similarly for $y\in X$.
This inequality implies that $d_*(x,y)$ is a well-defined number in $[0,\infty)$ for all $x,y\in \sx$.
The triangle inequality, the $\mathrm{Isom}(X)$--invariance and independence of~$x_0$ follow from those
of~$\dcross$. Also, $d_*(x,x)=0$.

It remains to show that $x\not=y$ implies $d_*(x,y)>0$.
Pick any $x\in [\![a,a']\!]$ and $y\in[\![b,b']\!]$ with $x\not= y$. 
If $a\not= b$, then, by Theorem~\ref{t_d^times}(c), $\dcross(a,b)>0$. 
Since $r^+x\to a$ and $r^+y\to b$ as $r\to -\infty$, 
there exists $r_0\in\rr$ such that for all $r\le r_0$,
$$\dcross(r^+x,a)\le d(a,b)/3\qquad\mbox{and}\qquad \dcross(r^+y,b)\le d(a,b)/3,$$ 
hence $\dcross(r^+x,r^+y)\ge \dcross(a,b)/3$ and 
$$d_*(x,y)\ge
\int_{-\infty}^{r_0}  \dcross(r^+x,r^+y) \frac{e^{-|r|}}{2}\, dr\ge
\int_{-\infty}^{r_0}  \dcross(a,b) \frac{e^{-|r|}}{6}\, dr>0.$$
The case $a'\not=b'$ is similar, so now we can assume that $a=b$ and $a'=b'$,
i.e.\ $x$ and $y$ lie on the same line $[\![a,a']\!]$, and $x\not= y$.
For all $r\in\rr$, 
$r^+x$ and $x^+y$ lie on the same line $[\![a,a']\!]$ and $r^+x\not=r^+y$;
then by Theorem~\ref{t_d^times}(b), $\dcross(r^+x,r^+y)>0$.
Since $\dcross(r^+x,r^+y)$ is continuous in $r$,
$$d_*(x,y)= 
\int_{-\infty}^{\infty} \dcross(r^+x,r^+y) \frac{e^{-|r|}}{2}\, dr>0.\eqno{\qed}$$

\noindent {\bf Remark}\qua (\ref{def-dstar}) resembles the formula used by
Gromov in~\cite[8.3.B]{Gr2}. He starts with the set of all bi-infinite geodesics
in a hyperbolic metric space $X$, i.e.\ each geodesic comes equipped with
an embedding into $X$. Given two points $x$ and $y$ lying on two bi-infinite geodesics, one can view them as
lying in~$X$, measure the distance between them and apply~(\ref{def-dstar}) to define a metric on
the disjoint union of all geodesics.
In general there are many geodesics connecting points $a,b\in\p X$, and the metric is used
to identify all of them into one by a quasiisometric homeomorphism (not necessarily an isometry).
In that construction, $\rr$ does not necessarily act on lines by isometries.

In the construction of this paper we start with an arbitrary metric space~$X$.
Geodesics are of finite length and are abstractly assigned to each pair of points;
no embedding into the space is given. 
(There might be no embedding at all, the space may be even discrete!)
This is why it was important
to construct $\dcross$ first, and this was done using $\becross$ and
the double difference in~$X$.
The advantage of this formal approach is that for each pair $a,b\in X$ there 
is a canonically associated line in $\sx$ that depends continuously 
on $a$ and $b$. Moreover, we will see that when~$X$ is a hyperbolic complex,
this construction extends continuously to the ideal boundary of~$X$ so that $\rr$ acts
by isometries on each ideal line. We will also define and use the $\rr$--action on both finite and
infinite lines.

\begin{ttt}\label{bi-lipschitz}
For each $r\in\rr$, the map $r^+\co (\sx,d_*)\to(\sx,d_*)$ is a bi-Lipschitz homeomorphism
with constant $e^{|r|}$.
\end{ttt}
\begin{proof}
For all $v,r\in\rr$,
$e^{-|v-r|}\le e^{|r|-|v|}= e^{|r|} e^{-|v|}$,
hence using the substitution $v=u+r$,
\begin{eqnarray*}
&& d_*(r^+x,r^+y)=
\int_{-\infty}^{\infty} \dcross(u^+r^+x,u^+r^+y)\frac{e^{-|u|}}{2}\,du\\
&& =\int_{-\infty}^{\infty} 
\dcross\big((u+r)^+x,(u+r)^+y\big)\frac{e^{-|u|}}{2}\,du
= \int_{-\infty}^{\infty} \dcross(v^+x,v^+y)\frac{e^{-|v-r|}}{2}\,dv\\
&& \le e^{|r|}\int_{-\infty}^{\infty} \dcross(v^+x,v^+y)\frac{e^{-|v|}}{2}\,dv
= e^{|r|} d_*(x,y),
\end{eqnarray*}
and similarly for the inverse map $(-r)^+$.
\end{proof}

\subsection{Relation between $\dcross$ and $d_*$}
Define a function $\varphi=\varphi_X\co \sx\to\sx$ by
\begin{equation}\label{d_varphi}
\varphi(x):=
\int_{-\infty}^\infty r^+ x\, \frac{e^{-|r|}}{2}\, dr.
\end{equation}
Here $x\in[\![a,a']\!]$ and we view $[\![a,a']\!]$ as a subinterval $[\al,\al']\se\rr$ as
in~\ref{ss_parametrizations} to make sense of the integral.

\begin{ttt}\label{t_surjection}
Let $X$ be any metric space and define $\dcross$, $d_*$ and $\varphi$ as above.
\begin{itemize}
\item [\rm(a)] $\varphi$ is a well-defined canonical surjection
$(\sx,d_*)\twoheadrightarrow (\sx,\dcross)$
whose restriction to each line $([\![a,a']\!],d_*)$ is an isometry onto
$([\![a,a']\!],\dcross)$. In particular,  each line $[\![a,a']\!]$ in $\sx$
can be parameterized to become a $d_*$--geodesic from $a$ to $a'$.
\item [\rm(b)] The restriction of $\varphi$ to~$X$ is the identity
map~$(X,d_*)\to (X,\dcross)$, and it is an isometry. 
\item [\rm(c)] $d_*$, $\dcross$ and $d$ coincide on~$X$,
i.e.\ the canonical embeddings $(X,d)\hookrightarrow (\sx,\dcross)$ and
$(X,d)\hookrightarrow (\sx,d_*)$ are isometric.
\end{itemize}
\end{ttt}
\begin{proof}
\noindent {\bf (a)}\qua
Each $\varphi(x)$ is well-defined as a real number because $[\![a,a';\cdot]\!]'$
is non-expanding.

If $x\in [\![a,a']\!]\cap X$, then $\rr$ fixes $x$ and
$$\varphi(x)=\int_{-\infty}^\infty x\, \frac{e^{-|r|}}{2}\, dr= x,$$
so $\varphi(x)$ is well-defined and equals $x$, regardless of the choice of $[\![a,a']\!]$.
Therefore $\varphi$ is identity on~$X$.

Now we assume $x\in [\![a,a']\!]\setminus X$, i.e.\ $x\in ]\!]a,a'[\![$ and $a\not=a'$.
Let the function $[\![a,a';\cdot]\!]''\co \rr\to [\![a,a']\!]$ be
the smooth-out of $[\![a,a';\cdot]\!]'$, i.e.\
\begin{equation}\label{e_[[]]''}
[\![a,a';s]\!]'':= \int_{-\infty}^\infty [\![a,a';r+s]\!]' \, \frac{e^{-|r|}}{2}\, dr.
\end{equation}
Since $a\not= a'$, then $[\![a,a';\cdot]\!]'\co  \rb\to [\![a,a']\!]$ is a homeomorphism,
and by definition
$\varphi|_{[\![a,a']\!]}$ equals the composition
$[\![a,a';\cdot]\!]''\circ ([\![a,a';\cdot]\!]')\inv$.
By Lemma~\ref{l_theta} applied to the function $[\![a,a';\cdot]\!]'$,
$[\![a,a';\cdot]\!]''$ is a homeomorphism from $\rb$ onto $[\![a,a']\!]$,
therefore $\varphi$ maps $[\![a,a']\!]$ homeomorphically onto itself.

Pick any $x,y\in[\![a,a']\!]$ with $x\le y$. By Lemma~\ref{l_r-action}, 
$r^+x\le r^+y$ for all $r\in\rr$, hence $\varphi(x)\le \varphi(y)$.
By Theorem~\ref{t_d^times}(b),
\begin{eqnarray*}
&& \dcross\big(\varphi(x), \varphi(y) \big)= 
\varphi(y)- \varphi(x)= 
\int_{-\infty}^\infty (r^+ y- r^+x)\, \frac{e^{-|r|}}{2}\, dr\\
&& = \int_{-\infty}^\infty \dcross(r^+ y, r^+x)\, \frac{e^{-|r|}}{2}\, dr=
  d_*(x,y),
\end{eqnarray*}
hence $\varphi\co  ([\![a,a']\!],d_*)\to ([\![a,a']\!],\dcross)$ is an isometry.

\noindent {\bf (b)}\qua
In particular, 
$\dcross(a,a')=\dcross(\varphi(a),\varphi(a'))= d_*(a,a')\quad\mbox{for all}\quad a,a'\in X$,
so the restriction $\varphi\co  (X,d_*)\to (X,\dcross)$ is an isometry.

{\bf (c)}\qua This follows from (b) and Theorem~\ref{t_d^times}(c).
\end{proof}

\begin{lemma}\label{l_aa'-dcross-aa'-dstar}
For all $a,a'\in X$, 
\begin{itemize}
\item [\rm(a)] the identity map $([\![a,a']\!],\dcross)\to ([\![a,a']\!],d_*)$ 
  is a homeomorphism and
\item [\rm(b)]  if $a\not= a'$, then $[\![a,a';\cdot]\!]'\co \rb\to ([\![a,a']\!],d_*)$
  is a homeomorphism.
\end{itemize}
\end{lemma}
\proof
{\bf (a)}\qua If $a=a'$, then the statement is obvious.
Now assume $a\not= a'$.
Consider the composition
$$\rb\overset{[\![a,a';\cdot]\!]'}{\to} ([\![a,a']\!],\dcross)
\overset{id}{\to} ([\![a,a']\!],d_*)
\overset{\varphi}{\to}([\![a,a']\!],\dcross),$$
where $\varphi$ is the isometry from Theorem~\ref{t_surjection}(a).
The first map is a homeomorphism, and one checks that the composition
is given by the formula
$$s\mapsto \int [\![a,a';r+s]\!]'\frac{e^{-|r|}}{2}\,dr,$$
hence it is the map $[\![a,a';\cdot]\!]''$ in~(\ref{e_[[]]''}),
which is, again, a homeomorphism by Lemma~\ref{l_theta} applied
to the function $[\![a,a';\cdot]\!]'$.
This implies that $([\![a,a']\!],\dcross)\overset{id}{\to} ([\![a,a']\!],d_*)$ 
is a homeomorphism.

\noindent {\bf (b)}\qua
The map $[\![a,a';\cdot]\!]'\co \rb\to ([\![a,a']\!],\dcross)$
is the composition of homeomorphisms
$$\rb\overset{[\![a,a';\cdot]\!]'}{\to} ([\![a,a']\!],\dcross)
\overset{id}{\to} ([\![a,a']\!],d_*).\eqno{\qed}$$

\begin{lemma}\label{l_d*dcross}
For all $x,y\in\sx$,  $|d_*(x,y)-\dcross(x,y)|\le 2.$
In particular, $d_*$ and $\dcross$ are $^+$equivalent.
\end{lemma}
\proof
By~(\ref{dtimesrxry}),
\begin{eqnarray*}
& & d_*(x,y)=
\int_{-\infty}^{\infty} \dcross(r^+x,r^+y)\frac{e^{-|r|}}{2}\,dr\\
& & \qquad\quad\le\int_{-\infty}^{\infty} \big(\dcross(x,y)+2|r|\big) \frac{e^{-|r|}}{2}\,dr
= \dcross(x,y)+ 2.
\end{eqnarray*}
Similarly,
\begin{align*}
& d_*(x,y)=
\int_{-\infty}^{\infty} \dcross(r^+x,r^+y)\frac{e^{-|r|}}{2}\,dr\\
& \qquad\quad\ge\int_{-\infty}^{\infty} \big(\dcross(x,y)-2|r|\big) \frac{e^{-|r|}}{2}\,dr
= \dcross(x,y)- 2.\tag*{\qed}
\end{align*}

\subsection{The functor $\sj$ and embeddings into geodesic spaces}\label{ss_functor}
$\sj$ is a functor on the category of topological spaces: 
to every topological space $X$ it associates the topological space $\sx$. 
$\sj$~is also a functor on the category of metric spaces:
to every metric space $(X,d)$ it associates the space $\sx$ with the metric
$\sj d$ as in~\ref{ss_metric-d*}. 
We are intentionally vague about the choice of morphisms here -- it is an interesting 
educational question how $\sj$ behaves under continuous, Lipschitz and other maps,
but this will not be addressed in this article.
At the very least, since the construction is canonical,
$\sj$ is functorial with respect to isometries.

As an illustration for the use of functor $\sj$ we prove the following fact.
\begin{ppp}
Every metric space $(X,d_X)$ isometrically embeds into a geodesic metric space~$(Y,d_Y)$.
Moreover, $Y$ can be chosen so that each $g\in\mathrm{Isom}(X)$ extends
to an isometry $g'$ of~$Y$, and the map $\mathrm{Isom}(X)\to \mathrm{Isom}(Y)$, $g\mapsto g'$,
is a group monomorphism.
\end{ppp}
\begin{proof}
Define $\sj^0 X:= X$ and inductively $\sj^{i} X:=\sj (\sj^{i-1} X)$. 
Then $\sj^{i} X$ is a metric space with the metric
$d_i:=\sj^{i} d_X$ defined inductively from the metric on $X$.
By Theorem~\ref{t_surjection}(c) there are canonical isometric embeddings
$\sj^{i-1} X\hookrightarrow\sj^{i} X$,
therefore the union 
$$Y:= \sj^{\infty}(X):= \bigcup_{i=0}^{\infty} \sj^{i} X$$
can be given a metric $d_Y$ which restricts to $d_i$ on each $\sj^{i} X$.
By Theorem~\ref{t_surjection}(b), each pair of points $a,b\in \sj^{i-1} X$ is
connected by a geodesic in $\sj^{i}X$, so $Y$ is geodesic. 
By Theorem~\ref{p_d*}, each isometry $g$ of $X$ induces an isometry
of $\sj X$. This gives a homomorphism $\mathrm{Isom}(X)\to \mathrm{Isom}(\sj X)$
which is clearly injective. Inductively, $g$ extends to an isometry of $\sj^{i} X$,
giving a monomorphism $\mathrm{Isom}(X)\hookrightarrow \mathrm{Isom}(\sj^{i} X)$ for each $i$, and therefore 
to an isometry of~$Y$ giving a monomorphism $\mathrm{Isom}(X)\hookrightarrow \mathrm{Isom}(Y)$. 
\end{proof}

\section{The metric $d_*$ and the topology of $\starx$}\label{s_metric-and-topology}
The goal of this section is to prove the following.

\begin{ppp}\label{p_Xhomeomorphic}
Let $(X,d)$ be any metric space. The metric $d_*$ from~\ref{ss_metric-d*}
induces the original topology on the open
symmetric join~$\starxxo$ described in~\ref{ss_two-models-starx-starxxo}.
Equivalently,
the map from~\ref{ss_two-models-starx-starxxo} viewed as
$$]\!]\cdot,\cdot\,;\cdot[\!['\co 
  (X^2\setminus\Delta)\times\rr\to (\starxxo,d_*)$$
is a homeomorphism,
where $\Delta$ is the diagonal of $X^2$.
\end{ppp}

The proof requires some technical arguments, the reader might want 
to skip this section at first reading.

\begin{lemma}\label{l_exy}
For all $x,y\in[0,\infty)$, \ $|e^{-y}-e^{-x}|\le |y-x|$.
\end{lemma}
\begin{proof}
One easily checks that $1-e^z\le -z$ for all $z$.
We can assume $x\le y$. Then
$\hfill |e^{-y}-e^{-x}|= e^{-x}|1-e^{x-y}|\le |1-e^{x-y}|= 1-e^{x-y}\le y-x= |y-x|.$
\end{proof}

\begin{lemma} \label{l_dhat-dcross-dstar}
Let $(X,d)$ be any metric space. Then for all $a,a',b,b'\in X$ and $s,t\in\rb$,
\begin{itemize}
\item [\rm(a)] $\quad \dcross ([\![a,a';s]\!],[\![b,b';t]\!])\le 
  |t-s|+ 2 d(a,b)+ 2d(a',b')$,
\item [\rm(b)] $\quad \dcross ([\![a,a';s]\!]',[\![b,b';t]\!]')\le |t-s|+ 4d(a,b)+ 4d(a',b')$,
\item [\rm(c)] $\quad d_* ([\![a,a';s]\!]',[\![b,b';t]\!]')\le |t-s|+ 4d(a,b)+ 4d(a',b')$.
\end{itemize}
\end{lemma}
\proof
\noindent {\bf (a)}\qua 
For any $u\in X$, by the definition of $\becross$ and triangle inequality,
\begin{align*}
&|\becross_u([\![a,a';s]\!],[\![b,b';t]\!])|=
  \big|\left<a|a'\right>_u+ \left|s-\left<a,a'|u\right>\right|-
        \left<b|b'\right>_u- \left|t-\left<b,b'|u\right>\right|\big|\\
&\le \big|\left|s-\left<a,a'|u\right>\right|-\left|t-\left<b,b'|u\right>\right|\big|+
\left|\left<a|a'\right>_u-\left<b|b'\right>_u\right|\\
&\le \big|s-\left<a,a'|u\right>-t+\left<b,b'|u\right>\big|+
\left|\left<a|a'\right>_u-\left<b|b'\right>_u\right| \\
&\le |t-s|+ \left|\left<a,a'|u\right>- \left<b,b'|u\right>\right|+
  \left|\left<a|a'\right>_u-\left<b|b'\right>_u\right| \\
&\le |t-s|+ \left|\left<a,b|u\right>\right|+
   \left|\left<a',b'|u\right>\right|+
   \left|\left<a|a'\right>_u-\left<b|a'\right>_u\right|+
   \left|\left<b|a'\right>_u-\left<b|b'\right>_u\right|\\
&\le |t-s|+ 2 d(a,b)+ 2 d(a',b').
\end{align*}
Then by the definition of~$\dcross$,
$$ \dcross([\![a,a';s]\!],[\![b,b';t]\!])=
 \sup_{u\in X} |\becross_u([\![a,a';s]\!],[\![b,b';t]\!])| \le |t-s|+ 2d(a,b)+ 2d(a',b').$$
\noindent {\bf (b)}\qua
Recall from~\ref{ss_parametrizations} that $[\![a,a']\!]$ is a copy of the interval
$[\al,\al']\se\rr$,
where $\al:= -\left<a'|x_0\right>_a$ and $\al':= \left<a|x_0\right>_{a'}$.
Similarly, $[\![b,b']\!]$ is a copy of $[\be,\be']\se\rr$ where
$\be:= -\left<b'|x_0\right>_b$ and $\be':= \left<b|x_0\right>_{b'}$.
By the definition of $\left<\cdot|\cdot\right>_\cdot$ and triangle inequality, 
\begin{eqnarray}\label{be-al}
&& |\be-\al|= |\left<a'|x_0\right>_a- \left<b'|x_0\right>_b|
   \le  \big| d(a,a')-  d(b,b')\big|/2 \\
\nonumber && + \big|d(a,x_0)-  d(b,x_0)\big|/2+
   \big| d(b',x_0)-  d(a',x_0)\big|/2 \le  d(a,b)+  d(a',b'),
\end{eqnarray}
and similarly $|\be'-\al'|\le  d(a,b)+  d(a',b')$.
It follows from the definition of $\theta$ that
\begin{equation}\label{theta-be-al}
\big|\theta[\be,\be';t]- \theta[\al,\al';t]\big|\le \max\{|\be-\al|, |\be'-\al'|\}\le
 d(a,b)+  d(a',b').
\end{equation}
Using Lemma~\ref{theta'theta} we denote
\begin{eqnarray*}
&& A:= \theta'[\al,\al';t]= \theta[\al,\al';t]+ (e^{-|t-\al|}- e^{-|t-\al'|})/2\quad\mbox{and}\\
&& B:= \theta'[\be,\be';t]= \theta[\be,\be';t]+ (e^{-|t-\be|}- e^{-|t-\be'|})/2,
\end{eqnarray*}
then by~(\ref{theta-be-al}), Lemma~\ref{l_exy} and~(\ref{be-al}),
\begin{eqnarray*}
&&|B-A|\\
&& \le\big|\theta[\be,\be';t]- \theta[\al,\al';t]\big|+
 \big|e^{-|t-\be|}- e^{-|t-\al|}\big|/2+ \big|e^{-|t-\al'|}- e^{-|t-\be'|}\big|/2\\
&& \le  d(a,b)+  d(a',b') + \big||t-\be|- |t-\al|\big|/2+ \big||t-\al'|- |t-\be'|\big|/2\\
&& \le  d(a,b)+  d(a',b')+ |\be-\al|/2+ |\be'-\al'|/2\le
  2 d(a,b)+ 2 d(a',b').
\end{eqnarray*}
By~(\ref{abt'}), part (a) and the above inequality,
\begin{eqnarray*}
&& \dcross\big([\![a,a';t]\!]', [\![b,b';t]\!]'\big)=
   \dcross\big([\![a,a';\theta'[\al,\al';t]]\!], [\![b,b';\theta'[\be,\be';t]]\!]\big)\\
&& = \dcross\big([\![a,a';A]\!],[\![b,b';B]\!]\big)\le
   |B-A|+ 2 d(a,b)+ 2 d(a',b')\\
&& \le 4 d(a,b)+ 4 d(a',b').
\end{eqnarray*}
Since the map $[\![a,a';\cdot]\!]'\co \rr\to ([\![a,a']\!],\dcross)$ is non-expanding,
\begin{eqnarray*}
&& \dcross\big([\![a,a';s]\!]', [\![b,b';t]\!]'\big)\le
   \dcross\big([\![a,a';s]\!]', [\![a,a';t]\!]'\big)+
   \dcross\big([\![a,a';t]\!]', [\![b,b';t]\!]'\big)\\
&& \le |t-s|+ 4 d(a,b)+ 4 d(a',b').
\end{eqnarray*}
{\bf(c)}\qua This follows from (b) and the definition of~$d_*$:
\begin{align*}
& d_* ([\![a,a';s]\!]',[\![b,b';t]\!]')=
\int_{-\infty}^{\infty} \dcross \big([\![a,a';r+s]\!]',[\![b,b';r+t]\!]'\big)\frac{e^{-|r|}}{2}\, dr\\
&\le\int_{-\infty}^{\infty}\Big(|(r+t)-(r+s)|+ 4d(a,b)+ 4d(a',b')\Big)
  \frac{e^{-|r|}}{2}\, dr\\
&= |t-s|+ 4d(a,b)+ 4d(a',b').\tag*{\qed}
\end{align*}

\begin{lemma}\label{l_omega}
The function $\omega\co [0,\infty)\to [0,\infty)$ defined
by $\omega(\tau):= \tau+2e^{-\tau/2}-2$ is a homeomorphism.
The obvious extension $\omega\co [0,\infty]\to[0,\infty]$ is also a homeomorphism.
\end{lemma}
\begin{proof}
This follows from the facts that $\omega(0)=0$, 
$\ds\frac{\p\omega}{\p \tau}(\tau)> 0$ for $\tau>0$, and
$\omega(\tau)\to\infty$ as $\tau\to\infty$.
\end{proof}

\begin{lemma}\label{l_dcross-dstar}
For all $x,y\in \sx$,  $\dcross(x,y)\le \omega\inv(d_*(x,y))$,
where $\omega$ is from Lemma~\ref{l_omega}.

Moreover, for all $x,y\in \sx$ and $r\in\rr$,
$\dcross(r^+x,r^+y)\le  \omega\inv\big(e^{|r|}d_*(x,y)\big)$.
\end{lemma}
\begin{proof}
Denote $\tau:=\dcross(x,y)$.
Since the map $[\![a,b;\cdot]\!]'$ is non-expanding, we have
$$\dcross(r^+x,r^+y)\ge \dcross(x,y)-\dcross(x,r^+x)-\dcross(y,r^+y)\ge
  \dcross(x,y)- 2|r|=\tau -2|r|$$
for all $r\in\rr$, then
\begin{eqnarray}\label{dstar-omega-dcross}
&& d_*(x,y)=\int_{-\infty}^\infty \dcross(r^+x,r^+y)
    \frac{e^{-|r|}}{2}\, dr\ge
  \int_{-\tau/2}^{\tau/2} (\tau- 2|r|)\frac{e^{-|r|}}{2}\, dr\\
\nonumber
&& = \tau+2e^{-\tau/2}-2=\omega(\tau)= \omega (\dcross(x,y)).
\end{eqnarray}
Since $\omega$ is increasing, $\dcross(x,y)\le \omega\inv(d_*(x,y))$
for any $x$ and $y$ in $\sx$. Applying 
the inequality~(\ref{dstar-omega-dcross}) to $r^+x$ and $r^+y$ and using 
Theorem~\ref{bi-lipschitz},
$$\omega\big(\dcross(r^+x,r^+y)\big)\le d_*(r^+x,r^+y)\le
e^{|r|}d_*(x,y),$$
therefore  $\dcross(r^+x,r^+y)\le \omega\inv\big(e^{|r|}d_*(x,y)\big)$.
\end{proof}

\begin{lemma}\label{l_aixor-a'ix}
Let $a\in X$ and $x_i=[\![a_i,a'_i;s_i]\!]'$ be a sequence in $\sx$
such that $d_*(x_i,a)\to 0 $ as $i\to\infty$.
Then $d(a_i,a)\to 0$ or $d(a'_i,a)\to 0$ as $i\to\infty$.
\end{lemma}
\begin{proof}
Suppose not, then after taking a subsequence there exists $\varepsilon> 0$
such that $\dcross(a_i,a)=d(a_i,a)>\varepsilon$ and 
$\dcross(a'_i,a)=d(a'_i,a)>\varepsilon$ 
for all~$i$. Since $d_*(x_i,a)\to 0$, then by Lemma~\ref{l_dcross-dstar},
$\dcross(x_i,a)\to 0$, so extracting a subsequence again
we can assume that 
$$\dcross(x_i,a)\le \varepsilon/4\qquad\mbox{for all }i.$$
For each $i$ there is a point $y_i=[\![a_i,a'_i;t_i]\!]'$ 
such that 
$$t_i>s_i\qquad\mbox{and}\qquad\dcross(x_i,y_i)=\varepsilon/2.$$
In particular,
\begin{equation}\label{dcross-xyixiyi}
\dcross(a,y_i)\ge \dcross(x_i,y_i)-\dcross(x_i,a)\ge
\varepsilon/2-\varepsilon/4=\varepsilon/4.
\end{equation}
Let $z\in [\![a_i,a'_i]\!]$ represent an arbitrary point that lies between
$x_i$ and $y_i$, then
\begin{eqnarray*}
&&\dcross(a,z)\le \dcross(a,x_i)+\dcross(x_i,z)\le
  \varepsilon/4+\varepsilon/2=3\varepsilon/4, \quad\mbox{hence}\\
&&\dcross(z,a_i)>\varepsilon/4\qquad\mbox{and}\qquad 
\dcross(z,a'_i)>\varepsilon/4.
\end{eqnarray*}
Take $\epsilon>0$ sufficiently small so that
$$\epsilon <\varepsilon/4\qquad\mbox{and}\qquad
  \epsilon< (1-e^{-\varepsilon})/2.$$
By the above,
\begin{equation}\label{dcross-zai-epsilon}
\dcross(z,a_i)>\varepsilon/4>\epsilon\qquad\mbox{and}\qquad 
\dcross(z,a'_i)>\varepsilon/4>\epsilon\qquad\mbox{for all }i.
\end{equation}
Let $\al_i\le \al'_i$ be the end-point coordinates of $[\![a_i,a'_i]\!]$.
Since 
\begin{eqnarray*}
&&\al'_i-\al_i=
\dcross(a_i,a'_i)= \dcross(a_i,x_i)+\dcross(x_i,a'_i)\\
&&\ge \big(\dcross(a_i,a)-\varepsilon/4\big)+
\big(\dcross(a'_i,a)-\varepsilon/4\big)
\ge2(\varepsilon-\varepsilon/4)>\varepsilon,
\end{eqnarray*}
we have
\begin{equation}\label{epsilon1-ealal'}
0< \epsilon< (1-e^{-\varepsilon})/2\le (1-e^{\al_i-\al'_i})/2\qquad\mbox{for all }i.
\end{equation}
Equations (\ref{dcross-zai-epsilon}) and (\ref{epsilon1-ealal'}) show that
Lemma~\ref{l_lower-bound} applies to the map 
$[\![a_i,a'_i;\cdot]\!]'\co \rr\to ([\![a_i,a'_i]\!],\dcross)$, therefore
the derivative of this map in the interval $[s_i,t_i]$  is at least~$\epsilon$. Then 
$$t_i-s_i\le \frac{\dcross(x_i,y_i)}{\epsilon}=\frac{\varepsilon}{2\epsilon}.$$
Denote $R:=\varepsilon/(2\epsilon)$ and $r_i:=t_i-s_i\in [0,R]$,
so we have $y_i={r_i}^+x_i$.
By Lemma~\ref{l_dcross-dstar},
\begin{eqnarray*}
&& \dcross(a,y_i)=
  \dcross(a,{r_i}^+x_i)= \dcross({r_i}^+a,{r_i}^+x_i)\\
&& \qquad\le  \omega\inv\big(e^{r_i}d_*(a,x_i)\big)
 \le \omega\inv\big(e^{R}d_*(a,x_i)\big)\underset{i\to\infty}{\to} 0.
\end{eqnarray*}
This contradicts~(\ref{dcross-xyixiyi}).
\end{proof}

\begin{proof}[Proof of Proposition~\ref{p_Xhomeomorphic}]
Recall from~\ref{ss_two-models-starx-starxxo}
that the topology on 
$\starxxo\se \sxxo$ is induced by the bijection
$]\!]\cdot,\cdot\,;\cdot[\!['\co  (X^2\setminus\Delta)\times\rr\twoheadrightarrow \starxxo$.
Since the set $(X^2\setminus\Delta)\times\rr$ is open in $X^2\times \rb$,
the topology of $(X^2\setminus\Delta)\times\rr$
is locally the product of the topologies on $X$, $X$, and $\rr$. Therefore a  typical neighborhood of
$x= ]\!]a,a';s[\!['$ in $\starxxo$
is of the form $N(x,\varepsilon):=]\!]B_\varepsilon,B'_\varepsilon;I_\varepsilon[\!['$,
where $B_\varepsilon$ and $B'_\varepsilon$ are {\em disjoint} open balls in~$(X, d)$ 
of radius $\varepsilon$ centered
at~$a$ and~$a'$, respectively, and $I_\varepsilon=(s-\varepsilon,s+\varepsilon)\se\rr$.

Our goal is to show that the identity map $\starxxo\to(\starxxo,d_*)$ 
is a homeomorphism.
Let  $y= ]\!]b,b';t[\!['$ represent an arbitrary point in $N(x,\varepsilon)$,
then $|t-s|<\varepsilon$, $d(a,b)<\varepsilon$, $ d(a',b')<\varepsilon$.

Given any $\epsilon>0$,
we let $\varepsilon:=\epsilon/9$, then Lemma~\ref{l_dhat-dcross-dstar}(c) implies
$$d_*(x,y)=
d_*\big(]\!]a,a';s[\![', ]\!]b,b';t[\!['\big)
\le |t-s|+ 4 d(a,b)+4 d(a',b')<\varepsilon+4\varepsilon+4\varepsilon=\epsilon$$
i.e.\ $N(x,\varepsilon)\se B_{d_*}\!(x,\epsilon)$.
This shows that the identity map $\starxxo\to(\starxxo,d_*)$ is continuous at $x$.

Suppose that the inverse identity map $(\starxxo,d_*)\to\starxxo$ is not continuous at 
$x=]\!]a,a';s[\!['$, then there exist
a product neigbourhood $N(x,\varepsilon)=]\!]B_\varepsilon,B'_\varepsilon;I_\varepsilon[\!['$
of~$x$ and a sequence $\{x_i\}$ such that 
\begin{equation}\label{dstarxix}
d_*(x_i,x)\underset{i\to\infty}\to 0\qquad\mbox{and}\qquad
 x_i\in(\starx)\setminus N(x,\varepsilon).
\end{equation}
 We have $x_i=]\!]a_i,a'_i;s_i[\!['$ for some
$a_i,a'_i\in X$ and $s_i\in\rr$.

Since $d_*(x_i,x)\to 0$, using Lemma~\ref{l_dcross-dstar},
for any $j> 0$ there exists $i=i(j)$ such that 
\begin{equation}\label{e_dstar-jxijx}
d_*(j^+x_i,j^+x)\le \omega\inv\big(e^{|j|}d_*(x_i,x)\big)\le 1/j.
\end{equation}
Therefore the subsequence $\{x_j:=x_{i(j)}\}$ satisfies 
$d_*(j^+x_j,j^+x)\to 0$ as $j\to\infty$.
By Lemma~\ref{l_aa'-dcross-aa'-dstar}(b),
$d_*(j^+x,a')\to 0$ as $j\to\infty$, hence
$$d_*(j^+x_j,a')\le d_*(j^+x_j,j^+x)+ d_*(j^+x,a')
\underset{j\to\infty}{\to} 0.$$
Then by Lemma~\ref{l_aixor-a'ix} applied to the sequence $\{j^+x_j\}$
we have
\begin{equation}\label{aia'a'ia'}
a_j\to a'\qquad\mbox{or}\qquad
a'_j\to a'\qquad\mbox{in }X.
\end{equation}
By a similar argument using $(-j)^+$ in~(\ref{e_dstar-jxijx})
and passing to a subsequence we also deduce that
\begin{equation}\label{aiaa'ia}
a_j\to a\qquad\mbox{or}\qquad
a'_j\to a \qquad\mbox{in }X.
\end{equation}
Suppose that $a_j\to a'$ and $a'_j\to a$.
Since $\rb$ is compact, passing to a subsequence we can assume
that $s_j\to t$ for some $t\in\rb$.
Lemma~\ref{l_dhat-dcross-dstar}(c) implies that 
\begin{eqnarray}\label{dstar-aja'jsj}
\nonumber &&  d_*(x_j,[\![a',a;t]\!]') =
  d_*([\![a_j,a'_j;s_j]\!]',[\![a',a;t]\!]')\\
&& \qquad\le 
  |t-s_j|+ 4d(a_j,a')+ 4d(a'_j,a)\underset{j\to\infty}\to 0.
\end{eqnarray}
(\ref{dstar-aja'jsj}) and (\ref{dstarxix}) say that $x_j$ converges
in metric $d_*$ both to $x\in]\!]a,a'[\![$ and to $[\![a',a;t]\!]'$.
This is impossible since $]\!]a,a'[\![$ and $[\![a',a]\!]$
are disjoint.

The only possibility left is that $a_j\to a$ and $a'_j\to a'$.
After passing to a subsequence we can assume that 
$a_j\in B_\varepsilon$ and $a'_j\in B'_\varepsilon$ for all $j$.
Since $x_i\notin N(x,\varepsilon)$ by~(\ref{dstarxix}), 
we must have $|s_j-s|\ge \varepsilon$ for all $j$.
After passing to a subsequence we can assume that
$$ s_j\to \bar{s}\quad\mbox{for some \ } 
  \bar{s}\in[-\infty,s-\varepsilon]\cup[s+\varepsilon,\infty].$$
Then by Lemma~\ref{l_dhat-dcross-dstar}(c),
\begin{eqnarray*}
&& d_*(x_j,[\![a',a;\bar{s}]\!]') =
  d_*([\![a_j,a'_j;s_j]\!]',[\![a',a;\bar{s}]\!]')\\
&& \qquad \le 
  |\bar{s}-s_j|+ 4d(a_j,a)+ 4d(a'_j,a')\underset{j\to\infty}\to 0.
\end{eqnarray*}
Then $x_j$ converges in the metric $d_*$ both to $x=]\!]a,a';s[\!['$
and to $[\![a',a;\bar{s}]\!]'$, which is impossible since the two 
points are distinct.
\end{proof}

\section{Metric complexes and hyperbolic complexes}
\label{s_metric-complexes}
Suppose $X$ is a simplicial complex and $d$ is any metric on
its 0-skeleton~$X^{(0)}$. Let $P(X^{(0)})$ be the power set of $X^{(0)}$.
Each simplex in~$X$ is uniquely determined by its vertices,
so the simplicial structure on~$X$ can described combinatorially as
a collection~$\mathcal{U}\se P(X^{(0)})$
consisting of finite subsets of~$X^{(0)}$ which is subset-closed:
$U\in \mathcal{U}$ and $U'\se U$ imply $U'\in \mathcal{U}$.
Moreover, each $U\in \mathcal{U}$ can be viewed as the convex hull of its vertices, that is 
each point in the topological simplex $\sigma_U$ corresponding to $U$ is described 
uniquely as a linear combination
$$\sum_{x\in U} \al_x x\qquad \mbox{where}\qquad\al_x\in [0,1]\quad \mbox{and}\quad
\sum_{x\in U} \al_x= 1.$$

It is easy to check that the formula
\begin{equation}\label{metr-ext}
\tilde{d}\left(\sum_{x\in U} \al_x x, \sum_{y\in U} \be_y y\right):=
  \sum_{x\in U}  \sum_{y\in U} \al_x \be_y  d(x,y)
\end{equation}
defines a metric~$\tilde{d}$ on~$X$ whose restriction to each simplex $\sigma_U$
is homeomorphic (even linearly isomorphic) to the standard simplex of dimension~$\# U-1$.
Moreover, the inclusion
$(X^{(0)}, d)\hookrightarrow (X,\tilde{d})$ is an isometric embedding, i.e.\
$\tilde{d}$ is an extension of $ d$. We will omit $\tilde{\ }$ from the notation
of extended metric.

\subsection{The functor $\Psi$ and the canonical word metric $d_X$}
\label{ss_word-metric}
Let $X$ be a simplicial complex and $d$ be an arbitrary generalized metric on $X$,
$X^{(1)}$ or $X^{(0)}$. Denote $\Psi(d)$ the result of the following procedure:
restrict~$d$ to~$X^{(0)}$ and extend to all of $X$ by linearity formula~(\ref{metr-ext}).
$\Psi(d)$ is a generalized metric, and $\Psi(d)$ is a metric if and only if $d$ is.
If the simplices in $(X,d)$ have uniformly bounded diameters, then the inclusions
$(X^{(0)},d)\hookrightarrow (X,d)$ and $(X^{(0)},d)\hookrightarrow (X,\Psi(d))$
and the identity map $(X,d)\to (X,\Psi(d))$
are quasiisometries.

For any simplicial complex $X$ there is a canonical choice of generalized metric
$d_X$, obtained as follows: take the generalized path metric on $X^{(1)}$
defined by assigning length~1 to each edge, and apply functor $\Psi$.
$d_X$ is a metric if and only if $X$ is connected; we will call it the {\sf word metric}
on $X$, since this generalizes the word metrics on groups. In general
$d_X$ is {\em not} intrinsic.

\subsection{Metric complexes}\label{ss_metric-complexes}
We will say that a metric $d$ on~$X$ is 
{\sf induced from the 0-skeleton}\label{i_induced-from-0skel}
if $d=\Psi(d)$.
A {\sf metric complex\label{metric-complex}} will be a pair $(X,d)$ where~$X$ is a simplicial complex
and $d$ is a metric that is induced from the 0-skeleton.
Examples:
\begin{itemize}
\item [\rm(0)] Any metric space is a 0-dimensional metric complex.
\item [\rm(1)] Let $\G$ be the Cayley graph of any finitely generated group
with respect to a finite generating set.
Subdivide $\G$ if needed to make it a simplicial complex.
For any metric $d$ on $\G$ quasiisometric to the word metric,
$(\G,\Psi(d))$ is a 1-dimensional metric complex,
where $\Psi(d)$ is quasiisometric to $d$.
\item [\rm(2)] Let $M$ be any compact triangulated
smooth $n$--manifold with $\pi_1(M)=\Gamma$.
Its universal cover~$\tilde{M}$ has the distance function~$d$ coming from
the Riemannian structure, then $(\tilde{M},\Psi(d))$ is an $n$--dimensional metric complex.
The identity map $(\tilde{M},d)\to (\tilde{M},\Psi(d))$ is a $\Gamma$--invariant 
homeomorphism and quasiisometry.
\item [\rm(3)] Let $M$ be any compact triangulated $n$--manifold with $\pi_1(M)=\Gamma$.
Let~$d$ be any $\Gamma$--invariant metric on $\tilde{M}$ quasiisometric to
the word metric $d_X$ defined in~\ref{ss_word-metric}, for example the word metric itself.
Then $(\tilde{M},\Psi(d))$ is an $n$--dimensional metric complex homeomorphic
to $\tilde{M}$. The metrics $d_X$ and $\Psi(d)$ are quasiisometric.
\end{itemize}

For a metric complex $(X,d)$, $\Isom(X,d)$ will denote the group
of simplicial automorphisms of $X$ preserving the given metric $d$.
If $d$ is the word metric~$d_X$ on~$X$, then we will use the notation $\Isom(X)$
for $\Isom(X,d)$. Since $d_X$ is canonical, $\Isom(X)$ is just the group of
all simplicial automorphisms of $X$.

\subsection{Hyperbolic complexes}\label{ss_hyp-compl}
First let $\G$ be any connected graph.
We equip $\G$ with the {\sf word metric}~$d$\label{i_word-metric2}
 which is by definition 
the path metric induced by assigning length 1 to each edge.
A {\sf geodesic in $\G$}\label{i_geodesic2} is an isometric embedding
$\al \co  I\to (\G,d)$ of a closed interval $I\se \rb$. 
Often we will use the same notation for a geodesic and its image.
For all $a,b\in \bar\G$ we denote $\Geod(a,b)$\label{i_geod} the set of all geodesics in $\bar\G$
connecting $a$ to $b$, and fix one arbitrary choice of a geodesic $[a,b]\in \Geod(a,b)$.
When $a\in \G$, we assume more precisely that $[a,b]$ is the image
of the isometric embedding $[a,b;\cdot]\co  [0,d(a,b)]\to\G$ with
$[a,b;0]=a$ and $[a,b;d(a,b)]=b$; that is
for $t\in [0,d(a,b)]$, $[a,b;t]$ is the unique point in $[a,b]$ with
$d(a,[a,b;t])=t$.

Geodesic triangles have canonically defined {\sf inscribed triples}:\label{i_inscr-trip}
for all $a,b,c\in \G$ and any choice of geodesics $\al\in \Geod(b,c)$, $\be\in \Geod(c,a)$,
$\gamma\in \Geod(a,b)$ there are unique $\bar{a}\in\al$, $\bar{b}\in\be$,
$\bar{c}\in\gamma$ with 
$$d(a,\bar{b})=d(a,\bar{c}),\quad d(b,\bar{a})=d(b,\bar{c}),\quad
  d(c,\bar{a})=d(c,\bar{b}).$$
Equivalently, $\bar{a}$ is the unique point in~$\al$ satisfying
$d(b,\bar{a})- d(\bar{a},c)=d(b,a)- d(a,c)$, and similarly for $\bar{b}$ and $\bar{c}$.
It is convenient to think of~$\bar{a}$ as of a projection of~$a$ to~$\al$.

A graph $\G$ is {\sf hyperbolic\label{hyperbolic-graph}} if it is connected and there exists $\de\in[0,\infty)$
so that all geodesic triangles in $\G$ are $\de$--fine:
for all $a,b,c\in\G$ and the inscribed triple $\bar{a},\bar{b},\bar{c}$
as above we have
$$\mbox{if}\ x\in \be,\ y\in\gamma\ \mbox{and}\ 
d(a,x)=d(a,y)\le d(a,\bar{b})=d(a,\bar{c}),\ \ \mbox{then}\ 
d(x,y)\le\de.$$

Each simplicial complex $X$ has the canonical word metric $d=d_X$ as in~\ref{ss_word-metric}.
A {\sf hyperbolic complex\label{hyperbolic-complex}} will be a uniformly locally finite metric complex
whose 1-skeleton $(\G,d)$ is a hyperbolic graph. Note that it is a part
of the definition that $X$ is connected and uniformly locally finite.

Each hyperbolic graph has ideal boundary (see \cite{Gr2,GHV1991,BH} for definitions).
Since each hyperbolic complex $X$
is quasiisometric to its 1-skeleton $\G$, this defines the boundary
$\p X:=\p\G$ and the compactification $\bar{X}:=X\sqcup\p X$.

We put the generalized metric~$|\cdot|$ on
$\bar{\rr}:=[-\infty,\infty]$ by declaring
$$|\pm\infty-r|:=\infty\ \mbox{for all}\ \ r\in\rr,\quad 
|(\pm\infty)-(\mp\infty)|:=\infty,\quad |(\pm\infty)-(\pm\infty)|:=0,$$ 
and put the topology on~$\rb$ that makes it a closed interval.
(This topology is {\em not} induced by the generalized metric.)
Now a {\sf geodesic in $\bar\G$} is a continuous isometric embedding
$\al \co  I\to \bar\G$ of an interval $I\se \bar\rr$, where isometry is understood
as preserving the distance that can take infinite values. For example,
the map $\{0\}\to \{a\}$ for any $a\in\p\G$ is trivially a geodesic in~$\bar\G$.
Another example is a {\sf geodesic ray in $\bar{\G}$} converging to $a\in\p\G$,
i.e.\ a continuous isometric embedding $\al\co [0,\infty]\to \bar\G$ with
$\al(0)\in\bar\G$ and $\al(\infty)=a$.

\section{The metric $\dhat$ and the double difference in~$\bar{X}$}\label{s_metric-and-dd}
For the rest of the paper $X$ will denote a hyperbolic metric complex
and we fix a positive integer $\de$ such that the geodesic triangles in~$\G$ are $\de$--fine.

The goal of this section is to construct
a double difference -- a continuous function on quadruples of points in~$\bar{X}$
(Theorem~\ref{dd}). 
It plays the role of the logarithm of the absolute value of the classical cross
ratio in $\cc\cup\{\infty\}=\p\hh^3$, but is going to be defined in a much more general
setting: both on an arbitrary hyperbolic complex and on its ideal boundary.
The double difference will be defined {\em precisely}, rather than
``up to a bounded amount", as often happened for various
notions in hyperbolic groups. First we provide several important auxiliary results.

\subsection{The extended metric~$\dhat$}\label{ss_dhat}
In~\cite{MY}, a metric $\dhat$ was constructed on any hyperbolic group~$\Gamma$.
It was shown that $\dhat$ 
is strongly bolic, $\Gamma$--invariant, quasiisometric to the word metric $d$, 
and satisfies
\begin{ttt}{\rm\cite{MY}}\qua
\label{Cmu}
There exist constants $C\in [0,\infty)$ and $\mu \in[0,1)$ with the following
property. If $a,a',b,b'\in \Gamma$, $ d(a,a')\le 1$, and $ d(b,b')\le 1$, then
$$\big|  \dhat(a,b)- \dhat(a',b)- \dhat(a,b')+ \dhat(a',b')\big|\le C \mu^{d(a,b)}.$$
\end{ttt}
The above theorem is not enough for the purposes of this paper,
and just having a strongly bolic invariant metric is not enough either.
We will show that actually $\dhat$ satisfies stronger properties (Theorem~\ref{t_dhat}).

Now let $X$ be {\em any} hyperbolic complex, and denote $\G$ it 1-skeleton with the path metric~$d$.
We will moreover extend the construction in~\cite{MY} in two ways
to provide a metric~$\dhat=\dhat_X$ which is more general in the following sense:
\begin{itemize}
\item $\dhat$ is defined on all of~$X$, rather than 
just on a discrete group~$\Gamma$.
\item $\dhat$ is invariant under the full isometry group~$\Isom(X)$, rather than
just under~$\Gamma$.
\item The $\Isom(X)$--action on~$X$ is not assumed to be cocompact.
\end{itemize}
The construction in~\cite{MY} utilizes a $\Gamma$--equivariant choice of geodesic
paths~$p[a,b]$ in the Cayley graph, viewed sometimes as a path and sometimes as a cellular 1-chain,
from~$a$ to~$b$ for each pair $(a,b)\in\Gamma^2$. This bicombing~$p[\cdot,\cdot]$
is not good enough for our purposes
since in general it cannot be chosen to be $\Isom(X)$--equivariant.
This problem is fixed by using the 1-chain
$$p'[a,b]:= (\# \Geod(a,b))\inv \sum_{s\in \Geod(a,b)} s,$$
where $\Geod(a,b)$ is the finite set of all geodesic paths in $\G$ from~$a$ to~$b$,
instead of $p[a,b]$ whenever $p[a,b]$ was meant to be a 1-chain
(and making no change when $p[a,b]$ was meant to be a path).
Since $p'[\cdot,\cdot]$ is $\Isom(X)$--equivariant, this provides a $\Isom(X)$--invariant
metric~$\dhat$ on $X^{(0)}$. The uniform local finiteness of~$X$ guarantees that all the arguments
of~\cite{MY} go through with very minor modifications.
(Namely, since the cardinality of balls of radius $7\de$ in $X^{(0)}$
might not be constant, though bounded
above, one needs to change the definition of $star(a)$ \cite[page 100]{MY} to
$$star(a):= \frac{1}{\# B(a,7\de)} \sum_{x\in B(a,7\de)} x.$$
Then the discussion of~\cite[page~111]{MY} should be changed as follows.
Let $\omega_{max}$ be the maximum of cardinalities of the balls of radius $7\de$
in $X^{(0)}$. Without loss of generality we can assume $\be'\le \be$, then
\begin{eqnarray*}
&& \big|star(f_0)+\al x_0- star(f'_0)- \al'x_0{\big|}_1\\
&& \le \big|star(f_0)+\al x_0- \be x_0{\big|}_1+ \big|-star(f'_0)- \al'x_0+ \be'x_0{\big|}_1+
   \big|(\be-\be') x_0{\big|}_1\\
&& = (1-\be)+ (1-\be')+ (\be-\be')= 2(1-\be')\le 2\left(1-\frac{1}{\omega_{max}}\right).
\end{eqnarray*}
The rest of the argument goes through.) 

In particular, Theorem~\ref{Cmu}
still holds for this new metric $\dhat$ on $X^{(0)}$, and $\dhat$ is quasiisometric to~$d$.
Finally, let $\dhat=\dhat_X$ be the extension of~$\dhat$ to all of~$X$
by formula~(\ref{metr-ext}).

The hyperbolicity constant $\de$ and the word-metric $d$ canonically depend
on $X$, and $\dhat$ canonically depends on $\de$ and $d$. This shows in particular
that $\dhat$ is $\Isom(X)$--invariant. Recall that $\Isom(X)$ is the group of
simplicial isometries of $(X,d)$.

\begin{lemma}
For a hyperbolic complex~$X$ and a simplicial map $g\co X\to X$,
the following are equivalent.
\begin{itemize}
\item [\rm(a)] $g$ is a simplicial automorphism of $X$.
\item [\rm(b)] $g$ is a simplicial automorphism of $X$ preserving the word metric $d$ on $X^{(1)}$.
\item [\rm(c)] $g$ is a simplicial automorphism of $X$ preserving $\dhat$ on $X^{(0)}$.
\item [\rm(d)] $g$ is a simplicial automorphism of $X$ preserving $\dhat$ on $X$.
\end{itemize}
\end{lemma}
In other words, there is no difference between simplicial isometries of $(X,d)$
and simplicial isometries of $(X,\dhat)$.

\subsection{Now use $\dhat$ instead of $d$ everywhere}\label{ss_dhat-everywhere}
From now on we deal with hyperbolic complexes and we
change the notations of~\ref{ss_dd-gp}, redefining 
everything in terms of the metric $\dhat$:
the double difference $\left<\cdot,\cdot|\cdot,\cdot\right>\co  X^4\to\rr$ is
\begin{equation}\label{dd-def-new}
\left<a,a'|b,b'\right>:=
  \frac{1}{2} \big( \dhat(a,b)- \dhat(a',b)- \dhat(a,b')+ \dhat(a',b') \big)
\end{equation}
and the Gromov product is
\begin{equation}\label{d_gr-pr-new}
\left<a|b\right>_c:= \left<a,c|c,b\right>=
\frac{1}{2} \big( \dhat(a,c)+ \dhat(b,c)- \dhat(a,b)\big).
\end{equation}
With definition~(\ref{dd-def-new}), Theorem~\ref{Cmu} says that 
there exist $C\in [0,\infty)$ and $\mu \in[0,1)$ such that for all
$a,a',b,b'\in X^{(0)}$ with $ d(a,a')\le 1$ and $ d(b,b')\le 1$, 
$\big| \left<a,a'|b,b'\right> \big|\le C \mu^{d(a,b)}$.

From now on the letter $d$ will always stand for the canonical word metric
as in~\ref{ss_word-metric}. Define the corresponding functions
with respect to the word metric $d$:
\begin{eqnarray}\label{dd-def-metric-d}
&& (a,a'|b,b'):=
  \frac{1}{2} \big( d(a,b)- d(a',b)- d(a,b')+ d(a',b') \big),\\
\nonumber && (a|b)_c:= (a,c|c,b)=
\frac{1}{2} \big( d(a,c)+ d(b,c)- d(a,b)\big).
\end{eqnarray}
These functions can be partially extended to infinity; the extension is 
usually not continuous and is defined up to an additive constant.

\noindent {\bf Remark}\qua  The linearity formula~(\ref{metr-ext}) implies that
the double difference in $X$ is determined by its values on the vertices of $X$.
Explicitly, if $a=\sum_x \al_x x$, $a'=\sum_{x'} \al_{x'} x'$,
$b=\sum_y \be_y y$, $b'=\sum_{y'} \be_{y'} y'$
are convex combinations of vertices, then
\begin{equation}\label{dd-ext}
\left<a,a'|b,b'\right>=
\sum_x\sum_{x'}\sum_y\sum_{y'} \al_x\al_{x'}\be_y\be_{y'}
\left<x,x'|y,y'\right>.
\end{equation}
Starting from the metric $\dhat$ (rather than from~$d$) on~$X$ we define
$\ell$, $\becross$, the pseudometric $\dcross$ and the metric $d_*$ on $\sx$ as in~\ref{ss_d^times}
and~\ref{ss_metric-d*}, i.e.\ we denote
$d_*:= \sj\dhat$ (rather than $d_*:= \sj d$).

\subsection{Examples of hyperbolic complexes}\label{ss_examples-of-hyp-c}
\begin{itemize}
\item [\rm(1)] As in~\ref{ss_metric-complexes}(1), the Cayley graph of a hyperbolic group can
be considered a hyperbolic complex with respect to the metric~$\dhat_\G$.
\item [\rm(2)] The universal cover of a compact triangulated
smooth manifold with hyperbolic
fundamental group with respect to the metric $\Psi(d)$ as in~\ref{ss_metric-complexes}(2),
where $d$ is the intrinsic metric coming from the Riemannian structure.
\item [\rm(3)] The universal cover of a compact triangulated manifold with hyperbolic
fundamental group. One can take either the word metric $d$ or the canonical
metric $\dhat$ defined in~\ref{ss_dhat}, both induced from the $0$--skeleton
(see~\ref{ss_metric-complexes}(3)). The two metrics are quasiisometric
but $\dhat$ behaves better at infinity.
\end{itemize}

\subsection{Geodesics and nearest points}\label{ss_geod-and-near}
``Geodesic'' will always refer to the word metric~$d$ in~$\G=X^{(1)}$.
Let $a,a',b,b'\in \bar{\G}$, $\al\in \Geod(a,a')$, 
$\be\in \Geod(b,b')$, so that $\al$ and $\be$ are isometric
embeddings of intervals:  $\al\co I\to\bar\G$, $\beta\co J\to\bar\G$.
Pick distance-minimizing vertices $a_0\in\al$ and
$b_0\in\beta$, i.e.\ such that $ d(a_0,b_0)= d(\al,\beta)$.
Choose the parameterizing intervals $I$ and $J$ of~$\al$ and $\be$ containing 0 so that
$\al(0)=a_0$, $\al(-M)=a$, $\al(M')=a'$,
$\beta(0)=b_0$, $\beta(-N)=b$, $\beta(N')=b'$ for some non-negative 
(possibly infinite) $M$, $M'$, $N$, $N'$.

\begin{lemma}\label{ab6de}$\phantom{99}$
\begin{itemize}
\item [\rm(a)] With the above notations, if~$ d(\al,\be)\ge 2\de$, then for all $i\in I$ and $j\in J$,
$$ d(\al(i),\beta(j))\ge |i|+  d(\al,\beta)+ |j| -6\de.$$
In particular, $d(a,b)\ge d(a,a_0)+ d(a_0,b_0)+ d(b_0,b)-6\de$.
\end{itemize}
Now assume $a=a'=a_0$, and $d(a_0,\be)$ and $j\in J$ are arbitrary. Then
\begin{itemize}
\item [\rm(b)] $ d(a_0,\be(j))\ge |j|+ d(a_0,\be)-2\de$,
\item [\rm(c)] $ d([a_0,b_0],\be(j)])\ge |j|-2\de$,
\item [\rm(d)] $ d([a_0,b],\be(j))\ge j-3\de$.
\end{itemize}
\end{lemma}
\proof
\noindent {\bf (a)}\qua 
From symmetry, it suffices to show the lemma when $i\ge 0$ and $j\ge 0$.
Draw geodesics $[\al(i), \be(j)]$ and $[a_0,\be(j)]$, and
inscribe triples of points in the triangles $\{a_0, \al(i), \be(j)\}$ and
$\{a_0, \be(j), b_0\}$ as shown on Figure~\ref{f_ab6de}.
\setlength{\unitlength}{1cm} 
\begin{figure}[ht!]
  \begin{center}
   \begin{picture}(10,4)

\qbezier(1.6,-0.3)(2.75,1.3)(0.5,4)
\qbezier(8,-0.5)(7.5,2.5)(9,4)
\qbezier(2,1)(4.8,1)(7.88,1)
\qbezier(1,3.35)(4.5,0.5)(8.355,3.1)
\qbezier(2,1)(6,1.3)(8.355,3.1)

\put(0.5,4){\circle*{0.1}}
\put(0.1,3.9){\tiny$a'$}

\put(1.6,-0.3){\circle*{0.1}}
\put(1.2,-0.3){\tiny$a$}

\put(2,1){\circle*{0.1}}
\put(.5,1){\tiny$a_0=\al(0)$}

\put(1,3.35){\circle*{0.1}}
\put(0.3,3.2){\tiny$\al(i)$}

\put(9,4){\circle*{0.1}}
\put(9.2,3.9){\tiny$b'$}

\put(8,-0.5){\circle*{0.1}}
\put(8.3,-0.5){\tiny$b$}

\put(7.88,1){\circle*{0.1}}
\put(8.1,1){\tiny$b_0=\be(0)$}

\put(8.355,3.1){\circle*{0.1}}
\put(8.6,3){\tiny$\be(j)$}

\put(2.8,1.09){\circle*{0.1}}
\put(2.6,1.28){\tiny$v_1$}

\put(1.874,1.7){\circle*{0.1}}
\put(2,1.7){\tiny$v_2$}

\put(2.82,2.25){\circle*{0.1}}
\put(2.63,1.9){\tiny$v_3$}

\put(2.81,0.995){\circle*{0.1}}
\put(2.7,0.7){\tiny$v'$}

\put(7.2,1){\circle*{0.1}}
\put(7,1.22){\tiny$w_1$}

\put(7.93,1.7){\circle*{0.1}}
\put(7.5,1.65){\tiny$w_2$}

\put(6.9,2.2){\circle*{0.1}}
\put(6.6,1.7){\tiny$w_3$}

   \end{picture}
  \end{center}
\caption{Illustration for Lemma~\ref{ab6de}}\label{f_ab6de}
\end{figure}
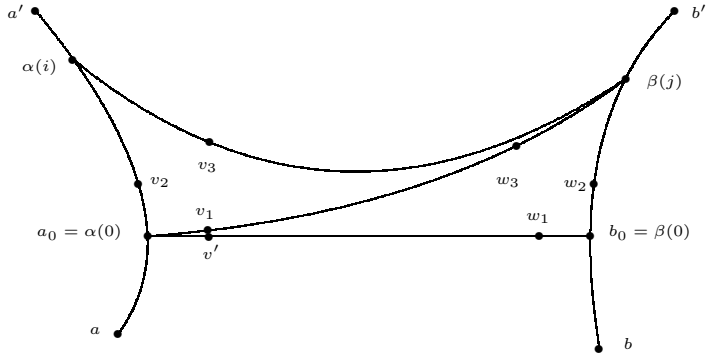
Pick $v'\in[a_0,b_0]$ with $ d(a_0,v')= d(a_0, v_1)$.
The vertex~$v_1$ is $\de$--close to $v_2$, and therefore to~$\al$.
If $ d(a_0,v')>  d(a_0,w_1)$ was true, then $v'$ would be $\de$--close to $b_0$,
so  $ d(\al,\be)< 2\de$, which contradicts our assumption.
Therefore $ d(a_0, v')\le  d(a_0, w_1)$.
It is an easy exercise to see from the figure that
$ d(b_0,w_2)=  d(b_0,w_1)\le \de$ and $ d(a_0, v_1)=  d(a_0,v_2)= d(a_0,v')\le 2\de$
(otherwise $[a_0,b_0]$ could be shortened), therefore
\begin{eqnarray*}
&&   d(\al(i),\be(j)) =   d(\al(i),v_2)+  d(v_1, \be(j))\\
&&  = d(\al(i),v_2)+  d(v',w_1)+  d(w_2, \be(j))\\
&&  \ge  \big[ d(\al(i),a_0)- 2\de\big]+
                \big[ d(a_0,b_0)-3\de\big]+\big[ d(b_0,\beta(j))-\de\big]\\
&&   = i+ d(\al,\beta)+j -6\de.
\end{eqnarray*}

\noindent {\bf (b)}\qua Since $a=a'=a_0$, we also have $a=a'=a_0=v_1=v_2=v_3=v'$,
and the inequality follows as in (a) by a similar, and simpler, argument.

\noindent {\bf (c)}\qua From symmetry it suffices to show the inequality in the case $j\ge 0$.
If $j-2\de\le 0$ then the inequality is obvious.
Now suppose to the contrary that there exists $j> 2\de$ and a point $x\in [a_0,b_0]$ such that
$ d(x,\be(j))< j-2\de$. If $ d(b_0,x)\le  d(b_0,w_1)$,
then $ d(x,\be(j))\ge  d(b_0,\be(j))- d(b_0,x)\ge j-\de$
which contradicts the assumption. Then $ d(b_0,x)\ge  d(b_0,w_1)$,
so by the fine-triangles condition, there exists $y\in [a_0,\be(j)]$
with
$$ d(a_0,y)= d(a_0,x)\le  d(a_0,w_3)\quad \mathrm{and}\quad  d(x,y)\le\de.$$
Then
$$ d(x,\be(j))\ge  d(y,\be(j))-\de\ge  d(w_3,\be(j))-\de=  d(w_2,\be(j))-\de\ge j-2\de$$
provides a contradiction.

\noindent {\bf (d)}\qua If $j-3\de\le 0$, then the inequality is obvious, so we assume
$j>3\de$. From (c), $ d([a_0,b_0],\be(j))\ge |j|-2\de=j-2\de$,
and since $\be$ is geodesic, $ d([b_0,b],\be(j))= |j|=j$. 
Since $[a_0,b]$ lies in the $\de$--neighborhood of $[a_0,b_0]\cup[b_0,b]$,
\begin{align*}
& d([a_0,b],\be(j))\ge d([a_0,b_0]\cup[b_0,b],\be(j))-\de \\
&\qquad\ge\min\{d([a_0,b_0],\be(j)), d([b_0,b'],\be(j))\}-\de\ge j-3\de.\tag*{\qed}
\end{align*}

Given $y\in\bar{\G}$ and geodesics~$\al$ and $\g$ in~$\bar\G$, $\g$ is called
a {\sf distance-minimizing}\label{i_distance-min-geod-y-al} geodesic from $y$ to $\al$ if $\g$ 
starts at $y$, terminates at a point $x\in\al$,
and for every $z\in\g$, $d(z,x)=d(z,\al)$. A~distance-minimizing
geodesic exists for every pair $(y,\al)$. The set of all such terminal points~$x$
over all distance-minimizing geodesics from $y$ to $\al$, denoted
$\mathrm{np}[\al|y]$,\label{i_nearest-point-proj} is called 
the {\sf nearest point projection of~$y$ to~$\al$}.
When $\al=[a,a']$ we will use the notation $\mathrm{np}[a,a'|y]$
for $\mathrm{np}[\al|y]$.
For example, if $a\in\p\G$ then $\{\mathrm{0}\}\to \{a\}$
is the only distance-minimizing geodesic from~$a$ to~$[a,a']$,
so $\mathrm{np}[a,a'|a]=\{a\}$.

Given two geodesics~$\al$ and~$\be$ in $\bar{\G}$,
$\g$ is called a {\sf distance-minimizing geodesic}\label{i_distance-min-geod-al-be} 
from $\al$ to $\be$
if $\g$ is a geodesic in~$\bar{\G}$, $\g$ starts at a point $x\in\al$,
terminates at a point $y\in\be$, and for every $z\in\g$, $d(z,x)=d(z,\al)$
and $d(z,y)=d(z,\be)$.
Such a pair~$(x,y)\in \al\times\be$ will be called 
a {\sf distance-minimizing pair}\label{i_distance-min-pair}
for $(\al,\be)$. Every pair $(\al,\be)$ of geodesics in~$\bar{\G}$ admits
a distance-minimizing pair.

\begin{lemma}\label{l_xyz}
Let $b,b'\in\bar{\G}$, $\be\in \Geod(b,b')$, $x,y\in\G$,
$x'\in \mathrm{np}[\be|x]$, $y'\in \mathrm{np}[\be|y]$ \ (see Figure~\ref{f_xyz}). 
\begin{itemize}
\item [\rm(a)] If $ d(x',y')\ge 2\de$, then
$d(x,y)\ge  d(x,x')+  d(x',y')+  d(y',y)- 6\de$.
\item [\rm(b)] For any positive integer $n$, if $d(x,y)< n\de$, then $d(x',y')<(n+6)\de$.
\item [\rm(c)] Take any $\al\in \Geod(x,y)$ and let $\be'$ be the part of $\be$ between
  $x'$ and $y'$. Then for any $z\in\al$ and $z'\in \mathrm{np}[\be|z]$,
    $d(\be',z')\le 8\de$.
\end{itemize}
\end{lemma}

\setlength{\unitlength}{1cm} 
\begin{figure}[ht!]
  \begin{center}
   \begin{picture}(11,5)

\put(0,2.2){$\scriptstyle (1)$}
\qbezier(0.5,0)(2.5,0)(4.5,0) 
\put(0.5,0){\circle*{0.1}}  \put(0.4,0.15){\tiny$b$}  
\put(4.5,0){\circle*{0.1}}  \put(4.45,0.15){\tiny$b'$}

\put(4,1.2){\footnotesize$\G$}

\put(1,-0.3){\tiny$\be$}
\put(2,0){\circle*{0.1}} \put(1.8,-0.3){\tiny$x'$}
\put(2.6,0){\circle*{0.1}}  \put(2.6,-0.3){\tiny$y'$}

\put(1.5,4.5){\circle*{0.1}} \put(1.1,4.4){\tiny$x$} \qbezier(1.5,4.5)(2,2.25)(2,0)
\put(3.5,5){\circle*{0.1}}  \put(3.7,5){\tiny$y$} \qbezier(3.5,5)(2.6,2.5)(2.6,0)
\qbezier(3.5,5)(2.5,2.5)(2,0)

\qbezier(1.5,4.5)(2,2.6)(2.25,2.6) \qbezier(2.25,2.6)(2.5,2.6)(3.5,5)

\put(2.25,2.62){\circle*{0.1}} \put(2.1,2.85){\tiny$w_1$}
\put(1.857,2.42){\circle*{0.1}} \put(1.3,2.4){\tiny$w_3$}
\put(2.59,2.36){\circle*{0.1}} \put(2.9,2.32){\tiny$w_2$}

\put(2.7,3.8){\tiny$\al$}

\put(2.32,0){\circle*{0.1}} \put(2.18,-0.3){\tiny$v_1$}
\put(2.075,0.32){\circle*{0.1}} \put(1.6,0.35){\tiny$v_2$}
\put(2.605,0.26){\circle*{0.1}} \put(2.8,0.25){\tiny$v_3$}

\put(5.5,2.2){$\scriptstyle (2)$}
\qbezier(5.5,0)(8.25,0)(11,0) 
\put(5.5,0){\circle*{0.1}}  \put(5.4,0.15){\tiny$b$}
\put(11,0){\circle*{0.1}}  \put(10.95,0.15){\tiny$b'$}

\put(11,1.2){\footnotesize$\G$}

\put(6,-0.3){\tiny$\be$}
\put(7,0){\circle*{0.1}}  \put(7,-0.3){\tiny$x'$}
\put(10,0){\circle*{0.1}} \put(10,-0.3){\tiny$y'$}

\put(6.5,4.5){\circle*{0.1}}  \put(6.1,4.5){\tiny$x$}   \qbezier(6.5,4.5)(7,2.1)(7,0)
\put(10.5,5){\circle*{0.1}}  \put(10.7,5){\tiny$y$} \qbezier(10.5,5)(10.1,2.5)(10,0)

\qbezier(7,0)(8.8,0.1)(9.4,0.9) \qbezier(9.4,0.9)(9.92,1.6)(10.5,5)

\qbezier(6.5,4.5)(7.08,1.32)(7.85,0.8) \qbezier(7.85,0.8)(8.6,0.3)(9.17,0.9)
  \qbezier(9.17,0.9)(9.8,1.5)(10.5,5)

\put(7.78,0.86){\circle*{0.1}} \put(7.75,1.1){\tiny$w_1$}
\put(7.8,0.086){\circle*{0.1}} \put(7.5,0.3){\tiny$w_2$}
\put(7.81,0){\circle*{0.1}} \put(7.7,-0.3){\tiny$w'_2$}
\put(6.989,0.68){\circle*{0.1}} \put(6.5,0.7){\tiny$w_3$}

\put(9.5,0){\circle*{0.1}} \put(9.4,-0.3){\tiny$v_1$}
\put(9.234,0.7){\circle*{0.1}} \put(9.32,0.45){\tiny$v_2$}
\put(9.105,0.823){\circle*{0.1}} \put(8.9,1.2){\tiny$v'_2$}
\put(10.03,0.5){\circle*{0.1}} \put(10.2,0.55){\tiny$v_3$}

\put(9.7,3){\tiny$\al$}

   \end{picture}
  \end{center}
\caption{Illustration for Lemma~\ref{l_xyz}}
\label{f_xyz}
\end{figure}
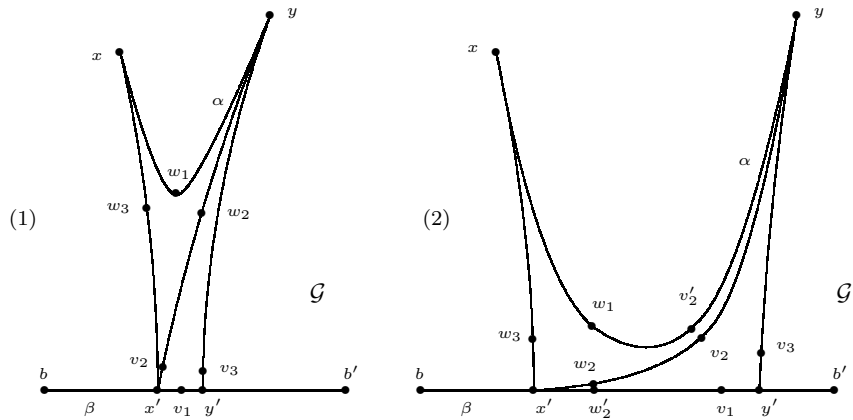

\begin{proof}
{\bf (a)}\qua
Set the notations as in Figure~\ref{f_xyz}(2).
If $d(x',v_1)< d(x',w'_2)$, then $w'_2$ would be $\de$--close to $y'$,
so $ d(x',y')< 2\de$, which contradicts our assumption. So $ d(x',w'_2)\le  d(x',v_1)$.
We have
$$d(x', w_2)=  d(x',w_3)= d(x',w'_2)\le 2\de \mbox{\quad and\quad}
    d(y',v_3)=  d(y',v_1)\le \de$$
(otherwise $[x,x']$ or $[y,y']$ could be shortened), hence
\begin{eqnarray*}
&&   d(x,y) =   d(x,w_3)+  d(w_2, y)
     = d(x,w_3)+  d(w'_2,v_1)+  d(v_3, y)\\
&&    \ge  \big[ d(x,x')- 2\de\big]
 + \big[ d(x',y')-3\de\big]+\big[ d(y',y)-\de\big]\\
&&  =d(x,x')+ d(x',y')+ d(y',y)-6\de.
\end{eqnarray*}
\noindent {\bf (b)}\qua
The proof is by contradiction: if $d(x',y')\ge (n+6)\de$, then $d(x',y')\ge 2\de$
hence by (a), $d(x,y)\ge d(x,x')+ d(x',y')+ d(y',y)- 6\de\ge d(x',y')-6\de\ge n\de$.

\noindent {\bf (c)}\qua
Draw the geodesic $[x',y]$ and inscribe triples of points in the triangles
$\{x,y,x'\}$ and $\{y,x',y'\}$.
There are two cases as on~Figure~\ref{f_xyz} depending on the order
of $v_2$ and $w_2$ along $[x',y]$. 
Using the $\de$--fine property it is easy to see from the figure that in either case
$\al$ lies in the $2\de$--neighborhood of $[x,x']\cup \be'\cup [y',y]$.
Pick any $z\in\al$ and $z'\in \mathrm{np}[\be|z]$, then there exists
$w\in [x,x']\cup \be'\cup [y',y]$ with $d(z,w)\le 2\de$.
But for any $w\in [x,x']\cup \be'\cup [y',y]$, there is a projection
$w'\in \mathrm{np}[\be|w]$ which also belongs to $\be'$.
By (b) applied to $z$ and $w$, $d(z',w')\le 8\de$, so $d(\be',z')\le 8\de$.
\end{proof}

\subsection{Neighborhoods in $\bar{\G}$ and $\bar{X}$}
\label{ss_neighborhoods-in-barG}
As before, $X$ will be a hyperbolic complex, and $\G:=X^{(1)}$ with the word metric~$d$.
We present a description of the topology on $\bar{X}=X\cup \p X$ in terms of half-spaces.
The idea of this presentation is due to Cannon, and the details were
written out by Swenson~\cite{Swenson1993}. It was proved in~\cite{Swenson1993}
that half-spaces generate the same topology as the one introduced by Gromov~\cite{Gr2};
see also~\cite{CS1998} for a list of properties.

Pick a basepoint $x_0\in\G^{(0)}$.
Given $a\in\p\G$ and $t\in\rr$, denote 
\begin{eqnarray*}\label{i_U+-at}
&&U^+(a,t):=\{x\in \bar{\G}\ |\  \exists\, x'\in\mathrm{np}[x_0,a|x]\quad d(x_0,x')\ge t \},\\
&&U^-(a,t):=\{x\in \bar{\G}\ |\  \exists\, x'\in\mathrm{np}[x_0,a|x]\quad d(x_0,x')\le t \}.
\end{eqnarray*}
Now let $V^\pm(a,t):= U^\pm(a,t)$
\label{i_V+-at} be the union
of simplices in $X$ whose vertices lie in $U^\pm(a,t)$.
The sets $V^+(a,t)\se \bar{X}$ form a fundamental system of closed neighborhoods
of $a\in\p X$ in $\bar{X}$.
\begin{lemma}\label{U-U+}
For any $a\in\p\G$ and any $s,t\in\rr$,
$$d\big( U^-(a,s), U^+(a,t)\big)\ge t-s-6\de.$$
\end{lemma}
\begin{proof}
If $t< s+6\de$ then the inequality is obvious.
Now pick arbitrary
$$s\in\rr,\quad t\in[s+6\de,\infty),\quad x\in U^-(a,s),\quad y\in U^+(a,t).$$
By the definition of $U^{\pm}$, we can choose 
$x'\in \mathrm{np}[x_0,a|x]$ and $y'\in \mathrm{np}[x_0,a|y]$
so that $d(x_0,x')\le s$ and $d(x_0,y')\ge t$.
Since $ d(x',y')\ge t-s\ge 6\de$, using Lemma~\ref{l_xyz},
$d(x,y)\ge  d(x', y')- 6\de\ge t-s-6\de$.
\end{proof}

For $S_1, S_2\se \bar{\G}$, $[S_1,S_2]$\label{i_[S1-S2]} will denote the union of the images of
all geodesics in~$\bar\G$ connecting points of~$S_1$ to points in~$S_2$.
\begin{lemma}\label{Uconvex}
Let $a\in\p\G$ and $s\in [0,\infty)$, then
$$[U^+(a,s), U^+(a,s)]\se U^+(a,s-8\de)\hbox{ and }[U^-(a,s), 
U^-(a,s)]\se U^-(a,s+8\de).$$
\end{lemma}
\begin{proof}
Pick any $x,y\in U^+(a,s)$, $\al\in \Geod(x,y)$.
By the definition of $U^+$, there exist projections
$x'\in\mathrm{np}[x_0,a|x]$ and $y'\in\mathrm{np}[x_0,a|y]$
with $d(x_0,x')\ge s$ and $d(x_0,y')\ge s$. 
Let $\be'$ be the interval in $[x_0,a]$ between $x'$ and $y'$.
By Lemma~\ref{l_xyz}(c), the projection of $\al$ into $[x_0,a]$ lies 
in the $8\de$--neighborhood of $\be'$, so $\al\se U^+(a,s-8\de)$.
The second inclusion is treated similarly.
\end{proof}

\subsection{Properties of~$\dhat$}\label{ss_prop-of-dhat}
\begin{ttt}\label{t_dhat}
For any hyperbolic complex~$X$ there is a canonically associated metric
$\dhat=\dhat_X$ with the following properties. 
\begin{itemize}
\item [\rm(a)] $\dhat$ is quasiisometric to the word metric~$d_X$, i.e there exist
$A\in[1,\infty)$ and $B\in [0,\infty)$ depending only on~$X$ such that for all $x,y\in X$,
$$d(x,y)/A-B\le \dhat(x,y)\le A d(x,y)+B.$$
\item [\rm(b)] $\dhat$ is $\Isom(X)$--invariant: $\dhat(gx,gy)=\dhat(x,y)$ for $g\in\Isom(X)$.
\item [\rm(c)] $(X,\dhat)$ is a metric complex, that is $\dhat$ is induced from
the 0-skeleton.
\item [\rm(d)] There exist constants $L\in[0,\infty)$ and $\mu\in [0,1)$ depending only on~$X$
with the following property.
If $a,a',b,b'\in X^{(0)}$, $\al\in \Geod(a,a')$, $\be\in \Geod(b,b')$ and
$d\big(\al,\be\big)\ge 2\de$, then
$\left| \left<a,a'|b,b'\right> \right|\le L\mu^{ d(\al,\be)}$.
\end{itemize}
\end{ttt}
\begin{proof} 
Take the metric $\dhat$ as described in~\ref{ss_dhat}.
(a) and (b) are shown in theorems~17 and~6(2) in~\cite{MY} for $\G$, this easily
implies (a) and (b) for the whole~$X$.
(c) holds by definition. It only remains to show (d).
Pick arbitrary vertices $a$, $a'$, $b$, $b'$ satisfying the hypotheses. 
Set notations as in Figure~\ref{f_ab6de}: pick a distance-minimizing pair
$(a_0,b_0)$ for $(\al,\be)$ and choose parametrizations
of $\al$ and $\be$ such that $a=\al(-M)$,
$a'=\al(M')$, $b=\be(-N)$, $b'=\be(N')$ for some non-negative integers
$M$, $M'$, $N$, $N'$. Then by Lemma~\ref{ab6de} and Theorem~\ref{Cmu},
taking the sums of geometric series,
\begin{eqnarray*}
&&\left|\left<a,a'| b, b'\right>\right|
 = \left|\sum_{i=-M}^{M'-1} \sum_{j=-N}^{N'-1}
  \big<\al(i), \al(i+1)| \be(j), \be(j+1) \big> \right|\\
&&   \le \sum_{i=-M}^{M'-1}
\sum_{j=-N}^{N'-1}
  C\mu^{ d(\al(i), \be(j))}
 \le  \sum_{i=-\infty}^{\infty} \sum_{j=-\infty}^{\infty}
  C\mu^{|i|+|j|+ d(\al, \be)-6\de}\\
&& \le 4\frac{C\mu^{ d(\al,\be)-6\de}}{(1-\mu)^2}= L\mu^{ d(\al,\be)},
\end{eqnarray*}
where we denoted $\ds L:= {4C\mu^{-6\de}}/{(1-\mu)^2}$.
\end{proof}

\begin{ppp}\label{p_uvw}
Let $\G$ be a hyperbolic graph of bounded valence.
There exists $C'\in[0,\infty)$ depending only on~$\G$ such that
if $u,v,w\in\G^{(0)}$ and $w$ lies on a geodesic in~$\G$ connecting~$u$ to~$v$,
then $|\dhat(u,v)-\dhat(u,w)-\dhat(w,v)|\le C'$.
\end{ppp}
\noindent This immediately follows from the definition of~$\dhat$ and Proposition~10(b)
in~\cite{MY}. The proposition says in other words that geodesics in $(\G,d)$
can be parameterized to become {\sf $^+$geodesic}\label{i_plus-geod}
in $(\G,\dhat)$, that is geodesic up to an additive constant,
where the constant depends only on $\G$. Note though that $d$ and $\dhat$
are {\em not} necessarily $^+$equivalent.

\begin{lemma}\label{b'a0a0b}
There exists $P\in[0,\infty)$ depending only on $X$ with the following property.
Let $a_0\in X^{(0)}$, $b,b'\in X^{(0)}\cup\p X$, $\be\in \Geod(b,b')$
and $b_0\in \mathrm{np}[\be|a_0]$. If $b,b'\in\p X$, assume $b\not= b'$.
Then $|\left<b|b'\right>_{a_0}- \dhat(a_0,b_0)|\le P$.
\end{lemma}
\noindent Note that the statement is about $\dhat$, but nearest points
are defined using $d$.
\begin{proof}
Assume first that $a,b,b'\in X^{(0)}$.
We use hyperbolicity and Proposition~\ref{p_uvw}:
\begin{eqnarray*}
&& \left<b|b'\right>_{a_0}- \dhat(a_0,b_0) \ \ssimplus\   
\frac{1}{2} \big( \dhat(a_0,b)+ \dhat(a_0,b')- \dhat(b,b')\big)-
  \dhat(a_0,b_0) \\
&& \ssimplus\ 
\frac{1}{2} \big( \dhat(a_0,b_0)+ \dhat(b_0,b)+ 
  \dhat(a_0,b_0) + \dhat(b_0,b')  - \dhat(b,b')\big)-
  \dhat(a_0,b_0)\\
&& =\ 
\frac{1}{2} \big(  \dhat(b_0,b)+ \dhat(b_0,b')  - \dhat(b,b')\big)
\ \ssimplus\  0,
\end{eqnarray*}
where $\simplus$ is $^+$equivalence with a constant that depends
only on $X$. This extends to the case when $b$ or $b'$ is
in $\p X$.
\end{proof}

\subsection{Extension of double difference}\label{ss_ext-of-dd}
Let $(a, a', b, b')\in \bar{X}^4$.  A {\sf $\p X$--triple in}\label{i_px-triple}
$(a, a', b, b')$ is a set of three distinct letters taken 
from $a, a', b, b'$ such that each letter represents
a point in $\p X\se \bar{X}$.
A $\p X$--triple is {\sf trivial}\label{i_trivial-px-triple} if the three letters represent the same
point in~$\p X$.
Denote 
$${\bar{X}}^{\diamond}:=\big\{(a,a',b,b')\in\bar{X}^4\ \big|\ 
   \mbox{each } \p X\mbox{-triple in } (a, a', b, b')
    \mbox{ is non-trivial} \big\}.$$
We have $X^4\se \bar{X}^{\diamond}\se \bar{X}^{4}$.
The topology on $\bar{X}^{\diamond}$ is induced by the last inclusion, and
$X^4$ is dense in both $\bar{X}^{\diamond}$ and $\bar{X}^{4}$.

Recall that $\bar{\rr}:=[-\infty, \infty]$ is the two-point compactification of
$\rr$. 

\begin{ttt}\label{dd}
Suppose $X$ is a hyperbolic complex, then the double difference
$\left<\cdot,\cdot|\cdot,\cdot\right>\co X^4\to \rr$ with respect to $\dhat$
defined in~(\ref{dd-def-new}) extends to a continuous $\Isom(X)$--invariant 
function
$\left<\cdot,\cdot|\cdot,\cdot\right>\co\bar{X}^{\diamond}\to \bar{\rr}$
with the following properties. 
\begin{itemize}
\item [\rm(a)]  $\left<a,a'| b,b'\right>= \left<b,b'| a,a'\right>$.
\item [\rm(b)]  $\left<a,a'|b,b'\right>= -\left<a',a| b,b'\right>= -\left<a,a'| b',b\right>$.
\item [\rm(c)]  $\left<a,a|b,b'\right>=0$, $\left<a,a'|b,b\right>=0$.
\item [\rm(d)] $\left<a,a'| b,b'\right>+ \left<a',a''| b,b'\right>= \left<a,a''| b,b'\right>$,
  where by convention we allow $\pm\infty\mp\infty=r$ and 
  $\pm\infty+ r=\pm\infty$ for any $r\in\rr$, and
  $\pm\infty\pm\infty=\pm\infty$.
\item [\rm(e)] $\left<a,b|c,x\right>+\left<b,c|a,x\right>+\left<c,a|b,x\right>=0$
  with the same convention.
\item [\rm(f)]  $\left<a,a'|b,b'\right>= \infty$ if and only if $a=b'\in \p X$ or $a'=b\in \p X$.
\item [\rm(g)]  $\left<a,a'|b,b'\right>= -\infty$ if and only if $a=b\in \p X$ or $a'=b'\in \p X.$
\item [\rm(h)] Let $(\cdot,\cdot|\cdot,\cdot)$ be the double difference with respect to
  the word metric $d$ as in~(\ref{dd-def-metric-d}). Then 
  $\left<\cdot,\cdot|\cdot,\cdot\right>$ and $(\cdot,\cdot|\cdot,\cdot)$ are
  $^{\times +}$equivalent as functions on $\bar{X}^{\diamond}$.
\end{itemize}
\end{ttt}
\begin{proof}
For $(a,a',b,b')\in \bar{X}^4$,
the pairs of letters $(a,b)$, $(a',b')$, $(a,b')$, $(a',b)$ and their inverses will be called
{\sf side pairs}. \label{i_side-pair}
A {\sf side $\p X$--pair} is a side pair in which each letter represents
a point in~$\p X$. A side $\p X$--pair is {\sf trivial} if the two letters represent the same
point. Consider the intermediate set
$$S:=\big\{(a,a',b,b')\in\bar{X}^4\ \big|\ 
   \mbox{each } \mbox{side \p X-pair in } (a, a', b, b')
    \mbox{ is non-trivial} \big\};$$
one checks that $X^4\se S\se \bar{X}^{\diamond}$.

We first extend
$\left<\cdot,\cdot|\cdot,\cdot\right>$ to a continuous function $S\to\rr$.
Since $X^4$ is dense in $S$ and $\rr$ is regular, 
it suffices to prove the existence of 
\begin{equation}
\label{lim_aa'bb'}
\lim \left<a,a'|b,b'\right>\quad\mbox{as}\quad
    (a,a',b,b')\to (\bar{a},\bar{a}',\bar{b},\bar{b}')\ \ \mbox{along}\ \ X^4
\end{equation}
in ~$\rr$ for each $(\bar{a},\bar{a}',\bar{b},\bar{b}')\in S\setminus X^4$ \ 
(see~\cite[I \S 8 \No5, Theorem~1]{Bourbaki}).
In each such $(\bar{a},\bar{a}',\bar{b},\bar{b}')$
at least one of $\bar{a},\bar{a}',\bar{b},\bar{b}'$ is in $\p X$ and there are
no trivial side $\p X$--pairs.
We will present the proof in the case when $\bar{a},\bar{a}',\bar{b},\bar{b}'$
are all in~$\p X$; the other cases are similar. 

Since $\bar{\G}$ is Hausdorff, there is a number $s\in\rr$ such that 
\begin{eqnarray*}
&&  U^+(\bar{a},s)\cap  U^+(\bar{b}',s)
   =U^+(\bar{b}',s)\cap U^+(\bar{a}',s)\\
&& =U^+(\bar{a}',s)\cap U^+(\bar{b},s)
   =U^+(\bar{b},s)\cap  U^+(\bar{a},s) =\emptyset,
\end{eqnarray*}
So in particular
$U^+(\bar{b},s)\cup U^+(\bar{b}',s)\se U^-(\bar{a},s)$, and
Lemma~\ref{Uconvex} implies 
$$[U^+(\bar{b},s+8\de), U^+(\bar{b}',s+8\de)]\se
[U^-(\bar{a},s),U^-(\bar{a},s)]\se U^-(\bar{a},s+8\de),$$
hence replacing $s$ with $s+6\de$ we can assume
$[U^+(\bar{b},s), U^+(\bar{b}',s)]\se U^-(\bar{a},s)$.
The same argument applies to the cyclic permutations of
$(\bar{a},\bar{b},\bar{a}',\bar{b}')$.
Denote 
$$\al:=[x_0,\bar{a}],\quad \al':=[x_0,\bar{a}'],\quad
   \be:=[x_0,\bar{b}],\quad \be':=[x_0,\bar{b}'].$$
Pick any $i\ge s+8\de$. We have
$\al(i),\al(i+1)\in U^+(\bar{a},i)$,
therefore the edge $[\al(i),\al(i+1)]$ also lies in $U^+(\bar{a},i)$. Also 
$\be(i)\in U^+(\bar{b},i)\se U^+(\bar{b},i)\se U^+(\bar{b},s)$
and similarly $\be'(i+1)\in U^+(\bar{b}',s)$, hence
$[\be(i),\be'(i+1)]\se [U^+(\bar{b},s), U^+(\bar{b}',s)]\se U^-(\bar{a},s)$
(see Figure~\ref{f_G4toS}).
By Lemma~\ref{U-U+},
$$ d\big([\be(i),\be'(i+1)], [\al(i+1),\al(i)]\big)
  \ge  d\big(U^-(\bar{a}, s), U^+(\bar{a}, i)\big)  \ge i-s-6\de \ge 2\de.$$
hence by Theorem~\ref{t_dhat}(d),
$$\big|\left<\al(i+1),\al(i)\, |\, \be(i),\be'(i+1)\right>\big|\le 
   L\mu^{d([\be(i),\be'(i+1)], [\al(i+1),\al(i)])}\le
  L\mu^{i-s-6\de},$$
and similarly for the cyclic permutations of $(\al,\be,\al',\be')$; then
\begin{eqnarray*}
&& \big|\left<\al(i+1),\al'(i+1)\,|\,\be(i+1),\be'(i+1)\right>
 - \left<\al(i),\al'(i)\,|\,\be(i),\be'(i)\right>\big|\\
&& = \big|\left<\al(i+1),\al(i)\,|\,\be(i),\be'(i+1)\right>
   - \left<\be(i+1),\be(i)\,|\,\al'(i),\al(i+1)\right>\\
&& + \left<\al'(i+1),\al'(i)\,|\,\be'(i),\be(i+1)\right>
   - \left<\be'(i+1),\be'(i)\,|\,\al(i),\al'(i+1)\right>\big|\\
&&\le 4L\mu^{i-s-6\de}
\end{eqnarray*}%
\setlength{\unitlength}{1cm}%
\begin{figure}[ht!]
  \begin{center}
   \begin{picture}(4,4)



\put(0,4){\circle*{0.1}} \put(-1.3,3.8){\tiny$\al(i+1)$}
\put(0,0){\circle*{0.1}} \put(-1.3,0){\tiny$\be(i+1)$}
\put(4,4){\circle*{0.1}} \put(4.25,3.8){\tiny$\be'(i+1)$}
\put(4,0){\circle*{0.1}} \put(4.25,0){\tiny$\al'(i+1)$}

\qbezier(0,0)(1,2)(0,4)  \qbezier(4,0)(3,2)(4,4) 
\qbezier(0,0)(2,1)(4,0)  \qbezier(0,4)(2,3)(4,4) 

\put(0.5,3.5){\circle*{0.1}} \put(0.9,2.85){\tiny$\al(i)$}
\put(0.5,0.5){\circle*{0.1}} \put(0.9,1){\tiny$\be(i)$}
\put(3.5,3.5){\circle*{0.1}} \put(2.5,2.85){\tiny$\be'(i)$}
\put(3.5,0.5){\circle*{0.1}} \put(2.5,1){\tiny$\al'(i)$}

\qbezier(0.5,0.5)(1.2,2)(0.5,3.5)  \qbezier(0.5,3.5)(2,2.8)(3.5,3.5) 
\qbezier(3.5,3.5)(2.8,2)(3.5,0.5)  \qbezier(3.5,0.5)(2,1.2)(0.5,0.5) 

\qbezier(0,0)(1.5,1)(3.5,0.5)  \qbezier(4,0)(3,1.5)(3.5,3.5) 
\qbezier(4,4)(2.5,3)(0.5,3.5)  \qbezier(0,4)(1,2.5)(0.5,0.5)

   \end{picture}
  \end{center}
\caption{Extending $\left<\cdot,\cdot|\cdot,\cdot\right>$ from~$X^4$ to~$S$}
\label{f_G4toS}
\end{figure}
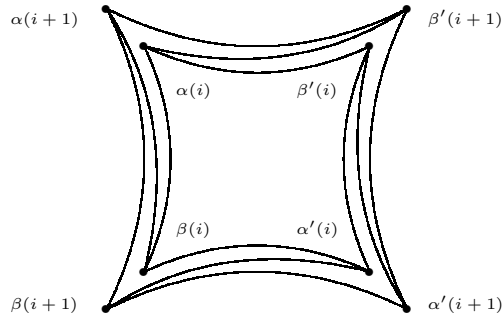%
and
\begin{eqnarray*}
&& \sum_{i=s}^\infty \big|\left<\al(i+1),\al'(i+1)\,|\,\be(i+1),\be'(i+1)\right>
 - \left<\al(i),\al'(i)\,|\,\be(i),\be'(i)\right>\big|\\
&&\le \sum_{i=s}^\infty  4L\mu^{i-s-6\de}= \frac{4L\mu^{-6\de}}{1-\mu}. 
\end{eqnarray*}
This shows that $\left<\al(i),\al'(i)\,|\,\be(i),\be'(i)\right>$ 
is a Cauchy sequence in~$\rr$ so it has a limit in~$\rr$ which we denote
$\left<\bar{a},\bar{a}'\,|\,\bar{b},\bar{b}'\right>$. We show that
$\left<\bar{a},\bar{a}'\,|\,\bar{b},\bar{b}'\right>$ is indeed the limit
in~(\ref{lim_aa'bb'}). For any $i\ge s+16\de$ and any
$(a,a',b,b')\in \big(U^+(\bar{a},i)\times U^+(\bar{a}',i)\times
  U^+(\bar{b},i)\times U^+(\bar{b}',i)\big)\cap (X^{(0)})^4$,
we have by lemmas~\ref{Uconvex} and~\ref{U-U+},
\begin{align*}
d([\be(i),b'],[a,\al(i)])&\ge 
d(U^-(\bar{a},s),U^+(\bar{a},i-8\de))\\
&\ge (i-8\de)-s-6\de= i-s-14\de\ge 2\de
\end{align*}
and the same for the cyclic permutations of $(\al,\al',\be,\be')$.
Then similarly to Figure~\ref{f_G4toS},
\begin{eqnarray*}
&& \left|\left<a,a'\,|\,b,b'\right>
 - \left<\al(i),\al'(i)\,|\,\be(i),\be'(i)\right>\right|\\
&& = \left|\left<a,\al(i)\,|\,\be(i),b'\right>\right.
   - \left<b,\be(i)\,|\,\al'(i),a\right>\\
&& \qquad  + \left<a',\al'(i)\,|\,\be'(i),b\right>
   - \left.\left<b',\be'(i)\,|\,\al(i),a'\right>\right|\\
&&\le 4L\mu^{i-s-14\de}\underset{i\to\infty}{\to} 0.
\end{eqnarray*}
The above inequality is proved for {\em vertices} $a, a', b, b'$, but
by linearity of double difference over simplices~(\ref{dd-ext})
it also holds in the case when
$$(a,a',b,b')\in \big(V^+(\bar{a},i)\times V^+(\bar{a}',i)\times
  V^+(\bar{b},i)\times V^+(\bar{b}',i)\big)\cap X^4.$$
This implies that $\left<\bar{a},\bar{a}'\,|\,\bar{b},\bar{b}'\right>$
is the limit in~(\ref{lim_aa'bb'}).

Now we want to extend $\left<\cdot,\cdot|\cdot,\cdot\right>$ to a continuous function
$\bar{X}^{\diamond}\to\bar\rr$.
Pick an arbitrary $(\bar{a},\bar{a}',\bar{b},\bar{b}')\in\bar{X}^{\diamond}\setminus S$.
This means that there is a trivial side $\p X$--pair in $(\bar{a},\bar{a}',\bar{b},\bar{b}')$,
for example $\bar{a}=\bar{b}'\in\p X$, then $\bar{a}'\not=\bar{a}=\bar{b}'\not=\bar{b}$.
Again, since $S$ is dense in $\bar{X}^{\diamond}$ and $\bar\rr$ is regular, 
it suffices to prove the existence of
\begin{equation}
\label{lim_aa'bb'S}
\lim \left<a,a'|b,b'\right>\quad\mbox{as}\quad 
   (a,a',b,b')\to (\bar{a},\bar{a}',\bar{b},\bar{b}')\ \ \mbox{along}\ \ S
\end{equation}
in~$\bar\rr$.
Fix any $a_0\in X^{(0)}$, then
we have $(\bar{a},\bar{a}',\bar{b},a_0), (a_0,\bar{a}',a_0,\bar{b}')\in S$.
Above we proved that $\left<\cdot,\cdot|\cdot,\cdot\right>$
is continuous on $S$ and takes values in~$\rr$, hence
\begin{eqnarray}\label{aa'ba0}
&&\lim\left<a,a'|b,a_0\right>=\left<\bar{a},\bar{a}'|\bar{b},a_0\right>\in\rr
    \qquad \mbox{and}\\
&& \nonumber \lim\left<a_0,a'|a_0,b'\right>=\left<a_0,\bar{a}'|a_0,\bar{b}'\right>\in\rr\\
\nonumber && \mbox{as}\ (a,a',b,b')\to(\bar{a},\bar{a}',\bar{b},\bar{b}')
\ \mbox{along}\ S.
\end{eqnarray}
\begin{eqnarray}\label{aa'bb'infty}
&&  \left<a,a'|b,b'\right>
  = \left<a,a'|b,a_0\right>+ \left<a,a'|a_0,b'\right>\\ 
\nonumber &&
\qquad = \left<a,a'|b,a_0\right>+\left<a,a_0|a_0,b'\right>+ 
    \left<a_0,a'|a_0,b'\right>
\end{eqnarray}
holds if all the terms are in $X^4$ and therefore, by continuity, in $S$.
Pick any $i\in\rr$ and 
$(a,a',b,b')\in \big(U^+(\bar{a},i)\times U^+(\bar{a}',i)\times
  U^+(\bar{b},i)\times U^+(\bar{b}',i)\big)\cap S$.
Let $b_0\in \mathrm{np}[a,b'|a_0]$. It follows from Lemma~\ref{U-U+} that 
$d(a_0,U^+(\bar{a},i-6\de))\to\infty$ as $i\to\infty$.
Since $d$ and $\dhat$ are quasiisometric (Theorem~\ref{t_dhat}), the same holds for 
$\dhat(a_0,U^+(\bar{a},i-6\de))$.
By Lemma~\ref{Uconvex},
$[a,b']\se [U^+(\bar{a},i),U^+(\bar{a},i)]\se U^+(\bar{a},i-8\de)$,
then by Lemma~\ref{b'a0a0b},
\begin{eqnarray*}
&& \left<a,a_0|a_0,b'\right>=\left<a|b'\right>_{a_0}\ge
\dhat(a_0,b_0)-P\ge  \dhat(a_0,[a,b'])-P\\
&& \ge  \dhat(a_0,U^+(\bar{a},i-6\de))-P\underset{i\to\infty}{\to}\infty.
\end{eqnarray*}
The above inequality holds when $a,a',b,b'$ are vertices, and by~(\ref{dd-ext})
it extends linearly to the case
$(a,a',b,b')\in \big(V^+(\bar{a},i)\times V^+(\bar{a}',i)\times
  V^+(\bar{b},i)\times V^+(\bar{a},i)\big)\cap S$.
(\ref{aa'ba0}) and~(\ref{aa'bb'infty}) also hold in this case, and 
they together imply that the limit in~(\ref{lim_aa'bb'S}) equals~$\infty$.

The same result is obtained when $(\bar{a}',\bar{b})$ is a trivial side $\p X$--pair.
When $(\bar{a},\bar{b})$ or $(\bar{a}',\bar{b}')$ is a trivial side $\p X$--pair, the limit
equals $-\infty$ by the same argument with inequalities reversed.
This implies~(f) and~(g) and the existence of continuous extension to 
a function $\bar{X}^\diamond\to\rb$. 
Parts (a) though~(e) now follow by continuity from
the properties of the double difference in~$X$.

It remains to show (h). Take $a,a',b,b'\in X^{(0)}$,
a geodesic $[b,b']$ in the 1-skeleton  and let $x$ and $x'$
be some nearest point-projections of $a$ and $a'$ to $[b,b']$,
respectively. Orient $[b,b']$ from $b'$ to $b$.
It is an exercise in triangle inequality to see that
$(a,a'|b,b')$ is $^+$equivalent to the signed distance $d(x,x')$,
according to the orientation.
Similarly, using Proposition~\ref{p_uvw},
$\left<a,a'|b,b'\right>$ is $^+$equivalent to the signed distance $\dhat(x,x')$.
But $d(x,x')$ and $\dhat(x,x')$ are $^{\times +}$equivalent.
This proves (h) for vertices in $X$. Now the $^{\times +}$equivalence
extends to all of $\bar{X}^{\diamond}$.
\end{proof}

The double difference is continuous in $\bar{X}^{\diamond}$ and discontinuous at every
quadruple in $\bar{X}^4\setminus \bar{X}^{\diamond}$. Theorem~\ref{dd}
immediately gives the continuous extension of the Gromov product $\left<a|b\right>_c$
to the case when $(a,c,c,b)\in \bar{X}^{\diamond\!}$. This is equivalent to 
$(a,b,c)\in\bar{X}^\triangleright$, where
$$ \bar{X}^\triangleright:=  \{(a,b,c)\in \bar{X}^3\ |\ 
  c\in\p X\ \to\ (a\not=c\ \mbox{and} \ b\not=c)\}.$$
We have $X^3\se \bar{X}^\triangleright\se \bar{X}^3$
and Theorem~\ref{dd} implies

\begin{ttt}\label{t_gr-pr-ext}
If $X$ is a hyperbolic complex, the Gromov product $\left<a|b\right>_c$
with respect to $\dhat$ given by~(\ref{d_gr-pr-new}) extends to a continuous function
$\left<\cdot|\cdot\right>_\cdot: \bar{X}^\triangleright\to [0,\infty]$
such that $\left<a|b\right>_c=\infty$
iff $c\in\p X$ or $a=b\in \p X$.
\end{ttt}

\subsection{More properties of double difference}\label{ss_more-prop-dd}
\begin{lemma}\label{l_dd-dhat}
Let $X$ be a hyperbolic complex.
There exist $A\in [1,\infty)$ and $C\in[0,\infty)$ depending only on~$X$ such that
for all $(a,a',b,b')\in (X^{(0)}\cup\p X)^4\cap\bar{X}^{\diamond}$, 
$\al\in \Geod(a,a')$ and $\be\in \Geod(b,b')$,
$$d(\al,\be)\ge
\max\{\left<b',a|a',b\right>\},\left<b',a'|a,b\right>\}/A-C.$$
\end{lemma}
\begin{proof}
Take a distance-minimizing pair $(a_0,b_0)$, $a_0\in\al$, $b_0\in\be$, and
set notations as in Figure~\ref{f_ab6de}. 
If $a=a'\in \p X$ or $b=b'\in\p X$ then the inequality obviously holds
because both sides are $\infty$. Now we assume otherwise; this implies that
$a_0,b_0\in X^{(0)}$.
By the triangle inequality,
$$\dhat(a_0,b_0)\ge \dhat(a,b)- \dhat(a,a_0)- \dhat(b,b_0)\ \mbox{and}\ 
 \dhat(a_0,b_0)\ge \dhat(a',b')- \dhat(a',a_0)- \dhat(b',b_0).$$
Since $a_0\in\al$ and $b_0\in\be$, by Proposition~\ref{p_uvw},
\begin{eqnarray*}
&& \dhat(a_0,b_0)\\
&& \ge \frac{1}{2} \left(\dhat(a,b)- \dhat(a,a_0)- \dhat(b,b_0)\right)+
     \frac{1}{2} \left(\dhat(a',b')- \dhat(a_0,a')- \dhat(b_0,b')\right)\\
&& =\frac{1}{2}\left(\dhat(a,b)+ \dhat(a',b')\right)-
     \frac{1}{2}\left(\dhat(a,a_0)+ \dhat(a_0,a')\right)-
     \frac{1}{2}\left(\dhat(b,b_0)+ \dhat(b_0,b')\right)\\
&& \ge  \frac{1}{2}\left(\dhat(a,b)+ \dhat(a',b')\right)-
     \frac{1}{2}\left(\dhat(a,a')+C'\right)+
     \frac{1}{2}\left(\dhat(b,b')+C' \right)\\
&& = \left<b',a|a',b\right>-C'.
\end{eqnarray*}
The same argument with $a$ and $a'$ interchanged yields\newline
$\dhat(a_0,b_0)\ge \left<b',a'|a,b\right>-C'$, so
$$\dhat(a_0,b_0)\ge \max\{\left<b',a|a',b\right>,\left<b',a'|a,b\right>\}- C'.$$
Since $\dhat$ and $d$ are quasiisometric (Theorem~\ref{t_dhat}(a)),
\begin{eqnarray*}
&& d(\al,\be)=d(a_0,b_0)\ge \dhat(a_0,b_0)/A-B\\
&& \quad\ge \max\{\left<b',a|a',b\right>,\left<b',a'|a,b\right>\}/A- C'/A- B,
\end{eqnarray*}
so we denote $C:=C'/A+B$.
\end{proof}

\begin{ppp}\label{abba}
Let $X$ be a hyperbolic complex.
There exist constants $T\in [0,\infty)$ and $\la\in [0,1)$ 
depending only on~$X$ such that
for all $(u,a,b,c)\in\bar{X}^{\diamond}$,
if $\left<u,a|b,c\right>\ge T$ or\, $\left<u,b|a,c\right>\ge T$, then
$$\left|\left<u,c|a,b \right>\right|\le
   \la^{\left<u,a|b,c\right>}\le 1\quad\mbox{and}\quad
  \left|\left<u,c|a,b \right>\right|\le \la^{\left<u,b|a,c\right>}\le 1\quad
    \hbox{(see~Figure~\ref{f_vcab1})}.$$
Equivalently, 
$$\mbox{if}\  \max\{\left<u,a|b,c\right>,\left<u,b|a,c\right>\}\ge T\  
\mbox{then}\  \left|\left<u,c|a,b \right>\right|\le 
  \la^{\max\{\left<u,a|b,c\right>, \left<u,b|a,c\right>\}}\le 1.$$
Moreover, $\la$ can be taken arbitrarily close to 1, with $T$ depending on $\la$.
\end{ppp}
\proof
First assume that $(u,a,b,c)\in (X^{(0)}\cup\p X)^4\cap \bar{X}^{\diamond}$.
Take any\newline
 $T\ge A(C+2\de)$, where $A$ and $C$ are from Lemma~\ref{l_dd-dhat}.
Denote $$m:=\max\{\left<u,a|b,c\right>, \left<u,b|a,c\right>\};$$ 
our assumption
is that $m\ge T$. By Lemma~\ref{l_dd-dhat} and our choice of~$T$,
$$d([a,b],[u,c])\ge m/A-C\ge T/A- C\ge 2\de,$$
then by Theorem~\ref{t_dhat}(d),
$$\left|\left<u,c|a,b \right>\right|\le L\mu^{d([a,b],[u,c])}\le 
   L\mu^{m/A-C}=
   (L\mu^{-C})(\mu^{1/A})^{m}.$$
The right hand side decreases exponentially in~$m$,
so by taking $\la\in[0,1)$ sufficiently close to 1 and taking $T$ sufficiently large
we can guarantee that the right hand side is at most
$\la^{m}$ whenever $m\ge T$.
 
\setlength{\unitlength}{1cm} 
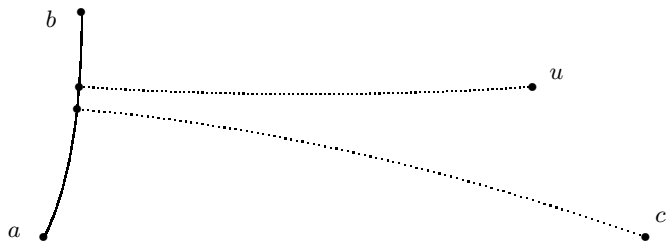
\begin{figure}[ht!]
  \begin{center}
   \begin{picture}(9,3)
\put(0.5,0){\circle*{0.1}}
\put(0,0){\footnotesize$a$}
\put(8.5,0){\circle*{0.1}}
\put(8.6,0.2){\footnotesize$c$}

\put(7,2){\circle*{0.1}}
\put(7.2,2.1){\footnotesize$u$}

\put(1,3){\circle*{0.1}}
\put(0.5,2.8){\footnotesize$b$}
\qbezier(1,3)(1,1)(0.5,0)

\put(0.944,1.7){\circle*{0.1}}
\qbezier[100](0.944,1.7)(4.5,1.4)(8.5,0)

\put(0.973,2){\circle*{0.1}}
\qbezier[100](0.973,2)(4,1.8)(7,2)
   \end{picture}
  \end{center}
\caption{$\left<u,c|a,b\right>$ is exponentially small}
\label{f_vcab1}
\end{figure}

Now consider the general case $(u,a,b,c)\in\bar{X}^\diamond$, and let $m$ be as above.
Let $D\in[0,\infty)$ be the maximal $\dhat$--diameter of a simplex in $X$, then $X^{(0)}\cup\p X$
lies in the D-neighborhood of~$\bar{X}$. If $u\in X$, then $u$ is in a simplex of $X$,
and we replace $u$ with an arbitrary vertex $u'$ of that simplex. If $u\in\p X$, let $u':=u$.
This replacement changes $\left<u,a|b,c\right>$ and $\left<u,b|a,c\right>$
by at most $D$. Doing the same for all four points, i.e.\
replacing $(u,a,b,c)$ with nearby
$(u',a',b',c')\in (X^{(0)}\cup\p X)^4\cap \bar{X}^{\diamond}$ changes
$\left<u,a|b,c\right>$ and $\left<u,b|a,c\right>$ by at most $4D$.
So by the above argument,
\begin{eqnarray*}
&& m \ge T+4D\quad   \Rightarrow\quad 
    \max\{\left<u',a'|b',c'\right>,\left<u',b'|a',c'\right>\}\ge T\\
&& \Rightarrow\quad \left|\left<u',c'|a',b' \right>\right|\le 
    \la^{\max\{\left<u',a'|b',c'\right>, \left<u',b'|a',c'\right>\}}\le 
    \la^{m-4D}= \la^{-4D}\la^m.
\end{eqnarray*}
This holds for all nearby points $(u',a',b',c')$, so by the linearity formula~(\ref{dd-ext}),
the same holds for $(u,a,b,c)$:
$$m \ge T+4D\quad \Rightarrow\quad \left|\left<u,c|a,b \right>\right|\le \la^{-4D}\la^m.$$
Again, since the right hand side decreases exponentially in $m$, we can increase
$\la\in[0,1)$ and $T\in[0,\infty)$ so that
$$m \ge T\quad  \Rightarrow\quad \left|\left<u,c|a,b \right>\right|\le 
\la^m.\eqno{\qed}$$

\section{The cross ratio in $\bar{X}$}\label{s_cross-ratio}
Consider the double difference $\left<\cdot,\cdot|\cdot,\cdot\right>$
given by Theorem~\ref{dd}.
\begin{ddd}
The {\sf  cross-ratio in $\bar{X}$}\label{i_cross-ratio} is the function
$\lbr \cdot, \cdot| \cdot,\cdot\rbr\co {\bar{X}}^{\diamond}\to [0,\infty]$ defined
by
\begin{equation}
\label{cr-def}
\lbr x, x' | y, y' \rbr:= e^{\left<x, x'| y, y'\right>},
\end{equation}
with the convention $e^{-\infty}=0$ and $e^{\infty}=\infty$.
\end{ddd}
Formulas (\ref{dd-def-new}) and (\ref{cr-def}) can be applied to
the standard metric on $\hh^3$ in place of~$\dhat$. In this case
$\lbr \cdot, \cdot | \cdot, \cdot \rbr$ on the boundary 
$\p \hh^3= \ss^2= \cc \cup \{\infty\}$ is the absolute value of the usual 
cross ratio in $\cc \cup \{\infty\}$, therefore the notation.
The following is immediate from Theorem~\ref{dd}.
\begin{ttt}\label{cr}
The cross ratio $\lbr \cdot, \cdot| \cdot,\cdot\rbr$ in $\bar{X}$ defined above
is continuous in ${\bar{X}}^{\diamond}$ and $\Isom(X)$--invariant.
\end{ttt}
This theorem generalizes the fact that M\"obius transformations of $\hh^n$,
and, more generally, isometries of $\mathrm{CAT}(-1)$--spaces, preserve
the cross-ratio on the ideal boundary.
In our case $|\![\cdot, \cdot|\cdot,\cdot]\!|$ is defined on $\bar{X}^\diamond$,
where $X$ is any hyperbolic complex. In particular, it is defined on all pairwise distinct
quadruples of points in $\bar{X}$.
Theorem~\ref{cr} is also a sharp version of~\cite[Proposition~4.5]{Paulin1996} where quasiinvariance
of a (non-continuous) cross-ratio (that is invariance up to an affine function under quasiisometries)
was proved; and also of~\cite{Furman2002} where a measurable
(non-continuous) invariant cross-ratio was constructed.

\section{The symmetric join of $\bar{X}$}\label{s_sj-of-barX}
Let $X$ be a hyperbolic complex. The functor $\sj$ defined
in section~\ref{s_symm-join} can be applied
to any topological space, in particular $\sxb\label{i_sxb}$ makes sense, at least as a set. 
$\sxb$ is called 
the {\sf symmetric join of} \label{i_symm-join-of-Xbar}
 $\bar{X}$.
In this section we extend the earlier constructions from $\sx$ to $\sxb$.

\subsection{Parametrizations of $\sxb$}
\label{ss_param-of-sxb}
Recall that in~\ref{ss_parametrizations} each line $[\![a,b]\!]$ of~$\sxxo$ was
parameterized by the interval
$[-\left<b|x_0\right>_a, \left<a|x_0\right>_b]\se \rr$, where
$x_0$ is a fixed basepoint in~$X$.
Now suppose in addition that $X$ is a hyperbolic complex
with the canonical metric $\dhat$ as in~\ref{ss_dhat}, and $a,b\in \bar{X}$.
Then $(b,x_0,a),(a,x_0,b)\in \bar{X}^\triangleright$,
hence by Theorem~\ref{t_gr-pr-ext},
   $-\left<b|x_0\right>_a$ and $\left<a|x_0\right>_b$ are well-defined elements
of~$\rb$, {\em except} for the case $a=b\in\p X$. In the case $a=b\in\p X$
we let $-\left<b|x_0\right>_a:=\left<a|x_0\right>_b:= 0$,
so that $[\![a,a]\!]$ is identified with the trivial interval $[0,0]\se\rb$.
Note that the function $\left<\cdot|\cdot\right>_\cdot$ is {\em not} continuous at the triples
$(b,x_0,a)$ with $a=b\in\p X$.

This extends the parametrization in~\ref{ss_parametrizations} to a parametrization
$\sxbxo$: lines $[\![a,b]\!]$\label{i_interval2} are identified with the closed intervals 
$[-\left<b|x_0\right>_a, \left<a|x_0\right>_b]$ which are now subintervals
of~$\rb$ rather than of~$\rr$. The lines connecting distinct points at infinity are copies
of~$\rb$. The maps  
$$[\![a,b;\cdot]\!]=[\![a,b;\cdot]\!]_{x_0}\qquad \text{and}\qquad
   [\![a,b;\cdot]\!]'= [\![a,b;\cdot]\!]'_{x_0}$$
are defined by the same formulas as in~\ref{ss_parametrizations}.

The projection function $[\![\cdot,\cdot|\cdot]\!]$ is defined by the same formula
as in Definition~\ref{d_projX}:
\begin{equation}
\left<a,a'|b\right>:= \left<a,a'|b,x_0\right>, \qquad
[\![a,a'|b]\!]:=[\![a,a';\left<a,a'|b\right>]\!].
\end{equation}
But now, by continuity, the projection makes sense for any $a,a',b\in \bar{X}$.

Take any $x_1\in X$. By the same argument as in~\ref{ss_pr-ch-of-basepoint}, 
$[\![a,b;\cdot]\!]_{x_1}$ is the isometric orientation-preserving
reparametrization of $[\![a,b]\!]$ whose origin $[\![a,b;0]\!]_{x_1}$
is the projection of~$x_1$ to~$[\![a,b]\!]$.

\subsection{Actions on $\sx$}\label{ss_actions-on-sxb}
The actions by $\rr$, $\zz_2$, and $\Isom(X)$\label{i_actions-on-sxb}
on $\sx$ extend to $\sxb$ by the same formulas as in~\ref{ss_r-action},
\ref{ss_z2-action}, \ref{ss_isom-action}.
The formula~(\ref{Isom-action}) for the $\Isom(X)$--action indeed makes sense because
by the triangle inequality $\left|\left<a,b|x_0,g\inv x_0\right>\right|\le \dhat(x_0,g\inv x_0)< \infty$
for all $a,b\in X$, and therefore by continuity for all $a,b\in\bar{X}$. For $a=b\in\bar{X}$
the formula implies
$$g\, [\![a,a;0]\!]= [\![ga,ga; 0+\left<a,a|x_0,g\inv x_0\right>]\!]=
    [\![ga,ga; 0]\!],$$
i.e.\ the $\Isom(X)$--action on $\sxb$ restricts to the usual $\Isom(X)$--action
on~$\bar{X}$.

These action satisfy Lemma~\ref{l_actions-prop} with $X$ replaced by $\bar{X}$,
in particular, the $\zz_2$ and $\rr$--actions fix $\bar{X}$
pointwise.

\subsection{The models $\sxb$, $\sxbxo$, $\starxb$, $\starxbxo$, 
$\raise0.12ex\hbox{$\scriptscriptstyle\smallsmile$}
   \mspace{-10.5mu}\raise0.3ex\hbox{$*$} \bar{X}$}
\label{ss_models-for-Xbar}
We use the same notations as in section~\ref{s_symm-join}
with $X$ replaced by $\bar{X}$.
In accordance with \ref{ss_sj-as-top} denote 
$$\starxb:= (\sxb)\setminus \bar{X}.\label{i_starxb}$$
$\starxb$ is called the {\sf open symmetric join of}\label{i_open symm-join-of-Xbar}
$\bar{X}$.

Just as in~\ref{ss_two-models-sx-sxxo} and~\ref{ss_two-models-starx-starxxo}, 
$[\![\cdot,\cdot\,;\cdot]\!]'$
induces a surjection $[\![\cdot,\cdot\,;\cdot]\!]'\co \bar{X}^2\times\rb\to\sxbxo$
and bijections 
\begin{eqnarray}
&&[\![\cdot,\cdot\,;\cdot]\!]'\co \sxb\to 
{
\raise-0.3ex\hbox{$\scriptscriptstyle x_0$}\mspace{-10mu}
   {\raise0.22ex\hbox{$\circ$}\mspace{-9mu}\raise0.22ex\hbox{$*$}
\bar{X}}
},\\
&&]\!]\cdot,\cdot\,;\cdot[\!['\co
 (\bar{X}^2\setminus\bar{\Delta})\times\rr\to\starxbxo\qquad\mbox{and}\qquad
  ]\!]\cdot,\cdot\,;\cdot{[\!['}\co \starxb\to\starxbxo,
\end{eqnarray}
where $\bar{\Delta}$ is the diagonal of $\bar{X}^2$.
The topology on $\sxbxo$ and $\starxbxo$
is defined by either of these maps.

Denote
$$\hsxb:= (\sxb)\setminus \p X=*\bar{X}\cup X,\label{i_hsxb}$$ 
so we have
$*\bar{X}\se \hsxb \se \sxb$.

\subsection{An extension of $\becross$}\label{ss_be-ext}
Pick $x,y\in\hsxb$, then they are of the form $x=[\![a,a';s]\!]$
and $y=[\![b,b';t]\!]$, $a,a',b,b'\in \bar{X}$, for some appropriate $s$ and $t$;
then necessarily $s,t\in\rr$ and
\begin{equation}\label{ainpX}
(a\in\p X\ \to\  a'\not=a)\quad \mbox{and}\quad (b\in\p X\ \to\  b'\not=b).
\end{equation}
Theorem~\ref{dd} and Theorem~\ref{t_gr-pr-ext} imply that the formula
\begin{equation}\label{e_becross-def-xb}
\be_u^{\scriptscriptstyle\times\!}(x,y):=
   \left<a|a'\right>_u + \left|s-\left<a,a'|u\right>\right|
   -\left<b|b'\right>_u - \left|t-\left<b,b'|u\right>\right|
\end{equation}
as in Definition~\ref{d_cocycle-inX} makes sense for all triples
$(u,x,y)\in X\times (\hsxb)^2$, and for such triples $\becross_u(x,y)\in\rr$.
For each fixed $(x,y)$, $\becross_u(x,y)$ is Lipschitz in $u\in X$:
this follows from Theorem~\ref{t_becross}(a) by the continuity of the double
difference, since $\sx$ is dense in $\hsxb$.
We will see later in Theorem~\ref{t_horofunction}
that $\becross$ further extends to a continuous horofunction.

\subsection{The map $\psi=\psi_X$}\label{ss_map-psi}
Recall that $[a,a']$ is an arbitrary fixed choice of geodesic in $X^{(1)}$
connecting vertices $a$ and $a'$. For each $a,a'\in X$ make an arbitrary choice
of projection point in $\mathrm{np}[a,a'|x_0]$
 denoted $[a,a'|x_0]$\label{i_[aa'x0]}.
Proposition~\ref{p_uvw} says that the (images of the) usual geodesics
in $(X^{(1)},d)$ can be parameterized to become
$^+$geodesic in $(X,\dhat)$, i.e.\ geodesic up to a uniform additive constant.

Define a map $\psi=\psi_X\co\sxb\to \bar{X}$ as follows.
\begin{itemize}
\item [\rm(a)] If $x\in \bar{X}$, let $\psi(x):=x$.
\item [\rm(b)] For $a,a'\in \bar{X}$, consider the open interval 
parameterized as in~\ref{ss_param-of-sxb}, so it has the usual metric
as a subinterval of $\rb$. Using Proposition~\ref{p_uvw}
 let $\psi$ map $]\!]a,a'[\![$ $^+$isometrically to $(\,]a,a'[,\dhat\,)$, with a uniform constant.
\item [\rm(c)] If both $a,a'\in\p X$, we additionally require the origins
  $[\![a,a';0]\!]$ to map uniformly close to $[a,a'|x_0]$.
\end{itemize}
The map $\psi$ {\sf $^+$commutes} with the $\Isom(X)$--action on $\bar{X}$ and $\sxb$,
\newline
i.e.\ $\dhat(\psi(gx),g\psi(x))$ is uniformly bounded over all $g\in\Isom(X)$
and $x\in\sxb$. $\psi$~maps $\hsxb$ to $X$.

\subsection{Extending $\dcross$ and $d_*$ to $\sxb$}\label{ss_ext-dcross-dstar}
Define $\dcross$ by the same formula 
\begin{equation*}\label{i_dcross-in-sxb}
\dcross(x,y):=
\sup_{u\in X} |\be_u^{\scriptscriptstyle\times\!}(x,y)|
\end{equation*}
as in~\ref{ss_d^times}, but now applied to all $x,y\in \sxb$, i.e.\
$\dcross\co (\sxb)^2\to[0,\infty]$. Define 
\begin{eqnarray}\label{sj-dhat-phi}
\nonumber 
&&  \sj\dhat(x,y):= \int_{-\infty}^{\infty} \dcross(r^+x,r^+y)\,\frac{e^{-|r|}}{2}\,dr\qquad
  \mbox{and}\\
&&  \varphi(x):= \varphi_X(x):=\int_{-\infty}^\infty r^+ x\, \frac{e^{-|r|}}{2}\, dr,
\end{eqnarray}
and denote for simplicity $d_*:=\sj\dhat$.
The formulas are the same as (\ref{def-dstar})
and~(\ref{d_varphi}) in~\ref{ss_metric-d*}, but 
$d_*\co (\sxb)^2\to[0,\infty]$ and $\varphi\co\sxb\to\sxb$.
In what follows we deal with $\dhat$, $\dcross$ and $d_*$ in the generalized
sense, with infinite values allowed.

\begin{ppp}\label{p_+equiv}
For any hyperbolic complex~$X$, the map $\psi$ in~\ref{ss_map-psi} 
viewed either as
$(\sxb,\dcross)\to (\bar{X},\dhat)$ or as $(\sxb,d_*)\to (\bar{X},\dhat)$,
is a $^+$map.
In particular, $\dcross$ and $d_*$ take finite values on $\hsxb$.
\end{ppp}

\begin{proof}
We can discard the boundary since the values of $\dhat$, $\dcross$ and $d_*$
are either $0$ or $\infty$ and $^+$equivalence is easily checked.
$\psi$ is surjective, 
therefore the conclusion of the proposition is equivalent
to saying that $(\hsxb,\dcross)\to (X,\dhat)$ and $(\hsxb,d_*)\to (X,\dhat)$
are $^+$maps. By Lemma~\ref{l_d*dcross} it suffices to show
that the first map $\psi\co (\hsxb,\dcross)\to (X,\dhat)$ is a $^+$map.

It is not hard to see from the definitions of $\psi$ and $\ell$ (Definition~\ref{d_cocycle-inX},
see also Figure~\ref{f_ell}), using hyperbolicity, Proposition~\ref{p_uvw}
and Lemma~\ref{b'a0a0b},
that $$\ell(u,x)\qquad\mbox{and}\qquad\dhat(u,\psi(x))$$
are $^+$equivalent as functions of $(u,x)\in X\times \hsxb$.
Then by Lemma~\ref{l_plus-eq},
$$|\becross_u(x,y)|=|\ell(u,x)-\ell(u,y)|\qquad\mbox{and}\qquad
  \big|\dhat(u,\psi(x))-\dhat(u,\psi(y))\big|$$
are $^+$equivalent as functions of $(u,x,y)\in X\times(\hsxb)^2$. 
Therefore
$$\dcross(x,y)=\sup_{u\in X} |\becross_u(x,y)|\qquad\mbox{and}\qquad
  \sup_{u\in X} \big|\dhat(u,\psi(x))-\dhat(u,\psi(y))\big|$$
are $^+$equivalent as functions of $(x,y)\in (\hsxb)^2$.
But since $\psi(x),\psi(y)\in X$ and by the triangle inequality,
the last supremum is achieved at $u=\psi(x)$ and it equals
$\dhat(\psi(x),\psi(y))$. So $\dcross(x,y)$ and $\dhat(\psi(x),\psi(y))$
are $^+$equivalent, i.e.\ $\psi$ is a $^+$map.
\end{proof}

A subset $S$ of a (pseudo)metric space $Y$ is {\sf cobounded} in $Y$
if there is $C\in [0,\infty)$ such that $Y$ is contained in the $C$--neighborhood of~$S$.
\begin{ppp}\label{+isom2}
Let $X$ be a hyperbolic complex
and $\psi$ be as in in~\ref{ss_map-psi},
then
$\dcross(y,\psi(y))$ and $d_*(y,\psi(y))$
are bounded uniformly over $y\in\sxb$.
In particular, lines $[\![a,a']\!]$ in $(\hsxb,\dcross)$ and 
in $(\hsxb,d_*)$ are uniformly close to geodesics $[a,a']$
in $X^{(1)}\cup\p X$. Also, 
$X$ is cobounded both in $(\hsxb,\dcross)$ and in $(\hsxb,d_*)$.
\end{ppp}
\begin{proof}
Let $B$ be the constant of the $^+$map $\psi\co (\hsxb,\dcross)\to (X,\dhat)$,
and pick any $y\in\sxb$. By definition $\psi$ is idempotent, i.e.\
$\psi^2(y)=\psi(y)\in \bar{X}$, hence
$$\dcross(y,\psi(y))\le \dhat(\psi(y),\psi^2(y))+B=
\dhat(\psi(y),\psi(y))+B =B.$$
The same proof for $d_*$.
\end{proof}
\begin{ppp}\label{+isom}
For any hyperbolic complex~$X$.
\begin{itemize}
\item [\rm(a)] $d_*$, $\dcross$ and $\dhat$ coincide on~$\bar{X}$,
i.e.\ the canonical embeddings\newline
$(\bar{X},\dhat)\hookrightarrow (\sxb,\dcross)$ and
$(\bar{X},\dhat)\hookrightarrow (\sxb,d_*)$ are isometric.
\item [\rm(b)] The map $\psi$ in~\ref{ss_map-psi} viewed either as
$(\sxb,\dcross)\to (\bar{X},\dhat)$ or as\newline
$(\sxb,d_*)\to (\bar{X},\dhat)$
is a $^+$isometry.
\end{itemize}
\end{ppp}
\begin{proof}
(a) follows from Theorem~\ref{t_surjection}(c)
(or~\ref{t_d^times}(c)).

(b) follows from (a),  propositions~\ref{p_+equiv}
and~\ref{+isom2}.
\end{proof}

\begin{ttt}\label{t_d^times-Xbar}
Let $X$ be a hyperbolic complex.
\begin{itemize}
\item [\rm(a)] The function $\dcross$ above is a well-defined
$\mathrm{Isom}(X)$--invariant pseudometric
on $\hsxb$ independent of $x_0$. It is a generalized pseudometric on $\sxb$.
\item [\rm(b)] The inclusion of each line 
  $([\![a,b]\!],|\cdot|)\hookrightarrow (\sxb,\dcross)$, $a,b\in \bar{X}$,
  is an isometric embedding.
\end{itemize}
\end{ttt}
\begin{proof} By Proposition~\ref{p_+equiv}
the values of $\dcross$ on $(\hsxb)^2$ are indeed in $[0,\infty)$. The rest
follows as in the proof of Theorem~\ref{t_d^times}.
\end{proof}

We summarize various properties that generalize from $X$ to $\bar{X}$.
\begin{ttt}\label{t_extended-d*}
Let $X$ be a hyperbolic complex.
\begin{itemize}
\item [\rm(a)] The function $d_*$ in~(\ref{sj-dhat-phi}) is a well-defined $\mathrm{Isom}(X)$--invariant
  metric on $\hsxb$ independent of~$x_0$. It is a generalized metric on $\sxb$.
\item [\rm(b)] For each $r\in\rr$, the map $r^+\co(\hsxb,d_*)\to(\hsxb,d_*)$ is a bi-Lipschitz
homeomorphism with constant $e^{|r|}$.
\item [\rm(c)] The map $\varphi$ in~\ref{ss_map-psi}  is a well-defined canonical
surjection\newline
$(\sxb,d_*)\twoheadrightarrow (\sxb,\dcross)$
whose restriction to each line $([\![a,a']\!],d_*)$ is an isometry onto
$([\![a,a']\!],\dcross)$. In particular,  each line $[\![a,a']\!]$ in $(\sx,d_*)$
can be parameterized to become a $d_*$--geodesic from $a$ to $a'$.
\item [\rm(d)] The restriction of $\varphi$ to~$\bar{X}$ is the identity
map~$(\bar{X},d_*)\to (\bar{X},\dcross)$, and it is an isometry. 
\end{itemize}
\end{ttt}
\begin{proof}
By Proposition~\ref{p_+equiv}
the values of $d_*$ on $(\hsxb)^2$ are in $[0,\infty)$. The rest
is shown as in theorems~\ref{p_d*}, \ref{bi-lipschitz},
\ref{t_surjection}.
\end{proof}

\section{The topology of $\sxb$}\label{s_top-of-sxb}
\subsection{The topology $\T_*$ on $\sxb$}\label{ss_the-top-on-sxb}
We define a topology $\T_*$\label{i_top-on-sxb} on $\sxb$ as follows.
Neighborhoods of a point $x\in\hsxb$ are
$d_*$--balls centered at $x$. Neighborhoods of a point $x\in\p X$ are
the preimages of neighborhoods of $x$ in $\bar{X}$ under
the $^+$isometry $\psi\co(\sxb,d_*)\to (\bar{X},\dhat)$ from
\ref{ss_map-psi} and Proposition~\ref{+isom}. $\T_*$ induces the original topology on 
$\bar{X}\se\sxb$.

Lemma~\ref{l_aa'-dcross-aa'-dstar}(b) immediately extends to 
\begin{lemma}\label{l_aa'-dcross-aa'-dstar-sxb}
For all $a,a'\in \bar{X}$ with $a\not= a'$, 
$[\![a,a';\cdot]\!]'\co\rb\to ([\![a,a']\!],\T_*)$ is a homeomorphism.
\end{lemma}

\begin{lemma}[Convexity of $\T_*$]\label{l_convexity-of-T*}
Let $X$ be a hyperbolic complex and $b\in\sxb$.
For any neighborhood $N$ of $b$ in $\T_*$ there is a neighborhood $N'$
of $b$ in $\T_*$ with the following properties.
\begin{itemize}
\item [\rm(a)]
If $a,a'\in \bar{X}$ and $x,x'\in[\![a,a']\!]\cap N'$, then
the subinterval of $[\![a,a']\!]$ between $x$ and $x'$ lies in $N$.
\item [\rm(b)]
If $a,a'\in \bar{X}$ and $x,x'\in[\![a,a']\!]\cap (\sxb\setminus N)$, then
the subinterval of $[\![a,a']\!]$ between $x$ and $x'$ lies in $\sxb\setminus N'$.
\end{itemize}
\end{lemma}
\begin{proof}
We will use the property that the lines, hence their subintervals,
in $\sxb$ can be parameterized to become $d_*$--geodesic
(Theorem~\ref{t_extended-d*}(c)).

Assume first that $b\in \hsxb$. (a) follows from the fact that neighborhoods
around $b$ can be taken to be  $d_*$--balls. Suppose (b) does not hold,
then there exist sequences $a_i$ and $a'_i$ in $\bar{X}$, points
$x_i,x'_i\in[\![a_i,a'_i]\!]\cap (\sxb\setminus N)$ and points $y_i$
between $x_i$ and $x'_i$ such that
$d_*(y_i,b)\to 0$.
Then by Lemma~\ref{l_aixor-a'ix}, $d_*(a_i, b)\to 0$ or $d_*(a'_i,b)\to 0$, 
for example the former. Since $x_i$ lie between $a_i$ and $y_i$,
both converging to $b$ in the metric $d_*$, then $x_i$ must also converge to $b$,
which is a contradiction with the choice of $x_i$.

Now assume $b\in\p X$. 
By Proposition~\ref{+isom2},
lines $[\![a,a']\!]$ in $(\hsxb,d_*)$ are uniformly close to geodesics $[a,a']$
in the 1-skeleton.
Then both (a) and (b) follow from lemmas~\ref{U-U+} and~\ref{Uconvex}.
\end{proof}

\subsection{The topology of $\starxb$}
Our goal is to prove the following.
\begin{ppp}\label{p_Xbar-homeomorphic}
Let $X$ be a hyperbolic complex with the standard metric~$\dhat$. 
The metric $d_*$ from~\ref{ss_ext-dcross-dstar}
induces the original topology on the open
symmetric join~$\starxbxo$ described in~\ref{ss_models-for-Xbar}.
Equivalently,
the map from~\ref{ss_models-for-Xbar} viewed as
$$]\!]\cdot,\cdot\,;\cdot[\!['\co
  (\bar{X}^2\setminus\bar{\Delta})\times\rr\to (\starxbxo,d_*)$$
is a homeomorphism,
where $\bar{\Delta}$ is the diagonal of $\bar{X}^2$.
\end{ppp}

\begin{lemma}\label{l_ui-to-b}
Suppose $b,a_i,u_i\in\sxb$, $a_i\to b$,
and $\left<a_i,b|u_i\right>\not\to 0$. Then $b\in \p X$ and
it is possible to take a subsequence so that $u_i\to b$.
\end{lemma}
\begin{proof}
If $b\in X$, then 
$\left<a_i,b|u_i\right>= \left<a_i,b|u_i,x_0\right>\le \dhat(a_i,b)\to 0$
which contradict our assumptions. Therefore $b$ must be in $\p X$.

Now suppose to the contrary that $u_i$ stays away from a neighborhood $V^+$ 
of $b$ for all $i$.
Using the definition of neighborhoods one checks that 
$(a_i,x_0|u_i,b)\to \infty$,
hence by Theorem~\ref{dd}(h), $\left<a_i,x_0|u_i,b\right>\to\infty$,
so by Theorem~\ref{abba},
$|\left<a_i,b|u_i\right>|$ $=|\left<a_i,b|u_i,x_0\right>|\le
  \la^{\left<a_i,x_0|u_i,b\right>}\to 0$,
which is a contradiction.
\end{proof}

\begin{lemma}\label{l_dcross-unif}
Let $X$ be a hyperbolic complex, 
$I\se \rr$ be a compact interval, and $a,a',b,b'\in\bar{X}$.
If $b,b'\in\p X$, assume $b\not= b'$.
Then 
$$\dcross([\![a,a';t]\!],[\![b,b';t]\!])\to 0\quad\mbox{as } 
a\to b\mbox{ and }
a'\to b'\mbox{ in }\bar{X}$$
uniformly on $t\in I$.
\end{lemma}
\begin{proof}
Suppose not, then there is $\epsilon> 0$ and there are sequences
$a_i\to b$, $a'_i\to b'$ in $\bar{X}$ and $t_i\in I$ such that
$$\dcross([\![a_i,a'_i;t_i]\!],[\![b,b';t_i]\!])\ge \epsilon\quad\mbox{for all \ }i.$$
Then by the definition of $\dcross$ there is a sequence $u_i\in X$ such that 
\begin{equation}\label{be-more-epsilon}
\becross_{u_i}([\![a_i,a'_i;t_i]\!],[\![b,b';t_i]\!])\ge \frac{\epsilon}2
\quad\mbox{for all \ }i.
\end{equation}
Taking subsequences we can assume that one of the following
four cases holds.

{\bf Case~1}\quad
{\sl $t_i\le \left<a_i,a'_i|u_i\right>$ and $t_i\le \left< b,b'|u_i\right>$ for all $i$}

By the definition of $\becross$ (\ref{e_becross-def-xb}),
\begin{align}\label{e_aa'uisi}
\nonumber &
\becross_{u_i}([\![a_i,a'_i;t_i]\!],[\![b,b';t_i]\!])=
 \left<a_i|a'_i\right>_{u_i}- t_i+ \left<a_i,a'_i|u_i\right>-
 \left<b|b'\right>_{u_i}+t_i- \left<b,b'|u_i\right>\\
&  =\left<a_i|a'_i\right>_{u_i}+ \left<a_i,a'_i|u_i\right>-
 \left<b|b'\right>_{u_i}- \left<b,b'|u_i\right>.
\end{align}
Denote the last expression $\be_i$. By computation,
$$\be_i= 2\left<a_i,b|u_i\right>- \left<a_i,b|a'_i\right>-
  \left<a'_i,b'|b\right>.$$
The last two terms approach 0 as $i\to\infty$.

We claim that  $\be_i\to 0$. If not then $\left<a_i,b|u_i\right>\not\to 0$,
and by Lemma~\ref{l_ui-to-b},
$b\in\p X$ and, after taking a subsequence,
$u_i\to b$. Then by Theorem~\ref{dd}(g),
$\left<a_i,a'_i|u_i\right>\to -\infty$, so by the assumptions of Case~1,
$t_i\to -\infty$. This is impossible since $t_i\in I$. This proves $\be_i\to 0$.

The condition $\be_i\to 0$ contradicts~(\ref{be-more-epsilon}).

{\bf Case~2}\quad
{\sl $\left<a_i,a'_i|u_i\right>\le t_i$ and $\left< b,b'|u_i\right>\le t_i$ for all $i$}

This is the same as Case 1 by interchanging $a_i\leftrightarrow a'_i$,
$b\leftrightarrow b'$, $t_i\leftrightarrow -t_i$.

{\bf Case~3}\quad
{\sl $\left< b,b'|u_i\right>\le t_i\le \left<a_i,a'_i|u_i\right>$ for all $i$}

By the definition of $\becross$ we have
\begin{eqnarray}\label{e_beaa'ui2}
\nonumber && \becross_{u_i}([\![a_i,a'_i;t_i]\!],[\![b,b';t_i]\!])\\
&& \qquad = \left<a_i|a'_i\right>_{u_i}- t_i+ \left<a_i,a'_i|u_i\right>-
 \left<b|b'\right>_{u_i}- t_i+ \left<b,b'|u_i\right>.
\end{eqnarray}
First we let $t_i:= \left< b,b'|u_i\right>$, then
\begin{equation*}
\becross_{u_i}([\![a_i,a'_i;t_i]\!],[\![b,b';t_i]\!])=
 \left<a_i|a'_i\right>_{u_i}- \left< b,b'|u_i\right>+ \left<a_i,a'_i|u_i\right>
 - \left<b|b'\right>_{u_i}.
\end{equation*}
This expression is the same as in~(\ref{e_aa'uisi}), so we call it $\be_i$.
The same argument as in Case~1 shows that $\be_i\to 0$.

Now let $t_i:= \left< a_i,a_i'|u_i\right>$, then
\begin{equation*}
\becross_{u_i}([\![a_i,a'_i;t_i]\!],[\![b,b';t_i]\!])=
 \left<a_i|a'_i\right>_{u_i}- \left< b|b'\right>_{u_i}- \left<a_i,a'_i|u_i\right>
 + \left<b,b'|u_i\right>.
\end{equation*}
We call the last expression $\be'_i$.
It is obtained from $\be_i$ by interchanging
$a_i\leftrightarrow a'_i$, $b\leftrightarrow b'$, and 
a similar argument shows that $\be'_i\to 0$.

By the assumptions of Case~3, (\ref{e_beaa'ui2}) lies between $\be_i$ and $\be'_i$, 
hence it converges to 0. This contradicts~(\ref{be-more-epsilon}).

{\bf Case~4}\quad
{\sl $\left< a_i,a'_i|u_i\right>\le t_i\le \left<b,b'|u_i\right>$ for all $i$}

This is the same as Case 3 by interchanging $a_i\leftrightarrow a'_i$,
$b\leftrightarrow b'$, $t_i\leftrightarrow -t_i$.
\end{proof}

\begin{lemma}\label{l_dcross'-unif}
Under the assumptions of Lemma~\ref{l_dcross-unif},
$$\dcross([\![a,a';t]\!]',[\![b,b';t]\!]')\to 0\quad\mbox{ as }
a\to b\mbox{ and } a'\to b'\mbox{ in }\bar{X}$$
uniformly on $t\in I$.
\end{lemma}
\begin{proof}
Recall from~\ref{ss_param-of-sxb} that $[\![a,a']\!]$ is a copy of the interval
$[\al,\al']\se\rb$,
where $\al:= -\left<a'|x_0\right>_a$ and $\al':= \left<a|x_0\right>_{a'}$.
Similarly, $[\![b,b']\!]$ is a copy of $[\be,\be']\se\rr$ where
$\be:= -\left<b'|x_0\right>_b$ and $\be':= \left<b|x_0\right>_{b'}$.
If $a\to b$ and $a'\to b'$ in $\bar{X}$, then $\al\to\be$ and $\al'\to\be'$ in~$\rb$.
(When some of $a,a',b,b'$ are in $\p X$, this is an exercise for the Gromov product
defined by the word metric $d$, and then use Theorem~\ref{dd}(h) to show the same
for $\left<\cdot|\cdot\right>_\cdot$.)

Using Lemma~\ref{theta'theta} we denote
\begin{eqnarray*}
&& A:= \theta'[\al,\al';t]= \theta[\al,\al';t]+ (e^{-|t-\al|}- e^{-|t-\al'|})/2\quad\mbox{and}\\
&& B:= \theta'[\be,\be';t]= \theta[\be,\be';t]+ (e^{-|t-\be|}- e^{-|t-\be'|})/2,
\end{eqnarray*}
then
$$|B-A|\le \big|\theta[\be,\be';t]- \theta[\al,\al';t]\big|+
 \big|e^{-|t-\be|}- e^{-|t-\al|}\big|/2+ \big|e^{-|t-\al'|}- e^{-|t-\be'|}\big|/2.$$
This implies that 
for each compact $J\se\rr$, 
\begin{equation}\label{e_B-A}
|B- A|\to 0 \quad\mbox{ as }
a\to b\mbox{ and } a'\to b'\mbox{ in }\bar{X}
\end{equation}
uniformly on $t\in J$.
By~(\ref{abt'}) and since the map
$[\![a,a';\cdot]\!]'\co\rr\to ([\![a,a']\!],\dcross)$ is non-expanding,
\begin{align*}
& \dcross\big([\![a,a';t]\!]', [\![b,b';t]\!]'\big)=
   \dcross\big([\![a,a';\theta'[\al,\al';t]]\!], [\![b,b';\theta'[\be,\be';t]]\!]\big)\\
& = \dcross\big([\![a,a';A]\!],[\![b,b';B]\!]\big)\le
  \dcross\big([\![a,a';A]\!],[\![a,a';B]\!]\big)+
  \dcross\big([\![a,a';B]\!],[\![b,b';B]\!]\big)\\
& \le |B-A|+\dcross\big([\![a,a';B]\!],[\![b,b';B]\!]\big).
\end{align*}
This, Lemma~\ref{l_dcross-unif} and~(\ref{e_B-A}) imply Lemma~\ref{l_dcross'-unif}.
\end{proof}

\begin{ppp}\label{l_dstar-to-0}
Let $X$ be a hyperbolic complex, $s,t\in\rb$, $a,a',b,b'\in\bar{X}$.
Then 
$$[\![a,a';s]\!]'\to [\![b,b';t]\!]'\quad\mbox{in }(\sxb,\T_*)\quad\mbox{as}\quad 
a\to b,\ a'\to b'\mbox{ in }\bar{X}\mbox{ and }s\to t\mbox{ in }\rb.$$
Equivalently, the map 
$[\![\cdot,\cdot\,;\cdot]\!]'\co\bar{X}^2\times\rb\to (\sxb,\T_*)$
is continuous.
\end{ppp}
\begin{proof}
First assume $b=b'$. Since $a\to b$ and $a'\to b'=b$, then 
by the convexity of $\T_*$  (Lemma~\ref{l_convexity-of-T*}),
all the points of $[\![a,a']\!]$ must converge to $b$ as well.
So now we assume $b\not= b'$.

Assume $t\in\rr$. For any $s\in\rr$,
since $[\![a,a';\cdot]\!]'\co\rr\to ([\![a,a']\!],\dcross)$ is non-expanding, 
\begin{eqnarray*}
&&\dcross \big([\![a,a';r+s]\!]',[\![b,b';r+t]\!]'\big)\\
&& \le
\dcross \big([\![a,a';r+s]\!]',[\![a,a';r+t]\!]'\big)+
\dcross \big([\![a,a';r+t]\!]',[\![b,b';r+t]\!]'\big)\\
&& \le |t-s|+ \dcross \big([\![a,a';r+t]\!]',[\![b,b';r+t]\!]'\big).
\end{eqnarray*}
Pick any $\epsilon >0$. By Lemma~\ref{l_dcross'-unif} there exist
neighborhoods $N_b$ of $b$ and $N_{b'}$ of $b'$ in $\bar{X}$ 
such that
$$\dcross \big([\![a,a';r+t]\!]',[\![b,b';r+t]\!]'\big)\le \epsilon$$
for all $a\in N_b$, $a'\in N_{b'}$ and
$r\in [-\frac{1}{\epsilon}, \frac{1}{\epsilon}]$.
Also for $r\in (-\infty, -\frac{1}{\epsilon}]\cup [\frac{1}{\epsilon},\infty)$,
we have
\begin{eqnarray*}
&& \dcross \big([\![a,a';r+t]\!]',[\![b,b';r+t]\!]'\big)\le
\dcross \big([\![a,a';r+t]\!]',[\![a,a';t]\!]'\big)+\\
&& \quad\dcross \big([\![a,a';t]\!]',[\![b,b';t]\!]'\big)+
\dcross \big([\![b,b';t]\!]',[\![b,b';r+t]\!]'\big)\le 2|r|+\epsilon.
\end{eqnarray*}
Then by the definition of~$d_*$,
\begin{eqnarray*}
&& d_* ([\![a,a';s]\!]',[\![b,b';t]\!]')=
\int_{-\infty}^{\infty} \dcross \big([\![a,a';r+s]\!]',[\![b,b';r+t]\!]'\big)
   \frac{e^{-|r|}}{2}\, dr\\
&&\le\int_{-\frac{1}{\epsilon}}^{\frac{1}{\epsilon}} 
   (|t-s|+ \epsilon) \frac{e^{-|r|}}{2}\, dr+ 
  \int_{\frac{1}{\epsilon}}^{\infty} 
    (2|r|+\epsilon)\frac{e^{-|r|}}{2}\, dr+
  \int_{-\infty}^{-\frac{1}{\epsilon}} 
    (2|r|+\epsilon)\frac{e^{-|r|}}{2}\, dr
\\
&&\le |t-s|+\epsilon+ 2(1/\epsilon+ 1)e^{-1/\epsilon}\to 0\quad\mbox{as }
s\to t \mbox{ and } \epsilon\searrow 0.
\end{eqnarray*}
This implies the statement of the proposition.

Assume $t=\infty$, then $[\![b,b';t]\!]'=[\![b,b';\infty]\!]'=b'\in\bar{X}$.
By Lemma~\ref{l_aa'-dcross-aa'-dstar-sxb},
$[\![b,b';\cdot]\!]'\colon$\break $\rb\to ([\![b,b']\!],\T_*)$ is a homeomorphism,
therefore by taking $s\in\rr$ close to $\infty$ we can make the whole interval
$[\![b,b';[s,\infty]]\!]'$ arbitrarily close to $b'$.
For a fixed $s\in\rr$, by the above argument, by taking $(a,a')$ close to $(b,b')$
in $T_*$ we can make $d_*([\![a,a';s]\!]',[\![b,b';s]\!]')$ arbitrarily small.
This implies that both $[\![a,a';s]\!]'$ and $a'=[\![a,a';\infty]\!]'$ are
arbitrarily close to $b'$, then by the convexity of $T_*$, the interval
$[\![a,a';[s,\infty]]\!]'$ can be made arbitrarily close to~$b'$.
The case $t=-\infty$ is similar.
\end{proof}

\begin{lemma}\label{l_aixor-a'ix-sxb}
Let $X$ be a hyperbolic complex.
Let $a\in \bar{X}$ and\newline
 $x_i=[\![a_i,a'_i;s_i]\!]'\in\sxb$ be a sequence 
such that $x_i\to a$ in $(\sxb,\T_*)$.
Then $a_i\to a$ or $a'_i\to a$ in~$\bar{X}$.
\end{lemma}
\begin{proof}
First consider the case $a\in X$, so our assumption is that
$d_*(x_i,a)\to 0$. The proof is  word-by-word as in Lemma~\ref{l_aixor-a'ix}.

Now we assume $a\in\p X$. 
If, to the contrary, some subsequences $a_j$ and $a'_j$ stay away from a neighborhood
$V$ of $a$ in $\bar{X}$, then 
by the convexity of $\T_*$
all the lines $[\![a_j,a'_j]\!]$, and hence all $x_j$,
must stay away from some smaller neighborhood of $a$, which contradicts $x_j\to a$.
\end{proof}

\begin{proof}[Proof of Proposition~\ref{p_Xbar-homeomorphic}]
We will show that the identity map\newline
 $\starxbxo\to(\starxbxo,d_*)$
and its inverse are continuous at every point $y=]\!]b,b';t[\!['$ 
in~$\starxbxo$.

Proposition~\ref{l_dstar-to-0} says that 
$\starxbxo\to(\starxbxo,d_*)$ is continuous at $y$.
Suppose that the inverse map $(\starxbxo,d_*)\to \starxbxo$
is not continuous at $y$.
We obtain a contradiction just as in the proof of 
Proposition~\ref{p_Xhomeomorphic}, but using $\bar{X}$,
Proposition~\ref{l_dstar-to-0} and Lemma~\ref{l_aixor-a'ix-sxb}.
\end{proof}

\subsection{Properness of $
{\raise0.12ex\hbox{$\scriptscriptstyle\smallsmile$}
   \mspace{-11.6mu}\raise0.25ex\hbox{$*$}
\bar{X}}$}\label{ss_proper}
A metric space $Y$ is called {\sf proper}\label{i_proper}
 if each closed ball in $Y$ is compact.
Recall that $\hsxb=\sxb\setminus\p X$. 
\begin{ppp}
\label{p_proper}
For any hyperbolic complex~$X$,
$(\hsxb,d_*)$ is proper.
\end{ppp}
\begin{proof}
Let $B_{d_*}\!(r)$  be the closed ball in $(\hsxb,d_*)$ of radius $r$
centered at the basepoint $x_0\in X$.
Since $X$ is cobounded in $(\hsxb,d_*)$, any ball lies 
in $B_{d_*}\!(r)$ for sufficiently large~$r$. 
Since $(\hsxb,d_*)$ is a metric space, it suffices to show the sequential
compactness of $B_{d_*}\!(r)$. 

Fix $r\ge 0$ and pick
$x=[\![a,a';s]\!]'\in B_{d_*}\!(r)$.
We have by the definitions of $\ell$ and $\dcross$ and Lemma~\ref{l_d*dcross},
\begin{eqnarray*}
&&\left<a|a'\right>_{x_0}\le \left<a|a'\right>_{x_0}+\left|s-\left<a,a'|x_0\right>\right| =
  \ell(x_0,x)= |\ell(x_0,x_0)-\ell(x_0,x)|\\
&& \le \sup_{u\in X}|\ell(u,x_0)-\ell(u,x)|
  =\dcross(x_0,x)\le d_*(x_0,x)+2\le r+2.
\end{eqnarray*}
Similarly $|s|=\left|s-\left<a,a'|x_0\right>\right|\le r+2$.
Let $x_i=[\![a_i,a'_i;s_i]\!]$ be  a sequence in $B_{d_*}\!(r)$.
We have 
$$\left<a_i|a'_i\right>_{x_0}\le r+2\quad\mbox{and}\quad |s_i|\le r+2
   \quad\mbox{for all }i.$$
After replacing $\{x_i\}$ with a subsequence, $s_i$ converges to some $\bar{s}$
with $|\bar{s}|\le r+2$. 
By Theorem~\ref{t_gr-pr-ext}, the function
$\left<\cdot|\cdot\right>_{x_0}\co\bar{X}^2\to[0,\infty]$ is continuous, 
hence the set
$$\{(a,a')\in\bar{X}^2\ |\ \left<a|a'\right>_{x_0}\le r+2\}$$
is closed in $\bar{X}^2$ and therefore 
compact, so after replacing $\{x_i\}$ with a subsequence again,
$a_i$ and $a'_i$ converge
to some $\bar{a}$ and $\bar{a}'$ in $\bar{X}$, respectively, satisfying
$\left<\bar{a}|\bar{a}'\right>_{x_0}\le r+2$. 

Denote $\bar{x}:=[\![\bar{a},\bar{a}';\bar{s}]\!]$.
By Proposition~\ref{l_dstar-to-0}, $d_*(x_i,\bar{x})\to 0$. Also
$$d_*(x_0,\bar{x})\le d_*(x_0,x_i)+ d_*(x_i,\bar{x})\le
r+ d_*(x_i,\bar{x})\ \underset{i\to\infty}\to \ r,$$
so $\bar{x}\in B_{d_*}\!(r)$.
\end{proof}

\section{Horofunctions and horospheres
  in~$\bar{X}$ and $\sxb$}\label{s_horofunctions-in-barX}
In the classical case of $\hh^n$, a horofunction, or Busemann function, $\be_u$ is defined
with respect to a point at infinity, $u\in \p\hh^n$.
A horofunction in a hyperbolic metric space
related to a point $u$ on the ideal boundary was usually defined either depending
on a geodesic ray converging to~$u$ (\cite[III.3.4]{BH}, \cite[7.5.C]{Gr2}),
or satisfied the natural identities only ``up to a constant" (\cite[Chapitre 8]{GH}),
or was a measurable function (\cite{Furman2002}). 

In this section we show that the metric $\dhat$ allows defining a continuous
horofunction $\be_u^{\times}$ on any hyperbolic complex~$X$;
the horofunction will satisfy sharp identities as in the classical $\hh^n$ case.
We will also allow $u\in X$; in this case $\be^{\times}_u(x,y)$
will be the {\sf distance cocycle}:
$\be^{\times}_u(x,y)= \dhat(u,x)- \dhat(u,y)$.
Moreover, $\be^{\times}_u(x,y)$ will be defined not only for $x,y\in X$
but also for $x,y\in \starx$, and more.

Let $X$ be a hyperbolic complex. Recall that
$\hsxb= (\sxb)\setminus \p X=\starxb\cup X$.

\begin{ttt}\label{t_horofunction}
Let $X$ be a hyperbolic complex and $\becross$ be the cocycle
from Definition~\ref{d_cocycle-inX}.
\begin{itemize}
\item [\rm(a)] Put the pseudometric $\dhat+\dcross+\dcross$ on~$X\times (\hsxb)^2$.
Then $\becross$ extends to a Lipschitz function $\becross\co X\times (\hsxb)^2\to\rr$
  independent of~$x_0$.
\item [\rm(a$'$)] Put the usual topology on $\bar{X}$ and the topology induced by~$\dcross$
  on $\hsxb$. Then $\becross$ further extends to a continuous function
  $\becross\co \bar{X}\times (\hsxb)^2\to\rr$ independent of~$x_0$.
\item [\rm(b)] $\be_u^{\scriptscriptstyle\times\!}$ is $\zz_2$--invariant in each variable:
$\be_u^{\scriptscriptstyle\times\!}(x,y)=\be_u^{\scriptscriptstyle\times\!}(x^\star,y)=\be_u^{\scriptscriptstyle\times\!}(x,y^{\star})$.
\item [\rm(c)] $\be_u^{\scriptscriptstyle\times\!}$ satisfies the cocycle condition \ 
$\be_u^{\scriptscriptstyle\times\!}(x,y)+ \be_u^{\scriptscriptstyle\times\!}(y,z)=\be_u^{\scriptscriptstyle\times\!}(x,z)$. 
\item [\rm(d)] $\becross$ is $\Isom(X)$--invariant:
  $\be_{gu}^{\scriptscriptstyle\times\!}(gx,gy)=\be_u^{\scriptscriptstyle\times\!}(x,y)$ for $g\in\Isom(X)$.
\item [\rm(e)] $\becross$ is isometric on lines: for all $a,b\in\bar{X}$ and
  $x,y\in [\![a,b]\!]\setminus\p X$,\newline
  $|\becross_a(x,y)|=|\becross_b(x,y)|=\dcross(x,y)$.
\end{itemize}
\end{ttt}
\noindent This extension $\becross\co \bar{X}\times (\hsxb)^2\to\rr$
will be called the {\sf horofunction in~$\sxb$}.\label{i_horofunction-in-sxb}
\begin{proof}
We only need to prove (a) and (a$'$), then
(b)-(d) follows from the properties of the original cocycle and (e) follows from
Lemma~\ref{l_beuxy}(a) and Theorem~\ref{t_d^times-Xbar}(b).

\noindent {\bf (a)}\qua
We have $x=[\![a,a';s]\!]$ and $y=[\![b,b';t]\!]$ as in~\ref{ss_be-ext}.
We saw in~\ref{ss_be-ext} that the formula
\begin{equation}\label{be-again}
\be_u^{\scriptscriptstyle\times\!}(x,y):=
   \left<a|a'\right>_u + \left|s-\left<a,a'|u\right>\right|
   -\left<b|b'\right>_u - \left|t-\left<b,b'|u\right>\right|
\end{equation}
makes sense for the triples $(u,x,y)\in X\times (\hsxb)^2$, and that 
for such triples $\becross_u(x,y)\in\rr$ and is Lipschitz in $u$.
The inequality 
\begin{eqnarray*}
&& \left|\becross_u(x',y')-\becross_u(x,y)\right|\le
 \left|\becross_u(x',x)\right|+\left|\becross_u(y',y)\right|\\
&& \le \sup_{u\in X}\left|\becross_u(x',x)\right|+\sup_{u\in X}\left|\becross_u(y',y)\right|=
   \dcross(x',x)+\dcross(y,y')\\
\end{eqnarray*}
shows that $\becross$ is also Lipschitz in $x$ and $y$, therefore in the three variables
simultaneously.

\noindent {\bf (a$'$)}\qua
We want to extend~$\becross$ to $\bar{X}\times(\hsxb)^2$ continuously.
Let $\bar{u}\in\p X$.
By~\cite[I \S 8 \No5, Theorem~1]{Bourbaki} it suffices to show the existence
of the limit
\begin{equation}\label{lim-be}
\lim \be_{u}^{\scriptscriptstyle\times\!}(x,y)\quad \mbox{as}
\ \ u\to \bar{u}\ \ \mbox{along}\ \,X
\end{equation}
in~$\rr$. First assume $\bar{u}\in{\p X}\setminus\{a,a',b,b'\}$, then
by Theorem~\ref{dd} and Theorem~\ref{t_gr-pr-ext}, as $u\to\bar{u}$ along~$X$,
\begin{eqnarray*}
&& \left<a,a'|u\right>\to \left<a,a'|\bar{u}\right>\in\rr,\qquad
\left<b,b'|u\right>\to\left<b,b'|\bar{u}\right>\in\rr,\\
&&\left<a|a'\right>_u- \left<b|b'\right>_u=
  \left<a|a'\right>_{x_0}- \left<b|b'\right>_{x_0}+
  \left<a,b|u,x_0\right>+ \left<a',b'|u,x_0\right>\\
&& \qquad\to   \left<a|a'\right>_{x_0}- \left<b|b'\right>_{x_0}+
  \left<a,b|\bar{u},x_0\right>+ \left<a',b'|\bar{u},x_0\right>\in\rr,
\end{eqnarray*}
because under our assumptions there are no trivial $\p X$--triples in the above terms.
Hence by~(\ref{be-again}) the limit~(\ref{lim-be}) exists in $\rr$.

Now assume $\bar{u}=a\in\p X\setminus\{a',b,b'\}$, then as $u\to a$ along~$X$,
\begin{equation}\label{bb'u}
\left<b,b'|u\right>\to \left<b,b'|a\right>\in\rr\qquad\mbox{and}\qquad
\left<a,a'|u\right>\to-\infty,
\end{equation}
hence for $u$ sufficiently close to~$a$, 
\begin{eqnarray}\label{aa'us-aa'}
&& \left<a|a'\right>_u + \left|s-\left<a,a'|u\right>\right|
   -\left<b|b'\right>_u = \left<a|a'\right>_u + s-\left<a,a'|u\right>
   -\left<b|b'\right>_u\\
\nonumber && = s+ \left<a',b'|u,x_0\right>+ \left<a',b|u,x_0\right>\ \to\ 
     s+ \left<a',b'|a,x_0\right>+ \left<a',b|a,x_0\right>\in\rr.
\end{eqnarray}
(\ref{bb'u}) and~(\ref{aa'us-aa'}) show that 
the limit~(\ref{lim-be}) exists in $\rr$. The similar argument works
for each of the cases $\bar{u}=a'\in\p X\setminus\{a,b,b'\}$, 
$\bar{u}=b\in\p X\setminus\{a,a',b'\}$, $\bar{u}=b'\in\p X\setminus\{a,a',b\}$.

Now assume $\bar{u}=a=b\in\p X\setminus\{a',b'\}$, then as $u\to a$ along~$X$,
$$\left<a,a'|u\right>\to -\infty\qquad\mbox{and}\qquad\left<b,b'|u\right>\to -\infty,$$
hence for $u$ sufficiently close to~$a$, 
\begin{eqnarray*}
&&\be_u^{\scriptscriptstyle\times\!}(x,y)=
   \left<a|a'\right>_u + \left|s-\left<a,a'|u\right>\right|
   -\left<b|b'\right>_u - \left|t-\left<b,b'|u\right>\right|\\
&& =\left<a|a'\right>_u + s-\left<a,a'|u\right>
   -\left<b|b'\right>_u - t+\left<b,b'|u\right>\\
&& =\left<a',b'|u,a\right>+ \left<a',b'|u,x_0\right>\ \to
   \left<a',b'|a,x_0\right>\in\rr.
\end{eqnarray*}
The cases $\bar{u}=a=b'\in\p X\setminus\{a',b\}$, 
$\bar{u}=a'=b\in\p X\setminus\{a,b'\}$, $\bar{u}=a'=b'\in\p X\setminus\{a,b\}$
are similar. Our condition~(\ref{ainpX}) implies that
there are no more cases to consider.
\end{proof}

Of interest are special cases of Theorem~\ref{t_horofunction}:
\begin{itemize}
\item When $u\in X$, $x\in[\![a,u]\!]$ and $y\in[\![b,u]\!]$, 
$\be_u^{\scriptscriptstyle\times\!}$ becomes the usual distance cocycle:
$\be_u^{\scriptscriptstyle\times\!}(x,y)=\dcross(u,y)-\dcross(u,x)$.
In particular, $\be_u^{\scriptscriptstyle\times\!}(x,y)= \dhat(u,y)- \dhat(u,x)$ for \mbox{$u,x,y\in X$}.
\item When $u\in\p X$, $\be_u^{\scriptscriptstyle\times\!}$ becomes 
a horofunction:\newline
$\be_u^{\scriptscriptstyle\times\!}(x,y)=
\lim \big(\dcross(v,y)-\dcross(v,x)\big)$ as $v\to u$ along $X$.
In particular, for $x,y\in X$, $\be_u^{\scriptscriptstyle\times\!}(x,y)=\lim \big(\dhat(v,y)- \dhat(v,x)\big)$
as $v\to u$ along $X$. Thus $\becross$ satisfies the usual definition of a horofunction;
the limits indeed exist.
\end{itemize}

For $x\in\hsxb$, $\mathcal{H}_u(x):= \{y\in\hsxb\ |\ \be_u^{\scriptscriptstyle\times\!}(x,y)=0\}$
is the {\sf horosphere}\label{i_horosphere} at $u$ containing~$x$.

\begin{lemma}
\label{l_horosphere}
Let $u, a, b\in \bar{X}$ and suppose that those of $u,a,b$ that lie in $\p X$ are
pairwise distinct. Then the four projections
$[\![a,u|b]\!]$, $[\![b,u|a]\!]$, $[\![u,a|b]\!]$ and $[\![u,b|a]\!]$
lie on the same horosphere at~$u$.
\end{lemma}
\begin{proof}
The assumptions on $u,a,b$ and Theorem~\ref{dd}(fg) 
imply that $\left<a,u|b\right>$, $\left<b,u|a\right>$, $\left<u,a|b\right>$,
$\left<u,b|a\right>$ are in~$\rr$, therefore 
all the four projections lie in $\hsxb$ and
the horospheres containing them are indeed well-defined.
If $u\in X$, using Definition~\ref{d_projX},
\begin{eqnarray*}
&& \be_u^{\scriptscriptstyle\times\!}\big([\![a,u|b]\!],[\![b,u|a]\!]\big)=
   \be_u^{\scriptscriptstyle\times\!}\big([\![a,u;\left<a,u|b\right>]\!],[\![b,u;\left<b,u|a\right>]\!]\big)\\
&& = \left<b|u\right>_u+ |\left<b,u|u\right>- \left<b,u|a\right>|-
     \left<a|u\right>_u- |\left<a,u|u\right>-\left<a,u|b\right>|\\
&& = |\left<b,u|u,a\right>|- |\left<a,u|u,b\right>|=
     \left<b|a\right>_u- \left<a|b\right>_u=0,
\end{eqnarray*}
i.e.\ $[\![a,u|b]\!]$ and $[\![b,u|a]\!]$
lie on the same horosphere at $u$. This extends by continuity
(Theorem~\ref{t_horofunction}(a)) to the case $u\in\bar{X}$.
The rest of lemma follows from the $\zz_2$--invariance of~$\becross$:
$$\be_u^{\scriptscriptstyle\times\!}([\![a,u|b]\!], [\![u,a|b]\!])= \be_u^{\scriptscriptstyle\times\!}([\![a,u|b]\!], [\![a,u|b{]\!]}^\star)= 
\be_u^{\scriptscriptstyle\times\!}([\![a,u|b]\!], [\![a,u|b]\!])=0,$$
and similarly for $[\![u,b|a]\!]$.
\end{proof}

\section{Synchronous exponential convergence of\newline
lines in~$\sxb$}\label{s_synchronous}
Recall from~\ref{ss_pr-ch-of-basepoint} and~\ref{ss_param-of-sxb} that 
$[\![b,c;\cdot]\!]_{a}$ is the isometric reparametrization of $[\![b,c]\!]$
whose origin $[\![b,c;0]\!]_{a}$ is the projection of~$a$ to $[\![b,c]\!]$.
\begin{ttt}[Exponential convergence in $\dcross$]\label{exp-conv}
Let $X$ be a hyperbolic complex.
There exist $N\in[0,\infty)$ and $\la\in (1/e,1)$ depending only on~$X$
such that for all~$t\in\rr$ and $a,b,c\in \bar{X}$,
$$\dcross\big([\![b,c;t]\!]_{a}, [\![a,c;t]\!]_{b}\big)\le N\la^t.$$
\end{ttt}
\noindent This theorem provides 
\begin{itemize}
\item an upper exponential bound that is independent of the choice
of $a,b,c$;
\item a {\em synchronous} exponential convergence: it is easy to see
  from Lemma~\ref{l_horosphere} that at each time~$t\ge 0$,
$[\![b,c;t]\!]_{a}$ and $[\![a,c;t]\!]_{b}$ lie on the same horosphere centered at~$c$;
\item the place where
the exponential convergence starts occurring, namely the projections
$[\![b,c|a]\!]$ and $[\![a,c|b]\!]$ corresponding to~$t=0$.
\end{itemize}
Note also that~$t$ is not assumed to be appropriate for $[\![b,c]\!]_{a}$ or $[\![a,c]\!]_{b}$.
\setlength{\unitlength}{1cm} 
\begin{figure}[ht!]
  \begin{center}
   \begin{picture}(12,4)
\put(0.5,0){\circle*{0.1}}
\put(0.2,0.2){\footnotesize$a$}
\put(11.5,0){\circle*{0.1}}
\put(11.6,0.2){\footnotesize$c$}
\qbezier(0.5,0)(6,0)(11.5,0)

\put(6.45,0.46){\circle*{0.1}}
\put(6.1,0.7){\tiny$[\![b,c;t]\!]_{a}$}

\put(6.4,0){\circle*{0.1}}
\put(5.9,-0.35){\tiny$[\![a,c;t]\!]_{b}$}

\put(3.8,1.68){\circle*{0.1}}
\put(3.7,1.95){\tiny$[\![b,c|a]\!]= [\![b,c;0]\!]_a$}
\qbezier[40](0.5,0)(2.5,0)(3.8,1.68)

\put(8.95,2.1){\footnotesize$u$}
\put(8.8,2){\circle*{0.1}}

\put(8.8,0.11){\circle*{0.1}}
\put(8.9,0.3){\tiny$[\![b,c|u]\!]$}
\qbezier[20](8.8,0.11)(8.8,1)(8.8,2)

\put(8.7,0){\circle*{0.1}}
\put(8.2,-0.35){\tiny$[\![a,c|u]\!]$}
\qbezier[20](8.7,0)(8.7,1)(8.8,2)

\put(3,3.5){\circle*{0.1}}
\put(3.2,3.6){\footnotesize$b$}
\qbezier(3,3.5)(3,0)(0.5,0)
\qbezier(3,3.5)(3,0)(11.5,0)

\put(3.2,0){\circle*{0.1}}
\put(2.3,-0.35){\tiny$[\![a,c|b]\!]= [\![a,c;0]\!]_b$}
\qbezier[40](3,3.5)(3.15,1)(3.2,0)
   \end{picture}
  \end{center}
\caption{Synchronous exponential convergence of lines}
\label{f_exp-conv}
\end{figure}
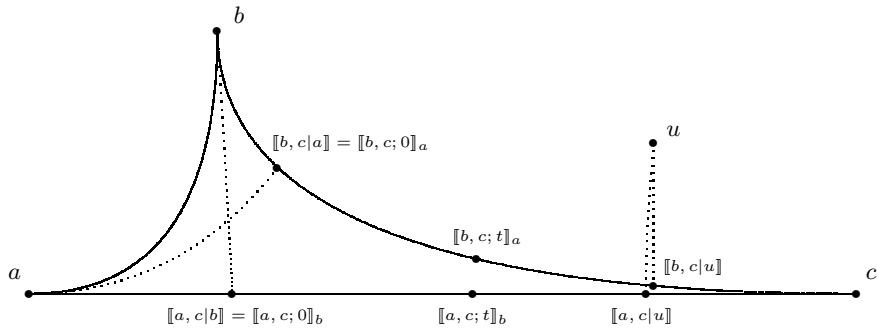
\proof
First note that if $a=b\in\p X$ or $a=c\in\p X$ or $b=c\in X$, then
$[\![b,c;t]\!]_{a}= [\![a,c;t]\!]_{b}\in\p X$ and the lemma obviously follows
(here we view $\dcross$ as a generalized pseudometric in $\sxb$).
So from now we will assume that those of $a,b,c$ that lie in~$\p X$ are pairwise distinct.

Let $T\in [0,\infty)$ and $\la\in[0,1)$ be the constants from Lemma~\ref{abba}.
Increase $\la$ if needed so that $\la\in(1/e,1)$. First we will show that
\begin{equation}\label{d2la}
  \dcross\big([\![b,c;t]\!]_{a}, [\![a,c;t]\!]_{b}\big)\le 2\la^t\quad
  \mbox{for all}\ t\in[T,\infty),
\end{equation}
i.e.\ that the lemma holds under the additional assumption $t\ge T$.
By the definition of $\dcross$, (\ref{d2la}) is equivalent to
\begin{equation*}
\left|\be_u^{\scriptscriptstyle\times\!}\big([\![b,c;t]\!]_{a}, [\![a,c;t]\!]_{b}\big)\right|\le 2\la^t\quad
  \mbox{for all}\ u\in X,\ t\in[T,\infty).
\end{equation*}
By the definition of~$\becross$,
\begin{eqnarray*}
&& \be_u^{\scriptscriptstyle\times\!}\big([\![b,c;t]\!]_a,[\![a,c;t]\!]_b\big)=
\be_u^{\scriptscriptstyle\times\!}\big([\![b,c;t+\left<b,c|a\right>]\!],[\![a,c;t+\left<a,c|b\right>]\!]\big)\\
&& =\left<b|c\right>_u+ \left|t+\left<b,c|a\right>-\left<b,c|u\right>\right|
  -\left<a|c\right>_u- \left|t+\left<a,c|b\right>-\left<a,c|u\right>\right|\\
&& = -\left<u,c|a,b\right>+ \left|t-\left<u,a|b,c\right>\right|-
  \left|t-\left<u,b|a,c\right>\right|.
\end{eqnarray*}
This equality will be used in the computations below. Pick an arbitrary $u\in X$.

If $\max\{\left<u,a|b,c\right>,\left<u,b|a,c\right>\}\le t$, then
$$\big|\be_u^{\scriptscriptstyle\times\!}\big([\![b,c;t]\!]_a,[\![a,c;t]\!]_b\big)\big|=
  \big| -\left<u,c|a,b\right>+ (t-\left<u,a|b,c\right>)- (t-\left<u,b|a,c\right>)\big|=0\le 2\la^t.$$

If $\max\{\left<u,a|b,c\right>,\left<u,b|a,c\right>\}\ge t$
  (see Figure~\ref{f_exp-conv}), then\newline
 $\max\{\left<u,a|b,c\right>,\left<u,b|a,c\right>\}\ge T$
and by Proposition~\ref{abba},
\begin{eqnarray*}
&& \left|\be_u^{\scriptscriptstyle\times\!}\big([\![b,c;t]\!]_a,[\![a,c;t]\!]_b\big)\right|
 \le \left|\left<u,c|a,b\right>\right|+
  \big| \left|t-\left<u,a|b,c\right>\right|-
  \left|t-\left<u,b|a,c\right>\right|\big|\\
&& \le \left|\left<u,c|a,b\right>\right|+\left|\left<u,b|a,c\right>- \left<u,a|b,c\right>\right|=
  2\left|\left<u,c|a,b\right>\right|\\
&& \le 2\la^{\max\{\left<u,a|b,c\right>,\left<u,b|a,c\right>\}}\le 2\la^t.
\end{eqnarray*}
This proves~(\ref{d2la}).
By calculus, there exists $N'\in[0,\infty)$ depending only on $T$ and $\la$ such that
$$2(\la^t+T-t)\le N'\la^t \quad\mbox{for all}\quad t\le T.$$
Let $N:= \max\{2, N'\}$, then~(\ref{d2la})
implies that the lemma holds for all $t\ge T$.
It remains to prove the lemma under the assumption $t\le T$.
By Theorem~\ref{t_d^times-Xbar}(b), (\ref{d2la}) and
since $[\![b,c;\cdot]\!]_a$ and $[\![a,c;\cdot]\!]_b$ are non-expanding,
\begin{align*}
& \dcross\big([\![b,c;t]\!]_a,[\![a,c;t]\!]_b\big)\\
&\le \dcross\big([\![b,c;t]\!]_a,[\![b,c;T]\!]_a\big)+
   \dcross\big([\![b,c;T]\!]_a,[\![a,c;T]\!]_b\big)+
   \dcross\big([\![a,c;T]\!]_b,[\![a,c;t]\!]_b\big)\\
& \le \dcross\big([\![b,c;T]\!]_a,[\![a,c;T]\!]_b\big)+ 2(T-t)
\le 2\la^T + 2(T-t)\\
& \le 2(\la^t+T-t)\le N'\la^t\le N\la^t.\tag*{\qed}
\end{align*}

\section{Translation length}
\label{s_translation}
For a hyperbolic complex $X$, we define 
the {\sf translation length}\label{i_translation length}
$\hat{l}(g)$ of\newline
$g\in \Isom(X)$ via the metric~$\dhat$, that is 
\begin{equation}\label{d_tr-l}
\hat{l}(g):=\lim_{n\to \infty} \dhat(x,g^nx)/n
\end{equation}
for some point $x\in X$. By the triangle inequality,
$\hat{l}(g)$ is independent of the choice of $x$.
In particular $\hat{l}(g)=\hat{l}(g')$ if $g$ and $g'$
are conjugate in $\Isom(X)$, so $\hat{l}$ can be viewed as
the {\sf length spectrum} of $\Isom(X)$, i.e.\ a function on the conjugacy classes.
This concept is of interest from the geometric point of view, for example
Otal~\cite{Ot1} and Croke~\cite{Croke1990} showed that each negatively curved metric on a surface 
is determined up to an isometry by its length spectrum.

An isometry of~$X$ is called {\sf elliptic} if it has a bounded orbit in~$X$;
this implies that each orbit is bounded.
An isometry $g$ is called {\sf hyperbolic} if 
its translation length (say with respect to the word metric)
is positive; this implies that $g$ fixes exactly 2 points at infinity.
We denote $g_-$\label{i_g+g-} the repelling point and $g_+$ the attracting point.
All isometries of a hyperbolic complex $X$ are either elliptic or hyperbolic.

Using $\dhat$ instead of the word metric enables us to express
translation length in terms of the double difference as follows
(cf \cite{HP} for the $\mathrm{CAT}(-1)$ case).

\begin{ppp}
$\hat{l}(g)$ in~(\ref{d_tr-l}) is a well-defined real number. Moreover, 
for any hyperbolic $g\in \Isom(X)$ and any 
$x\in \bar{X}\setminus \{g_-,g_+\}$,
$$\hat{l}(g)=  \left<g_-, g_+| gx, x\right>
=\ln |\![g_-,g_+|gx,x]\!|.$$
In particular, $\left<g_-, g_+| gx, x\right>$ 
and $|\![g_-,g_+|gx,x]\!|$ are independent of the choice of~$x$.
\end{ppp}
\proof
If $g$ is elliptic then the orbit of any $x\in X$ is bounded, hence $\hat{l}(g)=0$.
If $g$ is hyperbolic, pick a geodesic $\gamma$ from~$g_-$ to~$g_+$ and $y\in\gamma$.
All geodesics $g^n\gamma$ are $\de$--close to~$\gamma$, so  for each~$n$ we choose
a point $y_n\in\gamma$ which is $\de$--close to $g^n y$.
By \cite[Proposition 10(b)]{MY} and the definition of $\dhat$, there is a constant
$C\in[0,\infty)$ such that for all $n\ge 0$ and all $v\in\gamma$ between $g^n y$ and $g_+$,
$|\dhat(y,v)-\dhat(y_n,v)-\dhat(y,y_n)|\le C$. Since $g^n y$ and $y_n$ are $\de$--close,
we can also assume $|\dhat(y,v)-\dhat(y,g^ny)-\dhat(g^ny,v)|\le C$.
Similarly, for all $u\in\gamma$ between $g_-$ and $y$, 
$|\dhat(u,g^ny)-\dhat(u,y)-\dhat(y,g^ny)|\le C$.
Combining the two inequalities,
\begin{eqnarray*}
&&|\left< u,v|g^ny,y\right>-\dhat(y,g^ny)|
  \le |\dhat(y,v)-\dhat(y,g^ny)-\dhat(g^ny,v)|/2\\
&&+ |\dhat(u,g^ny)-\dhat(u,y)-\dhat(y,g^ny)|/2\le C/2+C/2=C.
\end{eqnarray*}
By continuity of double difference as $u\to g_-$ and $v\to g_+$,
\begin{equation}\label{g-g+}
|\left< g_-,g_+|g^ny,y\right>-\dhat(y,g^ny)|\le C.
\end{equation}
By the invariance of double difference under~$g$,
\begin{align*}
&\left< g_-,g_+|gy,y\right>= \left< g_-,g_+|g^2y,gy\right>=\ldots=
  \left< g_-,g_+|g^ny,g^{n-1}y\right>,\quad\mbox{hence}\\
&\left< g_-,g_+|g^ny,y\right>=\left< g_-,g_+|gy,y\right>+\ldots+
\left< g_-,g_+|g^ny,g^{n-1}y\right>=n\left< g_-,g_+|gy,y\right>,
\end{align*}
and~(\ref{g-g+}) rewrites as
$|\left< g_-,g_+|gy,y\right>-\dhat(y,g^ny)/n|\le C/n$.
Therefore\newline 
$\lim_{n\to\infty}\dhat(y,g^ny)/n$ exists and equals
$\left< g_-,g_+|gy,y\right>$, so $\hat{l}(g)$ is well-defined and 
$\hat{l}(g)=\left< g_-,g_+|gy,y\right>$.

For any $x\in \bar{X}\setminus \{g_-,g_+\}$,
by $g$--invariance, $\left< g_-,g_+|gx,gy\right>=\left< g_-,g_+|x,y\right>$, hence
\begin{align*}
& \left< g_-,g_+|gx,x\right>=\left< g_-,g_+|gx,gy\right>+
\left< g_-,g_+|gy,y\right>+ \left< g_-,g_+|y,x\right>\\
&= \left< g_-,g_+|gy,y\right>=\hat{l}(g).\tag*{\qed}
\end{align*}

\section{The geodesic flow of a hyperbolic complex}
\label{s_the-geod-flow}
Let $X$ be a hyperbolic complex. As an example, the reader might think
of a Cayley graph of a hyperbolic group, of the group itself.
Recall from~\ref{ss_dhat} that 
$X$ admits a nice canonical metric~$\dhat$, and we use it to define 
the (generalized) metric $d_*$ on the symmetric join $\sxb$. 
A part of the above symmetric join construction
is the {\sf flow space} of~$X$,
$$\ff(X):=*(\p X):= \sj(\p X)\setminus \p X\se \sxb.$$
As a set, $\ff(X)$ is a disjoint union
of open lines connecting disjoint ordered pairs of points at infinity.
We will use the same notation $d_*$ for the restriction of~$d_*$ to~$\ff(X)$.
$(\ff(X),d_*)$ plays the role of the total space of the unit tangent bundle on $X$
(though no bundle map is there), and it is canonically defined for any hyperbolic complex~$X$.

\begin{ppp}
\begin{itemize}
\item [\rm(a)] For $a,b\in\p X$, $[\![a,b;\cdot]\!]'=[\![a,b;\cdot]\!]$.
\item [\rm(b)] The restrictions of $\dcross$ and $d_*$ to each line in~$\ff(X)$ coincide
with the original metric on the line.
\item [\rm(c)] For each $x\in \fs(X)$, the orbit map $\rr\to(\fs(X),d_*)$, $r\mapsto r^+x$,
is an isometry onto the~$\rr$--orbit containing~$x$.
\end{itemize}
\end{ppp}
\proof
\noindent {\bf (a)}\qua Since $a,b\in\p X$,  $[\![a,b]\!]$ is a copy of
$[-\infty,\infty]$ and\newline
$[\![a,b;\cdot]\!]\co[-\infty,\infty]\to [\![a,b]\!]$ is the identity map.
By definitions in~\ref{ss_parametrizations},
$$[\![a,b;t]\!]'=
  \int_{-\infty}^\infty [\![a,b;r+t]\!]\frac{e^{-|r|}}{2}\, dr=
  \int_{-\infty}^\infty (r+t)\frac{e^{-|r|}}{2}\, dr= t= [\![a,b;t]\!].$$
\noindent {\bf (b)}\qua
By (a) and Theorem~\ref{t_d^times-Xbar}(b),
\begin{eqnarray*}
&& d_*\big([\![a,b;s]\!]',[\![a,b;t]\!]'\big)=
    \int_{-\infty}^\infty
    \dcross\big([\![a,b;r+s]\!]', [\![a,b;r+t]\!]'\big) \frac{e^{-|r|}}{2}\, dr\\
&& =\int_{-\infty}^\infty
     \dcross\big([\![a,b;r+s]\!], [\![a,b;r+t]\!]\big) \frac{e^{-|r|}}{2}\, dr
    = \int_{-\infty}^\infty
     |t-s| \frac{e^{-|r|}}{2}\, dr=|t-s|\\
&& =\dcross \big([\![a,b;s]\!],[\![a,b;t]\!]\big)=
   \dcross \big([\![a,b;s]\!]',[\![a,b;t]\!]'\big).
\end{eqnarray*}

\noindent {\bf (c)}\qua Let $x=[\![a,b;t]\!]'$, $a,b\in \p X$, then by (a), (b)
and Theorem~\ref{t_d^times-Xbar}(b),
\begin{align*}
& d_*(r_1^+x,r_2^+x)= d_*\big([\![a,b;r_1+t]\!]', [\![a,b;r_2+t]\!]'\big)\\
& \qquad= \dcross\big([\![a,b;r_1+t]\!], [\![a,b;r_2+t]\!]\big)=|r_2-r_1|.
\tag*{\qed}
\end{align*}

Let  $\p^2 X:= \{(a,b)\in (\p X)^2\ |\ a\not= b\}$.\label{i_d2x}
Recall from~\ref{ss_equivalence} that $^{\times +}$isometry is quasiisometry,
and $^+$isometry is  ``quasiisometry with multiplicative constant~1''.
We summarize the properties of the flow space (cf \cite[8.3.C]{Gr2}).

\begin{ttt}[Geodesic flow of~$X$]
\label{geodesic-flow}
Let $X$ be a hyperbolic complex, $d$ the word metric and $\dhat=\dhat_X$
the canonical metric as in~\ref{ss_dhat}. Then there is a metric
space~$(\fx, d_*)$ canonically associated to $X$ with the following properties.
\begin{itemize}
\item [\rm(a)] $(\ff(X),d_*)$ is  homeomorphic to $\p^2 X\times\rr$.
\item [\rm(b)] $(\ff(X),d_*)$ is proper.
\item [\rm(c)] If \,$\Isom(X)$ has a cobounded orbit in~$X$, for example if $X$ admits a cocompact
  isometric action, then $(\ff(X),d_*)$ is $^+$isometric to $(X,\dhat)$
  and $^{\times +}$isometric to $(X,d)$. 
\item [\rm(d)] There is a canonical isometric action of $\Isom(X)$ on $(\ff(X),d_*)$.
\item [\rm(e)] There is a canonical free $\rr$--action $(r,x)\mapsto r^+x$ 
  on $(\ff(X),d_*)$
  by bi-Lipschitz homeomorphisms which commutes with the $\Isom(X)$--action.
  For each $x\in\ff(X)$
  the orbit map $\rr\to (\fx,d_*)$, $r\mapsto r^+x$, is an isometry onto the $\rr$--orbit
  containing~$x$. In particular, $\rr$ acts by isometries on each $\rr$--orbit in~$\fx$,
  in the standard way.
\item [\rm(f)] There is a canonical free $\zz_2$--action on $(\ff(X),d_*)$,
  $x\mapsto x^\star$, by isometries which commutes with the $\Isom(X)$--action
  and anticommutes with the $\rr$--action.
  It moves every point a uniformly bounded distance and fixes $\p X$ pointwise.
\item [\rm(g)] There exists a horofunction $\becross_u\co\ff(X)^2\to \rr$ which is continuous
 in three variables $(u,x,y)\in \bar{X}\times\ff(X)^2$ and satisfies
\begin{itemize}
\item [\rm(i)] {\rm ($\zz_2$--invariance)}
   $\becross_u(x,y)=\becross_u(x^\star,y)=\becross_u(x,y^\star)$.
\item [\rm(ii)] {\rm (cocycle condition)} 
   $\becross_u(x,y)+ \becross_u(y,z)=\becross_u(x,z)$. 
\item [\rm(iii)] {\rm ($\Isom(X)$--invariance)}
   $\becross_{gu}(gx,gy)=\becross_u(x,y)$ for $g\in\Isom(X)$.
\item [\rm(iv)] {\rm (isometry on $\rr$--orbits)}
  For each $\rr$--orbit $]\!]a,b[\![$ in $\ff(X)$ and\newline
  $x,y\in ]\!]a,b[\![$, \ 
  $|\becross_a(x,y)|=|\becross_b(x,y)|=d_*(x,y)$.
\end{itemize}
\item [\rm(h)] {\rm (exponential convergence in $d_*$)} There exist $M\in [0,\infty)$ and
$\la\in [0,1)$ depending only on~$X$ with the following property.
For all $a,b,c\in \p X$, take the isometric parametrizations
$]\!]a,c;\cdot[\![_b\co\rr\to ]\!]a,c[\![$ and $]\!]b,c;\cdot[\![_a\co\rr\to ]\!]b,c[\![$ 
of the $\rr$--orbits  $]\!]a,c[\![$ and 
$]\!]b,c[\![$ as described in~\ref{ss_param-of-sxb}.
Then $[\![b,c;t]\!]_a$ and $[\![a,c;t]\!]_b$ lie on the same
horosphere at $c$ and
$$d_*\big([\![a,c;t]\!]_b, [\![b,c;t]\!]_a\big)\le M\la^t.$$
Similarly for any pair of lines among $]\!]a,c[\![$, $]\!]b,c[\![$, $]\!]c,a[\![$, $]\!]c,b[\![$.
\item [\rm(i)] Let $\hat{l}$ be the translation length in $X$ defined in section~\ref{s_translation}.
  Then for any hyperbolic $g\in \Isom(X)$ and any $y\in\ff(X)$,
  $$\hat{l}(g)=\lim_{n\to\infty} d_*(y,g^ny)/n=
     \inf\, \{d_*(y,gy)\,|\,y\in \ff(X)\}.$$
  Moreover, if $z\in ]\!]g_-,g_+[\![$ then $d_*(z,gz)=\hat{l}(g)$.
\end{itemize}
\end{ttt}

\noindent {\bf Remark}\qua
Note that the symmetric join $\sxb$ from section~\ref{s_sj-of-barX} gives even more structure:
both $X$ and $\ff(X)$ isometrically embed into $\sxb$, so $\sxb$ can be thought of
as a ``filling'' between $X$ and $\ff(X)$. This is a metric analogue of the following
geometric situation: if $Y$ is a smooth manifold,
then the space $B_1 Y$ of tangent vectors of length at most~1 contains both
the unit sphere bundle $S_1Y$ and $Y$, as the 0-section, 
so $B_1 Y$ is a ``filling'' between~$S_1 Y$ and~$Y$.

\begin{proof}[Proof of Theorem~\ref{geodesic-flow}]
\noindent {\bf (a)}\qua
The homeomorphism
$$]\!]\cdot,\cdot\,;\cdot[\!['\co 
  (\bar{X}^2\setminus\bar{\Delta})\times\rr\to (\starxbxo,d_*)$$
from  Proposition~\ref{p_Xbar-homeomorphic} maps
$\p^2 X\times\rr$ onto $\ff(X)$.

\noindent {\bf (b)}\qua
The same proof as in Proposition~\ref{p_proper}
with $(\ff(X),d_*)$ instead of $(\hsxb,d_*)$.

\noindent {\bf (c)}\qua The proof of $^+$isometry between 
$(\ff(X),d_*)$ and$(X,\dhat)$ is the same as in Proposition~\ref{+isom}.
The map $\psi\co\ff(X)\to X$ is not surjective, but the assumptions imply that it has
cobounded image in~$X$; this is sufficient to run the argument.
The $^{\times +}$isometry of $(\ff(X),d_*)$ and $(X,d)$
follows since $d$ and $\dhat$ are $^{\times +}$equivalent by
Theorem~\ref{t_dhat}(a).
(d)--(g) were proved in earlier sections.

\noindent {\bf (h)}\qua
Let $\la\in(1/e,1)$ and $N\in[0,\infty)$ be the constants from Theorem~\ref{exp-conv}.
Using the identity $[\![\cdot,\cdot\,;\cdot]\!]'= [\![\cdot,\cdot\,;\cdot]\!]$ in $\ff(X)$
and the definition of $d_*$,
\begin{eqnarray*}
&& d_*\big([\![a,c;t]\!]_b,[\![b,c;t]\!]_a\big)=
   d_*\big([\![a,c;t+\left<a,c|b\right>]\!],[\![b,c;t+\left<b,c|a\right>]\!]\big)\\
&& =\int_{-\infty}^\infty 
  \dcross\big([\![a,c;t+\left<a,c|b\right>]\!],[\![b,c;t+\left<b,c|a\right>]\!]\big)
  \frac{e^{-|r|}}{2}\, dr\\
&& = \int_{-\infty}^\infty \dcross \big([\![a,c;r+t]\!]_b,[\![b,c;r+t]\!]_a\big)
  \frac{e^{-|r|}}{2}\, dr\\
&& \le \int_{-\infty}^\infty N\la^{r+t}\,\frac{e^{-|r|}}{2}\, dr=
  \frac{N \la^t}{1-(\ln \la)^2},
\end{eqnarray*}
so we denote $M:=N/(1-(\ln \la)^2)$.

\noindent {\bf (i)}\qua $g\in\Isom(X)$ is an isometry of $(\hsxb,d_*)$.
For any $y\in\hsxb$, the limit $\lim_{n\to\infty} d_*(y,g^ny)/n$ exists,
and it is independent of $y$ (see \cite[II.6.6(1)]{BH}). Since $d_*$ coincides with
$\dhat$ on $X$, for any $x\in X$ we have
$$\lim_{n\to\infty} d_*(y,g^ny)/n= \lim_{n\to\infty} d_*(x,g^nx)/n=
    \lim_{n\to\infty} \dhat(x,g^nx)/n=\hat{l}(g).$$
This proves the first equality. For all $y\in\hsxb$,
$$d_*(y,g^ny)/n\le \big(d_*(y,gy)+\ldots+d_*(g^{n-1}y,g^ny)\big)/n=d_*(y,gy),$$
hence $\hat{l}(g)=\lim_{n\to\infty} d_*(y,g^ny)/n\le d_*(y,gy)$, so
\begin{equation}\label{lhatgle}
\hat{l}(g)\le \inf\, \{d_*(y,gy)\, |\, y\in \ff(X)\}.
\end{equation}
Let $g_-,g_+\in\p X$ be the fixed points of~$g$ and $z\in ]\!]g_-,g_+[\![$, then by 
the last two properties in~(g),
\begin{eqnarray*}
&& d_*(z,gz)=|\becross_{g_+}(z,gz)|=
  |\becross_{g_+}(z,gz)+\ldots+\becross_{g_+}(g^{n-1}z,g^nz)|/n\\
&& =|\becross_{g_+}(z,g^nz)|/n= d_*(z,g^nz)/n
\quad \mbox{for any}\ n,
\end{eqnarray*}
therefore $d_*(z,gz)= \hat{l}(g)$ and~(\ref{lhatgle}) becomes an equality.
\end{proof}

\section{Asymmetric join, the Borel conjecture and general remarks}
\label{s_asym-join}
\subsection{The definition of asymmetric join}\label{ss_asym-join}
If $(Y_1,d_1)$ and $(Y_2,d_2)$ are metric spaces, let $X:=Y_1\sqcup Y_2$ and pick a metric
$d$ on $X$ which induces the original topologies on~$Y_1$ and~$Y_2$.
Then $d$ canonically extends to the metric $d_*=\sj d$ on $\sx$ as in~\ref{ss_metric-d*},
and we define the {\sf asymmetric join}\label{i_asym-join}
of $Y_1$ and $Y_2$ to be the subspace $Y_1\sj\, Y_2\se \sx$ which is
the union of all lines in $\sx$ going from points in $Y_1$ to points in $Y_2$,
with the restricted metric~$d_*$.

For arbitrary~$Y_1$ and~$Y_2$, there is a choice involved in defining
$d$ on the union $Y_1\sqcup Y_2$, but if $Y_1$ and $Y_2$ are 
pointed metric spaces with isometric actions by the same group~$G$,
then $d$ and $d_*$ can be defined
canonically, as described below in~\ref{ss_metric-on-union}.
First let us mention some important examples
when this situation arises:
\begin{itemize}
\item [\rm(1)] Under the hypotheses of the Borel conjecture (in the PL setting), 
  $M_1$ and $M_2$ are closed triangulated manifolds with the same fundamental group~$\Gamma$.
  Consider the universal cover $Y_i:=\tilde M_i$
  with any $\Gamma$--invariant metrics $d_i$ induced from 0-skeleton
  (cf~\ref{ss_metric-complexes}(3)), and let $G:=\Gamma$.
  If $\Gamma$ is hyperbolic, put the canonical metrics $\dhat_{Y_i}$
  on $Y_i$ (cf~\ref{ss_examples-of-hyp-c}(3)).
\item [\rm(2)] If $M_1$ and $M_2$ are smooth manifolds with the same fundamental group
  $\Gamma$, let $d_i$ be the $\Gamma$--invariant intrinsic metric on the universal
  cover $Y_i:=\tilde M_i$, and $G:=\Gamma$.
\item [\rm(3)] For arbitrary metric complexes $(Y_1,d_1)$ and $(Y_2,d_2)$ one could just let $G$ be
  the trivial group.
\end{itemize}

\subsection{A metric on $Y_1\sqcup Y_2$}
\label{ss_metric-on-union}
Suppose that $(Y_1,d_1)$ and $(Y_2,d_2)$ are metric spaces
with isometric actions by the same group~$G$.
Let $X:=Y_1\sqcup Y_2$.
Pick basepoints $y_1\in Y_1$ and $y_2\in Y_2$. A pair $(y,z)$ of points
in $X$ is called {\sf admissible} if
\begin{itemize}
\item [\rm(a)] both $y$ and $z$ belong to the same
$Y_i$, or 
\item [\rm(b)] $y=gy_1$ and $z=gy_2$ for some $g\in G$, or
\item [\rm(c)] vise versa, $y=gy_2$ and $z=gy_1$ for some $g\in G$.
\end{itemize}
The {\sf length} of an admissible pair $(y,z)$, $l(y,z)$, is defined to be $d_i(y,z)$ in case $(a)$
and 1 in cases $(b)$ and $(c)$.
A finite sequence $x_1, \ldots, x_n$ of points in $X$ is called {\sf admissible}
if each consecutive pair $(x_j,x_{j+1})$ is admissible. The {\sf length}
of an admissible sequence, $l(x_1,\ldots,x_n)$, is $\sum_{j=1}^{n-1} l(x_j,x_{j+1})$.
We define a metric $d$ on $X$ by
$$d(a,b):= \inf \ l(x_1,\ldots,x_n)$$
over all admissible sequences $x_1,\ldots,x_n$ in $X$ with $x_1=a$ and $x_n=b$.
If $a,b\in Y_i$ and $d_i(a,b)\le 1$, then $d(a,b)=d_i(a,b)$,
i.e.\ $d$ and $d_i$ locally coincide on $Y_i$. Therefore~$d$ induces
the original topologies on $Y_1$ and $Y_2$.
If $G$ has a cobounded orbit in $Y_i$, for example if the $G$--action
on $Y_i$ is cocompact, then the embedding
$(Y_i,d_i)\hookrightarrow (X,d)$ is a quasiisometry.

The above definition of~$d$ on~$X=Y_1\sqcup Y_2$ makes the construction of the asymmetric join
$Y_1\sj Y_2$ and of the metric $d_*=\sj d$ canonical. In the case (1) above,
if~$M_1$ and~$M_2$ are closed manifolds and $\pi_1(M_1)=\pi_1(M_2)$
is hyperbolic, then the join $\bar{Y}_1\sj \bar{Y}_2$ of the compactifications
is also well-defined, equipped with the generalized metric $d_*$.
This allows for the use of both local and global structures of the manifolds.

The asymmetric join is another example of
a metric join\label{i_metric_join2}.

\subsection{Various model spaces for hyperbolic groups}
Given a hyperbolic group $\Gamma$, one can take
$X$ to be either 
\begin{itemize}
\item the group $\Gamma$ itself, 
\item or a Cayley graph of $\Gamma$,
\item or any other simplicial complex on which $\Gamma$ acts 
  (say, cocompactly).
\end{itemize}
We can put the metric $\dhat$ on $X$ as in~\ref{ss_dhat},
then the constructions of this paper provide many model
spaces, each equipped with the (generalized) 
metric $d_*$ and a $\Gamma$--action, for example
$\sx$, $\sxb$, $\ff(X)$, $X\sj X$, $X\sj \bar{X}$, $\bar{X}\sj \bar{X}$.
The symmetric and asymmetric join functors can be iterated, so 
for example, $\sj((\sj(X\sj \bar{X}))\sj( \starx))$ is a legitimate model space with
an isometric $\Gamma$--action.

\subsection{More general metric spaces}
We started with an arbitrary metric space $(X,d)$ and defined $\dcross$ and $d_*$
on $\sx$. When $X$ is a hyperbolic complex, the metric $\dhat$ on $X$
was used to define $\dcross$ and $d_*$ on $\sxb$ (with infinite values allowed).
The metric $\dhat$ was used because of its strong 
properties at infinity. The definition of $\dcross$ and $d_*$ can be
carried out if one starts with a more general hyperbolic metric space
$(X,\dhat)$, where $\dhat$ satisfies the properties 
of~\ref{ss_prop-of-dhat} through~\ref{ss_more-prop-dd}.
For example this would work for $CAT(-1)$ spaces.

\subsection{The sweep-out and the Borel conjecture}\label{ss_sweep-out}
Each line $[\![a,b]\!]$ in $Y_1\sj Y_2$ with $a\in Y_1$ and $b\in Y_2$, 
is given the canonical parametrization $[\![a,b;\cdot]\!]'\co\rb\to [\![a,b]\!]$.
Each $t\in\rb$ gives the point $[\![a,b;t]\!]'$ on each line $[\![a,b]\!]$
in $Y_1\sj Y_2$; the union $S_t$ of these points
for a fixed~$t$ will be called a {\sf slice}, and the set of all slices is the {\sf sweep-out}
from~$Y_1$ to~$Y_2$. For each $t\in\rr$, $S_t$ is homeomorphic
to $Y_1\times Y_2$ (this follows from Proposition~\ref{p_Xhomeomorphic}),
but $S_t$ converges to $Y_1$ 
as $t\to -\infty$ and to $Y_2$ as $t\to \infty$, in a metric sense
that can be made precise (Gromov--Hausdorff convergence on bounded subsets).

\setlength{\unitlength}{1cm}
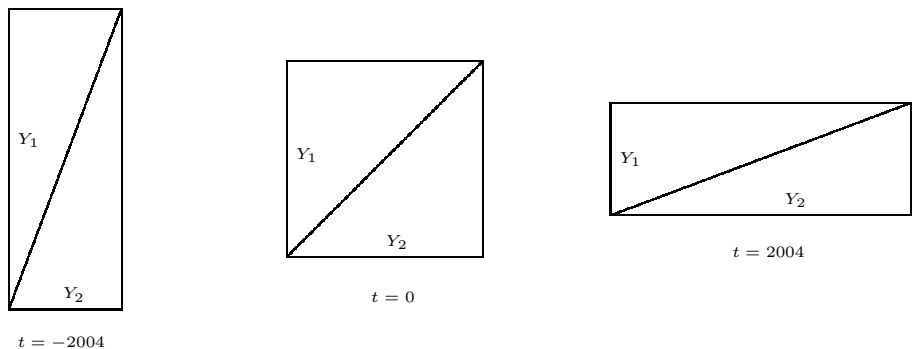
\begin{figure}[ht!]
  \begin{center}
   \begin{picture}(12,4)
\qbezier(0,0)(0.75,0)(1.5,0)  \qbezier(0,0)(0,2)(0,4)
\qbezier(1.5,4)(1.5,2)(1.5,0)  \qbezier(1.5,4)(0.75,4)(0,4)
\qbezier(0,0)(0.75,2)(1.5,4)
\put(0.7,0.15){\tiny$Y_2$}   \put(0.1,2.2){\tiny$Y_1$}
\put(0.1,-0.5){\tiny$t=-2004$}

\qbezier(3.7,0.7)(3.7,2)(3.7,3.3)  \qbezier(6.3,0.7)(6.3,2)(6.3,3.3)
\qbezier(3.7,0.7)(5,0.7)(6.3,0.7)  \qbezier(3.7,3.3)(5,3.3)(6.3,3.3)
\qbezier(3.7,0.7)(5,2)(6.3,3.3)
\put(5,0.85){\tiny$Y_2$}  \put(3.8,2){\tiny$Y_1$}
\put(4.8,0.1){\tiny$t=0$}

\qbezier(8,1.25)(8,2)(8,2.75)     \qbezier(12,1.25)(12,2)(12,2.75)
\qbezier(8,1.25)(10,1.25)(12,1.25)  \qbezier(8,2.75)(10,2.75)(12,2.75)
\qbezier(8,1.25)(10,2)(12,2.75)
\put(10.3,1.4){\tiny$Y_2$}  \put(8.1,1.95){\tiny$Y_1$}
\put(9.6,0.7){\tiny$t=2004$}

   \end{picture}
  \end{center}\vspace{1mm}
\caption{The sweep-out $S_t$ from $Y_1$ to $Y_2$, at different times}
\label{f_sweep-out}
\end{figure}

$X$ has the natural $\Gamma$--action induced from the $\Gamma$--actions on $Y_i$,
this provides an isometric $\Gamma$--action on $Y_1\sj Y_2$.
This action preserves slices and is the diagonal action on each slice $S_t=Y_1\times Y_2$.

When $M_1$ and $M_2$ are closed aspherical manifolds with the same 
fundamental group, the Borel conjecture asserts that
they are homeomorphic. If there is indeed a homeomorphism, then each
slice in $Y_1\sj Y_2$ contains a homeomorphic and $\Gamma$--invariant 
copy of $Y_1=Y_2$. Thus the asymmetric join is a place to look for homeomorphisms,
and the metric $d_*$ on it should allow for analytic and geometric tools to be used.

\subsection{Group-theoretic rigidity and the Poincar\'e conjecture}
The Mostow rigidity theorem~\cite{Most} implies that the Borel conjecture holds
in the case of hyperbolic manifolds: two closed hyperbolic manifolds with the same
fundamental group are homeomorphic. This can be viewed as an example
of a topological theorem where geometric assumptions are necessary to run the proof.

There is an interesting conjecture which is a topological version of the Mostow rigidity
theorem: if two closed aspherical manifolds $M$ and $N$ have the same
fundamental group and $N$ is hyperbolic, then $M$ and $N$ are homeomorphic.
This  was proved by Farrell and Jones~\cite{FJ1989, FJ1989-PNASUSA, FJ1990}
in the dimensions other than~3 and~4. Again, geometric assumptions were important
for the proof. Farrell and Jones use the dynamics of the geodesic flow on
a hyperbolic manifold: the flow shrinks certain paths in the unit tangent bundle.
Gabai, Meyerhoff and Thurston~\cite{GMT2003} showed this conjecture
for 3-manifolds under the additional assumption that $M$ is irreducible;
it is this assumption that prevents the Poincar\'e conjecture to be deduced
from the result.

From our metric (i.e.\ non-Riemannian) point of view, the following {\em group-theoretic}
version of the conjecture is of interest: if two closed aspherical manifolds
have the same fundamental group which is Gromov hyperbolic, then they
are homeomorphic. This conjecture is intermediate: it follows from
the Borel conjecture and, if true in dimension 3, it implies the Poincar\'e conjecture.
Theorem~\ref{geodesic-flow} provides a construction of
a geodesic flow $\ff(X)$ with the properties needed: the flow (i.e.\ the $\rr$--action)
indeed shrinks distances exponentially. What is missing here is the bundle structure,
for the obvious reason: we did not have a manifold to start with,
$X$ could be any metric space.
It is also worth mentioning that in our construction $\rr$ acts by bi-Lipschitz
homeomorphisms, so the topology and geometry of the space are preserved.

\newpage

\subsection{Index of symbols and terminology}\label{index}

\small
Symbols and terminology are indexed by section/subsection.  Readers interested
in the \LaTeX\ codes which produce the symbols should consult 
{\tt symbols.tex} or {\tt symbols.pdf} in the directory:

\url{http://www.maths.warwick.ac.uk/gt/ftp/aux/2005-13/}

\vspace{3mm}

\begin{minipage}{.3\hsize}
$X\!\Join\! X$,  \ref{ss_sj-as-top}\newline
$\diamondx$,  \ref{ss_sj-as-top}\newline
$\diamondxxo$, \ref{ss_parametrizations}\newline
$\sx$,  \ref{ss_sj-as-top}\newline
$\sxxo$,  \ref{ss_parametrizations}\newline
$Y_1\sj Y_2$,  \ref{ss_asym-join}\newline
$\sxb$,  \ref{s_sj-of-barX}\newline
$\starx$,  \ref{ss_sj-as-top}\newline
$\starxxo$, \ref{ss_two-models-starx-starxxo}\newline
$\starxb$,  \ref{ss_models-for-Xbar}\newline
$\raise0.12ex\hbox{$\scriptscriptstyle\smallsmile$}
   \mspace{-11.26mu}\raise0.2ex\hbox{$*$}\bar{X}$,
           \ref{ss_models-for-Xbar}\newline
${\bar{X}}^{\diamond}$,  \ref{ss_ext-of-dd}\newline
$ \bar{X}^\triangleright$, \ref{ss_ext-of-dd}\newline 
\noindent $\simplus$,  \ref{ss_equivalence}\newline
$^+$equivalence,  \ref{ss_equivalence}\newline
$^\times$equivalence,   \ref{ss_equivalence}\newline
$^{\times +}$equivalence,   \ref{ss_equivalence}\newline
$^+$geodesic,  \ref{ss_prop-of-dhat}\newline
$^+$isometry,  \ref{ss_equivalence}\newline
$^{\times+}$isometry, \ref{ss_equivalence}\newline
$^+$map,  \ref{ss_equivalence}\newline
$^\times$map,  \ref{ss_equivalence}\newline
$^{\times +}$map, \ref{ss_equivalence}\newline 
\noindent $\left<\cdot,\cdot|\cdot,\cdot\right>$, 
    \ref{ss_dd-gp}, \ref{ss_dhat-everywhere}\newline
$(\cdot,\cdot|\cdot,\cdot)$, \ref{ss_dhat-everywhere}\newline
$|\![\cdot,\cdot|\cdot,\cdot]\!|$,  \ref{s_cross-ratio}\newline
$\left<a|b\right>_c$, 
     \ref{ss_dd-gp}, \ref{ss_dhat-everywhere}\newline
$(a|b)_c$,     \ref{ss_dhat-everywhere}\newline
$\left<a,a'|b\right>$,  \ref{ss_projX}\newline
$[\![a,a'|b]\!]$,  \ref{ss_projX}\newline
$[S_1,S_2]$,  \ref{ss_neighborhoods-in-barG}\newline
$[a,a'|x_0]$,  \ref{ss_map-psi}\newline
$[\![a,b]\!]=[\![a,b]\!]_{x_0}$\newline
  $\phantom{999}$ in $\sx$, \ref{ss_parametrizations}\newline
  $\phantom{999}$ in $\sxb$, \ref{ss_param-of-sxb}\newline
$[\![a,b;t]\!]$ and $[\![a,b;t]\!]'$\newline
  $\phantom{999}$ in $\sx$,  \ref{ss_parametrizations}\newline
  $\phantom{999}$ in $\sxb$, \ref{ss_param-of-sxb}
\vfill
\end{minipage}
\begin{minipage}{.33\hsize}
$\becross$, 
    \ref{ss_the-cocycle}, \ref{ss_be-ext}, \ref{s_horofunctions-in-barX}\newline
$\theta$,  \ref{ss_aux-fun}\newline
$\theta'$,  \ref{ss_aux-fun}\newline
$\psi=\psi_X$,  \ref{ss_map-psi}\newline 
\noindent $d=d_X$, \ref{ss_word-metric}\newline
$\dhat=\dhat_X$,  \ref{ss_dhat}\newline
$\dcross$ in $\sx$,  \ref{ss_d^times}\newline
$\dcross$ in $\sxb$,  \ref{ss_ext-dcross-dstar}\newline
$d_*=\sj d$  in $X$,  \ref{ss_metric-d*}\newline
$d_*=\sj \dhat$  in $\bar{X}$, \ref{ss_ext-dcross-dstar}\newline
$\p X$--triple,  \ref{ss_ext-of-dd}\newline
  $\phantom{999}$ trivial, \ref{ss_ext-of-dd}\newline
$\p^2 X$, \ref{s_the-geod-flow}\newline
$\ell(u,x)$,  \ref{ss_the-cocycle}\newline
$\ff(X)$,  \ref{s_the-geod-flow}\newline
$\Geod(a,b)$,  \ref{ss_hyp-compl}\newline
$g_+$, $g_-$,  \ref{s_translation}\newline
$\hat{l}(g)$,  \ref{s_translation}\newline
$\mathrm{np}[\al|y]$,  $\mathrm{np}[a,a'|y]$, 
     \ref{ss_geod-and-near}\newline
$r^+$,  \ref{ss_r-action}\newline
$\rb$,  \ref{ss_dealing-with-rb} \newline
$\T_*$,  \ref{ss_the-top-on-sxb}\newline
$U^+(a,t)$, $U^-(a,t)$,  \ref{ss_neighborhoods-in-barG}\newline
$V^+(a,t)$, $V^-(a,t)$, \ref{ss_neighborhoods-in-barG}\newline
$[\![\cdot,\cdot\,;\cdot]\!]'$,  \ref{ss_two-models-sx-sxxo}\newline
 $\phantom{999}$\newline
asymmetric join,  \ref{s_asym-join}\newline
actions on $\sx$\newline
  $\phantom{999}$ $\Isom(X)$--,  \ref{ss_isom-action}\newline
  $\phantom{999}$ $\rr$--,  \ref{ss_r-action}\newline
  $\phantom{999}$ $\zz_2$--, \ref{ss_z2-action}\newline
actions on $\sxb$,  \ref{ss_actions-on-sxb}\newline
appropriate number,  \ref{ss_parametrizations}\newline
coordinate\newline
  $\phantom{999}$ $x_0$--coordinate, 
    \ref{ss_parametrizations}\newline
  $\phantom{999}$ of the projection,  \ref{ss_projX}\newline
convexity of $\T_*$,  \ref{ss_the-top-on-sxb}\newline
  $\phantom{999}$
\vfill
\end{minipage}
\begin{minipage}{.33\hsize}
cross-ratio,  \ref{s_cross-ratio}\newline
distance-minimizing\newline 
  $\phantom{999}$ geodesic, 
      \ref{ss_geod-and-near}\newline
  $\phantom{999}$ pair, 
    \ref{ss_geod-and-near}\newline
double difference, 
    \ref{ss_dd-gp}, \ref{ss_dhat-everywhere}\newline
end-point coordinates, 
     \ref{ss_parametrizations}\newline
flow space,  \ref{s_the-geod-flow}\newline
generalized metric,  \ref{ss_gener-metr}\newline
geodesic, \ref{ss_hyp-compl}\newline
Gromov product, 
   \ref{ss_dd-gp}, \ref{ss_dhat-everywhere}\newline
horosphere,  \ref{s_horofunctions-in-barX}\newline
horofunction in $\sxb$, 
    \ref{s_horofunctions-in-barX}\newline
hyperbolic complex, 
    \ref{ss_hyp-compl}\newline
hyperbolic graph, 
    \ref{ss_hyp-compl}\newline
induced from 0-skeleton, 
   \ref{ss_metric-complexes}\newline
inscribed triple, 
   \ref{ss_hyp-compl}\newline
isometric embedding, 
   \ref{ss_gener-metr}\newline
metric complex,
   \ref{ss_metric-complexes}\newline
metric join, 
      \ref{ss_sj-as-top}, \ref{ss_metric-on-union}\newline
nearest-point projection,
   \ref{ss_geod-and-near}\newline
open symmetric join\newline
    $\phantom{999}$ of $X$,  \ref{ss_sj-as-top}\newline
    $\phantom{999}$ of $\bar{X}$, 
      \ref{ss_models-for-Xbar}\newline
projection,  \ref{ss_projX}\newline
proper metric space,  \ref{ss_proper}\newline
quasiisometry, 
   \ref{ss_equivalence}\newline
shift-invariance\newline 
    $\phantom{999}$ of $\theta$ and $\theta'$, 
   \ref{ss_aux-fun}\newline
side pair,  \ref{ss_ext-of-dd}\newline
side $\p X$--pair,  \ref{ss_ext-of-dd}\newline
sweep-out,  \ref{ss_sweep-out}\newline
symmetric join\newline
  $\phantom{999}$ of $X$,  \ref{ss_sj-as-top}\newline
  $\phantom{999}$ of $\bar{X}$,  \ref{s_sj-of-barX}\newline
translation length,  \ref{s_translation}\newline
trivial $\p X$--triple,  \ref{ss_ext-of-dd}\newline
word metric $d$, 
    \ref{ss_word-metric}, \ref{ss_hyp-compl}\newline
  $\phantom{999}$
\vfill
\end{minipage}

\end{document}